\newif\iffinal
\finalfalse	
\finaltrue	

\documentclass[letterpaper,12pt,reqno]{amsart}

\usepackage{mathtools}
\mathtoolsset{showonlyrefs=true}

\RequirePackage[utf8]{inputenc}
\usepackage[portrait,margin=1.265in]{geometry}
\iffinal\else\usepackage[notref,notcite]{showkeys}\fi
\usepackage{mathrsfs,soul}
\usepackage{hyperref}

\usepackage[foot]{amsaddr}
\usepackage{amssymb,amsmath,amsthm,amsfonts,amsbsy,latexsym,dsfont}
\usepackage{leftidx}
\usepackage[textsize=small]{todonotes}       
\usepackage{graphicx}
\usepackage[numeric,initials,nobysame]{amsrefs}
\usepackage{upref,setspace}
\usepackage{xcolor,colortbl}
\usepackage{enumerate}
\newenvironment{enumeratei}{\begin{enumerate}[\upshape (i)]}{\end{enumerate}}
\newenvironment{enumeratea}{\begin{enumerate}[\upshape (a)]}{\end{enumerate}}

\newenvironment{enumerateA}{\begin{enumerate}[\upshape (A)]}{\end{enumerate}}

\usepackage{tikz}
\usetikzlibrary{calc,intersections,through,backgrounds,shapes.geometric}
\usetikzlibrary{graphs}
 	\usepackage{pgfplots}

\usepackage{paralist}

\newenvironment{inparaenumi}{\begin{inparaenum}[\upshape \bfseries (i) ]}{\end{inparaenum}}
\newenvironment{inparaenumii}{\begin{inparaenum}[\upshape (i) ]}{\end{inparaenum}}

\newenvironment{inparaenumaa}{\begin{inparaenum}[\upshape(a)]}{\end{inparaenum}}

\numberwithin{equation}{section}
\numberwithin{figure}{section}
\numberwithin{table}{section}

\usepackage{caption}
\usepackage{subcaption}

\newtheorem{thm}{Theorem}[section]
\newtheorem{lem}[thm]{Lemma}

\newtheorem{prop}[thm]{Proposition}
\newtheorem{defn}[thm]{Definition}

\newtheorem{ass}[thm]{Assumption}

\newtheorem{conj}[thm]{Conjecture}
\newtheorem{lemma}[thm]{Lemma}
\newtheorem{constr}[thm]{Construction}

\theoremstyle{definition}
\newtheorem{rem}{Remark}
\newtheorem{obs}{Observation}

\numberwithin{obs}{section}

\renewcommand{\leq}{\le}
\renewcommand{\geq}{\ge}

\newcommand{\ind}{\mathds{1}}
\newcommand{\eps}{\varepsilon}

\newcommand{\set}[1]{\left\{#1\right\}}

\newcommand{\equald}{\stackrel{\mathrm{d}}{=}}
\newcommand{\probc}{\stackrel{\mathrm{P}}{\longrightarrow}}

\newcommand{\weakc}{\stackrel{\mathrm{d}}{\longrightarrow}}
\newcommand{\convas}{\stackrel{\mathrm{a.s.}}{\longrightarrow}}

\newcommand{\hght}{\mathrm{ht}}

\newcommand{\cb}{\mathrm{CB}}
\newcommand{\cbd}{\mathrm{CBD}}

\newcommand{\cmnd}{{\cG}_{n,\vd}}
\newcommand{\cmnthree}{{\cG}_{n,3}}

\newcommand{\cmnthreeexp}{{\cG}_{n,3}^{\exp}}
\newcommand{\cmnthreem}{{\cG}_{n,3}^{(m)}}
\newcommand{\cmndm}{{\cG}_{n,\vd}^{(m)}}
\newcommand{\cgnd}{\overline{\cG}_{n,\vd}}
\newcommand{\cgnthree}{\overline{\cG}_{n,3}}

\newcommand{\remove}{\mathrm{Rem}}
\newcommand{\shape}{\mathrm{Shape}}
\newcommand{\diam}{\mathrm{diam}}
\newcommand{\perc}{\mathrm{Perc}}
\newcommand{\lhms}{\overline{\cH}_{m,s}}

\newcommand{\ltms}{\overline {T}_m}

\newcommand{\lu}{\overline {u}}

\newcommand{\lv}{\overline {v}}

\def\qed{ \hfill $\blacksquare$}

\newcommand{\cA}{\mathcal{A}}\newcommand{\cC}{\mathcal{C}}
\newcommand{\cD}{\mathcal{D}}\newcommand{\cE}{\mathcal{E}}
\newcommand{\cG}{\mathcal{G}}\newcommand{\cH}{\mathcal{H}}
\newcommand{\cK}{\mathcal{K}}\newcommand{\cL}{\mathcal{L}}
\newcommand{\cM}{\mathcal{M}}
\newcommand{\cP}{\mathcal{P}}\newcommand{\cQ}{\mathcal{Q}}
\newcommand{\cS}{\mathcal{S}}\newcommand{\cT}{\mathcal{T}}\newcommand{\cU}{\mathcal{U}}
\newcommand{\cV}{\mathcal{V}}
\newcommand{\cZ}{\mathcal{Z}}

\newcommand{\vd}{\mathbf{d}}\newcommand{\ve}{\mathbf{e}}

\newcommand{\vt}{\mathbf{t}}

\newcommand{\mvthree}{\boldsymbol{3}}

\newcommand{\mvS}{\boldsymbol{S}}
\newcommand{\mvX}{\boldsymbol{X}}

\newcommand{\mvk}{\boldsymbol{k}}
\newcommand{\mvp}{\boldsymbol{p}}

\newcommand{\mvnu}{\boldsymbol{\nu}}

\newcommand{\mvzeta}{\boldsymbol{\zeta}}

\newcommand{\fB}{\mathfrak{B}}

\newcommand{\fG}{\mathfrak{G}}

\newcommand{\fM}{\mathfrak{M}}\newcommand{\fN}{\mathcal{H}}
\newcommand{\fQ}{\mathfrak{Q}}\newcommand{\fR}{\mathfrak{R}}
\newcommand{\fS}{\mathfrak{S}}\newcommand{\fT}{\mathfrak{T}}
\newcommand{\fX}{\mathfrak{X}}

\newcommand{\fm}{\mathfrak{m}}
\newcommand{\fp}{\mathfrak{p}}

\newcommand{\bE}{\mathbb{E}}

\newcommand{\bN}{\mathbb{N}}
\newcommand{\bP}{\mathbb{P}}\newcommand{\bR}{\mathbb{R}}
\newcommand{\bT}{\mathbb{T}}

\newcommand{\bZ}{\mathbb{Z}}

\newcommand{\dV}{\mathds{V}}

\newcommand{\rH}{\mathrm{H}}

\DeclareMathOperator{\pr}{\mathbb{P}}

\DeclareMathOperator{\var}{Var}
\DeclareMathOperator{\cov}{Cov}
\DeclareMathOperator{\dis}{dis}

\DeclareMathOperator{\GH}{GH}
\DeclareMathOperator{\GHP}{GHP}
 \DeclareMathOperator{\height}{ht}
  \DeclareMathOperator{\Ht}{Ht}

 \DeclareMathOperator{\ER}{ER}

\DeclareMathOperator{\er}{er}

\DeclareMathOperator{\conn}{Conn}
\DeclareMathOperator{\conne}{Conne}
\DeclareMathOperator{\core}{Core}

\DeclareMathOperator{\avail}{avail}
\DeclareMathOperator{\attach}{attach}

\DeclareMathOperator{\modi}{modi}

\DeclareMathOperator{\prev}{prev}
\DeclareMathOperator{\nxt}{next}

\newcommand{\sss}{\scriptscriptstyle}

\newcommand{\ch}[1]{\textcolor{black}{#1}}
\newcommand{\chl}[1]{\textcolor{black}{#1}}

\newcommand{\erdos}{Erd\H{o}s-R\'enyi }

\newcommand{\len}{\mathrm{len}}

\usepackage{fourier}

\definecolor{aqua}{rgb}{0.0, 1.0, 1.0}
\definecolor{webbrown}{rgb}{.6,0,0}
\definecolor{pinegreen}{rgb}{0.0, 0.47, 0.44}
\definecolor{ultramarineblue}{rgb}{0.25, 0.4, 0.96}
\definecolor{jrnl}{rgb}{0.0, 0.5, 0.0}
\definecolor{lincolngreen}{rgb}{0.11, 0.35, 0.02}
\definecolor{green(html/cssgreen)}{rgb}{0.0, 0.5, 0.0}
\definecolor{airforceblue}{rgb}{0.36, 0.54, 0.66}
\definecolor{azure}{rgb}{0.0, 0.5, 1.0}
\definecolor{bleudefrance}{rgb}{0.19, 0.55, 0.91}
\definecolor{cobalt}{rgb}{0.0, 0.28, 0.67}

\hypersetup{colorlinks=true, linktocpage=true, pdfstartpage=1, pdfstartview=FitV,
    breaklinks=true, pdfpagemode=UseNone, pageanchor=true, pdfpagemode=UseOutlines,
    plainpages=false, bookmarksnumbered, bookmarksopen=true, bookmarksopenlevel=1,
    hypertexnames=true, pdfhighlight=/O,
    urlcolor=webbrown, linkcolor=cobalt, citecolor=jrnl
    }

\begin{document}

\title[MST of random $3$-regular graphs]{Geometry of the minimal spanning tree of a random $3$-regular graph}

\date{}
\subjclass[2010]{Primary: 60C05, 05C80. }
\keywords{Minimal spanning tree, Gromov-Hausdorff distance, critical percolation, real tree, random regular graphs, graphs with prescribed degree sequence, configuration model}

\author[Addario-Berry]{Louigi Addario-Berry$^1$}
\address{\hskip-15pt $^1$Department of Mathematics and Statistics, McGill  University, Montreal, Canada}
\author[Sen]{Sanchayan Sen$^2$}
\address{\hskip-15pt $^2$Department of Mathematics, Indian Institute of Science, Bangalore, India}
\email{louigi.addario@mcgill.ca, sanchayan.sen1@gmail.com}

\maketitle
\begin{abstract}
The global structure of the minimal spanning tree (MST) is expected to be universal for a large class of underlying random discrete structures. 
\ch{However,} very little is known about the intrinsic geometry of MSTs of most standard models, and so far the scaling limit of the MST viewed as a metric measure space has only been identified in the case of the complete graph \cite{AddBroGolMie13}.

In this work, we show that the MST constructed by assigning i.i.d. continuous edge weights to either the random (simple) $3$-regular graph or the $3$-regular configuration model on $n$ vertices, endowed with the tree distance scaled by $n^{-1/3}$ and the uniform probability measure on the vertices, converges in distribution with respect to Gromov-Hausdorff-Prokhorov topology to a random compact metric measure space. Further, this limiting space has the same law as the scaling limit of the MST of the complete graph identified in \cite{AddBroGolMie13} up to a scaling factor of $6^{1/3}$. \ch{Our proof relies on a novel argument that proceeds via a comparison between a $3$-regular configuration model and the largest component in the critical \erdos random graph.}
The techniques of this paper can be used to establish the scaling limit of the MST in the setting of general random graphs with given degree sequences provided two additional technical conditions are verified.
\end{abstract}

\tableofcontents

\section{Introduction}
\label{sec:intro}
Consider a finite, connected, and weighted graph $(V,E,w)$, where $(V,E)$ is the underlying graph and $w:E\to [0,\infty)$ is the weight function.
A spanning tree of $(V,E)$ is a tree that is a subgraph of $(V,E)$ with vertex set $V$.
A minimal spanning tree (MST) $T$ of $(V,E,w)$ satisfies
\begin{align}\label{eqn:def-mst}
\sum_{e\in T}w(e)=\min\bigg\{\sum_{e\in T'}w(e):\ T'\text{ is a spanning tree of }(V,E)\bigg\}.
\end{align}
The two natural choices for the underlying weighted graph are
(i) a deterministic graph (e.g., the complete graph on $n$ vertices or the hypercube) or a random graph (e.g., \erdos random graph, random regular graph, or \ch{inhomogeneous random graphs \cite{BJR07}}) with i.i.d. continuous edge weights assigned to them, and
(ii) the complete graph on a finite set of random points in $\bR^d$ (e.g., $n$ i.i.d. points or a Poisson point process in the unit cube) where the edge weights are some function of the Euclidean length of the edges.
The MST in the latter case is sometimes called the Euclidean MST.

The MST is one of the most studied \ch{objects} in combinatorial optimization and geometric probability and has inspired a large body of work.
For an account of law of large numbers and related asymptotics in the Euclidean setting, see e.g., \cite{beardwood-halton-hammersley, aldous-steele, alexanderI, avram-bertsimas, steele}.
Central limit theorems (CLT) for the total weight of Euclidean MSTs were first proved by Kesten and Lee \cite{kesten-lee} and by Alexander \cite{alexander} in 1996.
This was a long-standing open question at the time of its solution.
Later certain other CLTs related to MSTs were proved in \cite{leeI, leeII}.
A question raised in \cite{kesten-lee} about the convergence rate in the CLT for the total weight of the Euclidean MST was answered in \cite{chatterjee-sen}.

Studies related to \ch{the MST} in several other directions were undertaken in \cite{leeIII, bhatt, penroseI, penroseII, penroseIII}.
An account of certain structural and connectivity properties of minimal spanning forests can be found in
\cite{alexander+molchanov, alexander-forest, lyons,newman2017critical} and the references therein.
For an account of the scaling limit of minimal spanning trees in subsets of $\bZ^2$ with respect to the topology introduced by Aizenman, Burchard, Newman, and Wilson, see, e.g., \cite{Aizenman-Burchard-Newman-Wilson, petegabor}.

The MST of $K_n$-the complete graph on $[n]:=\{1,\ldots,n\}$ has been studied extensively as well.
A celebrated theorem of Frieze \cite{frieze1985value} shows that under some assumptions on the weight distributions, the total weight of the MST of $K_n$ converges in expectation to $\zeta(3)$.
\ch{Various extensions of this result were proved in \cite{beveridge-frieze-mcdiarmid, frieze-mcdiarmid, frieze-ruszink-thoma, penrose1998random, aldous1990random}.}
The central limit theorem for the total weight of the MST of $K_n$ constructed using i.i.d. $\mathrm{Uniform}[0,1]$ edge weights was proved in \cite{janson1995minimal}.

The global geometric properties of the MST, e.g., the diameter and the typical distance, have also been of considerable \ch{interest, but until very recently, there were few rigorous mathematical results on this problem.}
Frieze and McDiarmid asked a question \cite[Research Problem 23]{frieze-mcdiarmid-algorithmic} about the `likely shape of a minimum spanning tree' and the order of the diameter of the MST.
In the statistical physics literature, paths in the MST correspond to optimal paths in the so-called strong disorder regime for complex networks.
\ch{Using empirical observations, it was predicted in \cite{braunstein2003optimal} (see also \cite{braunstein2007optimal}) that in the strong disorder regime, the length of optimal paths in complex networks should scale like $n^{1/3}$ if the degree distribution of the network has finite third moment, although a rigorous justification of this claim was missing in this work.}

An upper bound of the order $n^{1/3}$ on the diameter of the MST of $K_n$ was proved in \cite{addarioberry-broutin-reed}:
Let $M_{\infty}^{n,\er}$ be\footnote{
\ch{Here, the superscript `$\er$' is being used to refer to the \erdos random graph. 
The reason behind using this notation will become clear in Section~\ref{sec:alternate-description-M}.}} 
the MST of $K_n$ constructed using i.i.d. continuous edge weights, and denote by $\diam(M_{\infty}^{n,\er})$ the maximum tree distance between vertices of $M_{\infty}^{n,\er}$.
Then $\diam(M_{\infty}^{n,\er})=O_P(n^{1/3})$.
Nachmias and Peres \cite{nachmias-peres-diameter} showed that the diameter of the largest component of the critical \erdos random graph is $\Theta_P(n^{1/3})$.
There is a natural coupling between MSTs and percolation (see Observation~\ref{observation:percolation}), which together with the above result gives a matching lower bound :
\begin{align}\label{eqn:461}
\diam(M_{\infty}^{n,\er})=\Omega_P(n^{1/3}).
\end{align}
Then a stronger result was proved in \cite{AddBroGolMie13}, where the scaling limit of $M_{\infty}^{n,\er}$ viewed as a metric measure space was obtained.
We state this result in the following theorem.
\ch{We refer the reader to Section~\ref{sec:gh-mc} for the definition of the Gromov-Hausdorff-Prokhorov topology.}

\begin{thm}[Scaling limit of the MST of the complete graph \cite{AddBroGolMie13}]\label{thm:mst_complete}
View $M_{\infty}^{n,\er}$ as a random metric measure space by endowing it with the tree distance and the uniform probability measure on its vertices.
Then there exists a random compact metric measure space $\cM$ such that
\[n^{-1/3}M_{\infty}^{n,\er}\weakc\cM\]
w.r.t. Gromov-Hausdorff-Prokhorov topology.
Further, almost surely, the space $\cM$ is a binary real tree and its Minkowski dimension exists and equals $3$.
\end{thm}

Theorem~\ref{thm:mst_complete} appears to be one of the first scaling limits to be identified for any problem from combinatorial optimization, and so far, the above theorem gives the only result where the metric space scaling limit of the MST has been identified.
Several questions about the geometry of $\cM$ remain open.
For instance, what is the distribution of the typical distance in $\cM$?
More generally, is there a stick-breaking construction of $\cM$?
Is the support of the mass measure $\mu$ on $\cM$ the whole of $\cM$?
Since $\cM$ is a compact real tree, by \cite[Corollary 1.2]{duquesne2006coding}, the metric space $\cM$ (without the measure) is encoded by a random continuous function (see Section~\ref{sec:mm-spaces}). 
What can we say about the distribution of this function?

The limiting space $\cM$ is expected to be a universal object in the following sense:
For a wide array of random discrete structures that exhibit mean-field behavior, the MST constructed using i.i.d. continuous edge weights should have a rescaled version of $\cM$  as its scaling limit.
Examples of such models include the high-dimensional discrete torus, the hypercube, random regular graphs or more generally random graphs with given degree sequence (under finite third moment assumption on the degrees), various models of inhomogeneous random graphs (under appropriate assumptions), bounded-size rules, and the quantum random graph model.
See Section~\ref{sec:disc} for a more detailed discussion.

In this work, we take a first step in this broader program of establishing universality of the MST by showing that the above claim is true for the random simple $3$-regular graph and the $3$-regular configuration model.
The core of the largest component of the  \erdos random graph, in the critical window and also in the barely-supercritical regime up to a certain threshold, can be described by a $3$-regular configuration model on a random number of vertices and having random edge lengths (see \cite{janson1993birth}).
This makes the $3$-regular case special.
We use an indirect approach by exploiting the above coupling between the $3$-regular configuration model and the \erdos random graph.
However, with two additional technical estimates, our arguments can be extended to establish the scaling limit of the MST of general random graphs with given degree sequences.
We refer the reader to Section~\ref{sec:disc} for details.

\subsection{Organization of the paper}
\label{sec:org}
In Section~\ref{sec:rg-models}, we describe the random graph models considered in this paper. 
Section~\ref{sec:res} contains precise statements of our main results. 
We have deferred many definitions to Section~\ref{sec:not}, where we also give the necessary background on results on scaling limits of critical random graph models.
\ch{The proofs of two results (Theorems~\ref{thm:scaling-connected} and~\ref{thm:scaling-given-degree-sequence}) stated in Section~\ref{sec:not} are outlined in Appendix~\ref{sec:appendix}.}
In Section~\ref{sec:mst-properties} we list several properties of MSTs, and describe the connection between MSTs and percolation and the so-called cycle-breaking algorithm.
We also state a result (Theorem~\ref{thm:derandomization}) central to our argument.
In Section~\ref{sec:outline-of-proof}, we describe some of the ideas used in the proof of Theorem~\ref{thm:mst_CM}.
The proofs of our main results are given in Section~\ref{sec:proofs}.
In Section~\ref{sec:disc} we discuss the relevance of this work and related open problems.

%
%

 \subsection{Random graph models}
 \label{sec:rg-models}
 First we define the classical \erdos random graph model.
 Recall that $K_n$ denotes the complete graph on $[n]$.

\begin{defn}[The \erdos process]\label{def:erdos-renyi-random graph}
The \erdos process $\big(\ER(n, \lambda), \lambda\in\bR\big)$ is a stochastic process taking values in the space of subgraphs of $K_n$ defined as follows:
Assign a random variable $U_{ij}$ to each edge $(i, j)$ of $K_n$, where $U_{ij}$, $1\leq i<j\leq n$, are i.i.d. $\mathrm{Uniform}[0,1]$ random variables.
Set $\ER(n, \lambda)$ to be the subgraph of $K_n$ whose vertex set is $[n]$ and edge set is $\big\{(i, j)\, :\, U_{ij}\leq n^{-1}+\lambda n^{-4/3}\big\}$.
\end{defn}

\begin{rem}
\ch{
The \erdos process is often defined as a random graph process that is indexed by a parameter $p\in[0, 1]$ and takes values in the space of subgraphs of $K_n$, where the graph at parameter value $p$ has edge set $\big\{(i, j)\, :\, U_{ij}\leq p\big\}$.
We instead work with the parametrization of Definition~\ref{def:erdos-renyi-random graph}, as this will be particularly convenient for us.
}
\end{rem}

Now fix a collection of $n$ vertices labeled by $[n] := \set{1,2,\ldots, n}$ and an associated degree sequence $\vd = (d_v,\ v\in [n])$ where $\ell_n := \sum_{v\in [n]}d_v$ is assumed even. There are two natural constructions resulting in a random graph on $[n]$ with the prescribed degree sequence.

\begin{defn}[Uniformly distributed simple graphs]\label{def:uniform-simple-random-graph}
Suppose $\vd=(d_v, v\in [n])$ is a given degree sequence.
Consider the set of all simple graphs with vertex set $[n]$ where vertex $v$ has degree $d_v$, and write $\cgnd$ for the random graph having uniform distribution over this set.
 	
When $d_v=3$ for all $v\in[n]$, we will denote the corresponding random graph by $\cgnthree$.
In this case, we assume that $n$ is even.
\end{defn}
Recall that a multigraph is a graph where we allow multiple edges and self-loops.
\begin{defn}[Configuration model \cite{Boll-book,molloy1995critical,bender1978asymptotic}]\label{def:configuration-model}
Let $\cmnd$ be the random multigraph with degree sequence $\vd$ constructed sequentially as follows:
Equip each vertex $v\in [n]$ with $d_v$ half-edges or stubs.
Initially all half-edges are unpaired, and then sequentially at each step, pick two half-edges uniformly from the set of half-edges that have not yet been paired, and pair them to form a full edge.
Repeat till all half-edges have been paired.

When $d_v=3$ for all $v\in[n]$, we will denote the corresponding random multigraph by $\cmnthree$.
In this case, we assume that $n$ is even.
\end{defn}

Note that $\cmnd$ is not uniformly distributed over the set of multigraphs with degree sequence $\vd$.
We record the distribution of $\cmnd$ here for later use.
Let $G$ be a multigraph on vertex set $[n]$ in which there are $x_{ij}$ many edges between $i$ and $j$, $1\leq i<j\leq n$, and vertex $i$ has $x_{ii}$ many loops, so that
$d_i=x_{ii}+\sum_{j=1}^n x_{ij}$
is the total degree of $i$ (note that a loop contributes two to the degree). Let
$\ell_n=\sum_{i=1}^n d_i$.
Then
\begin{align}\label{eqn:cm-distribution}
\pr\big(\cmnd=G\big)=
\frac{1}{(\ell_n-1)!!}\times
\frac{\prod_{i\in[n]}d_i!}{\prod_{i\in[n]}2^{x_{ii}}\prod_{1\leq i\leq j\leq n}x_{ij}!}.
\end{align}
The proof of \eqref{eqn:cm-distribution} can be found in \cite[Proposition 7.7]{Hofs13}.

\section{Main results}
\label{sec:res}
In this section we will describe our main results. We first fix some conventions that we will follow throughout this paper.

\medskip

\noindent{\bf Convention.}  (i) For any metric measure space $\mvX=(X,d,\mu)$ and $\alpha>0$, $\alpha\mvX$ will denote the metric measure space $(X,\alpha d,\mu)$, i.e, the space where the metric has been multiplied by $\alpha$ and the measure $\mu$ has remained unchanged. Precise definitions of metric space convergence including the Gromov-Hausdorff-Prokhorov (GHP) topology are deferred to Section~\ref{sec:not}.

\vskip4pt

\noindent (ii) For any finite (not necessarily connected) graph $G$, unless the edge weights are specified,
the ``MST of $G$'' will mean the (random) minimal spanning tree of the {\bf largest component of} $G$ obtained by assigning i.i.d. continuous weights to the edges of $G$.
It is a standard fact (see Observation~\ref{observation:only-rank-matters}) that the law of the MST constructed using exchangeable edge weights that are almost surely pairwise distinct does not depend on the distribution of the underlying weights.
So the above definition of MST of $G$ makes sense.

\medskip

Recall the definitions of $\cgnthree$ and $\cmnthree$ from Section~\ref{sec:rg-models}.
Our first main result concerns the scaling limit of the MST of $\cmnthree$.
\begin{thm}[Scaling limit of the MST of the $3$-regular configuration model]\label{thm:mst_CM}
For $n$ even, let $M_n$ denote the MST of $\cmnthree$.
Think of $M_n$ as a metric measure space by using the tree distance and the uniform probability measure on the vertices.
Let $\cM$ be as in Theorem~\ref{thm:mst_complete}.
Then
\[n^{-1/3}\cdot M_n\weakc 6^{1/3}\cdot\cM\ \ \text{ as }\ \ n\to\infty\]
with respect to the Gromov-Hausdorff-Prokhorov topology.
\end{thm}

Our next main result concerns the scaling limit of the MST of $\cgnthree$.

\begin{thm}[Scaling limit of the MST of the simple $3$-regular graph]\label{thm:mst_simple}
For $n$ even, let $\overline M_n$ denote the MST of $\ \cgnthree$.
Then the result in Theorem~\ref{thm:mst_CM} continues to hold with $\overline M_n$ in place of $M_n$, i.e.,
\[n^{-1/3}\cdot\overline M_n\weakc 6^{1/3}\cdot\cM\ \ \text{ as }\ \ n\to\infty\]
with respect to the Gromov-Hausdorff-Prokhorov topology.
\end{thm}

\begin{rem}\label{rem:conditioned-to-be-connected}
Let
\begin{align}\label{eqn:def-F-n}
B_n:=\big\{\cmnthree\text{ is connected}\big\}.
\end{align}
By the results of \cite{luczak1989sparse, federico-hofstad},
\begin{align}\label{eqn:connected-whp}
\lim_{n\to\infty}\pr(B_n)
=1
=\lim_{n\to\infty}\pr\big(\cgnthree\text{ is connected}\big).
\end{align}
Thus the conclusions of the two theorems above also hold for $\cmnthree$ and $\cgnthree$  conditioned to be connected.
Further, the results of Theorems~\ref{thm:mst_CM} and~\ref{thm:mst_simple} remain true if the MST were constructed using exchangeable edge weights that are almost surely pairwise distinct.
\end{rem}

\ch{Our next result, which is interesting in its own right, is a crucial ingredient in the proofs of the above two theorems.}

\begin{thm}\label{thm:non-atomic}
	Almost surely the mass measure $\mu$ on $\cM$ is non-atomic, i.e.,
	\[\pr\big(\mu(\{x\})=0\ \text{ for every }\ x\in\cM\big)=1.\]
\end{thm}

\section{Definitions and various scaling limits}
\label{sec:not}
\subsection{Notation and conventions}\label{sec:notation}

For any set $A$, we write $|A|$ or $\# A$ for its cardinality and $\ind\set{A}$ for the associated indicator function.
For any graph $H$, we write $V(H)$ and $E(H)$ for the set of vertices and the set of edges of $H$ respectively.
We write $|H|$ for the number of vertices in $H$, i.e., $|H|=|V(H)|$.
For any finite connected graph $H=(V, E)$, we write $\mathrm{sp}(H)$ for the number of surplus edges in $H$, i.e.,
\begin{align}\label{eqn:80}
\mathrm{sp}(H):=|E|-|V|+1\, .
\end{align}

For any finite multigraph $H=(V, E)$ and $e_1,\ldots,e_k\in E$, let $H\setminus\{e_1,\ldots,e_k\}:=(V, E\setminus\{e_1,\ldots,e_k\})$.
While removing a single edge $e$ we will simply write $H\setminus e$ instead of $H\setminus\{e\}$.
Further, denote by $\conne(H)$ the set of all edges $e\in E$ such that $H\setminus e$ is connected.
For any finite multigraph $H=(V, E)$ and edges $f_1,\ldots,f_{k}$ in the complete graph on $V$, let $H\cup\{f_1,\ldots, f_k\}:=(V, E\cup\{f_1,\ldots, f_k\})$.
For two multigraphs $H_i=(V_i, E_i)$, $i=1,2$, we write $H_1\cup H_2$ for the multigraph $(V_1\cup V_2, E_1\cup E_2)$.
If $H_2$ is a connected component of $H_1$, then we write $H_1\setminus H_2$ for the multigraph $(V_1\setminus V_2, E_1\setminus E_2)$.

For any $u>0$, $\Gamma_{u}$ will denote a $\mathrm{Gamma}(u, 1)$ random variable.
We will write $\Gamma_{u}^{(\alpha)}$, $\alpha\in\Lambda$, to denote i.i.d. $\mathrm{Gamma}(u, 1)$ random variables indexed by the set $\Lambda$.

For any metric space $(X,d)$ and $U\subseteq X$, we define $\diam(U; X):=\sup\,\{d(x_1, x_2)\, :\, x_1, x_2\in U\}$.
We simply write $\diam(U)$ when there is no scope of confusion.
For any $\delta>0$ and $x\in X$, we let $B(x,\delta)=\big\{y\in X : d(x,y)\leq\delta\big\}$.
For any metric measure space $(X,d,\mu)$, we define $\fm(\delta; X):=\sup\big\{\mu\big(B(x,\delta)\big)\ :\ x\in X\big\}$.

For any tree $\vt$ on $[m]$ rooted at $\rho$, we write
\[
\height(u, \vt):=\ch{d_{\vt}}(\rho, u)\ \text{ for }\ u\in \vt,\ \text{ and }\
\height(\vt)=\max_{u\in [m]}\ \height(u, \vt)\, ,
\]
\ch{where $d_{\vt}(\cdot\, ,\, \cdot)$ denotes the tree distance on $\vt$.}
If $u\neq\rho$, we write $\stackrel{\leftarrow {\sss (1)}}{u\ \ \ }$ or simply $\overleftarrow{u}$ for the parent of $u$ in $\vt$.
If $\overleftarrow{u}\neq \rho$, then $\stackrel{\leftarrow {\sss (2)}}{u\ \ \ }$ will denote the parent of $\overleftarrow{u}$. Similarly define $\stackrel{\leftarrow {\sss (k)}}{u\ \ \ }$ for $1\leq k\leq\height(u, \vt)$.
We set $\stackrel{\leftarrow {\sss (0)}}{u\ \ \ }=u$.

For any set $A$ and function $f:A\to\bR$, we let $\|f\|_{\infty}:=\sup_{x\in A}|f(x)|$.
We use the standard Landau notation of $o(\cdot)$, $O(\cdot)$ and the corresponding \emph{order in probability} notation $o_P(\cdot)$, $O_P(\cdot)$, and $\Theta_P(\cdot)$.
We use $\probc$, $\weakc$, and $\convas$ to denote convergence in probability, convergence in distribution, and almost sure convergence respectively.
We write $\equald$ to mean equality in distribution.

\ch{When a graph with edge lengths is viewed as a metric space, the underlying set will be the collection of vertices in the graph joined by line segments (that represent the edges in the graph) of lengths specified by the edge lengths.
When not specified, all edge lengths are taken to be one.}
When a finite connected graph is viewed as a metric measure space, the measure, unless specified otherwise, will be the uniform probability measure on the vertices.

We will work with edge lengths as well as edge weights.
To avoid confusion, we make a note here that their roles are completely different.
When a graph with edge lengths is viewed as a metric space, the distances are calculated using the edge lengths. 
In Section~\ref{sec:cycle-breaking}, we will define the `cycle-breaking' process, and edge lengths will be used to perform cycle breaking.
On the other hand, edge weights are used to construct the MST (as in \eqref{eqn:def-mst}).

Throughout this paper, $C, C', c, c'$ will denote positive universal constants, and their values may change from line to line.
Special constants will be indexed as $c_1, c_2$ etc.
We freely omit ceilings and floors when there is little risk of confusion in doing so.

\subsection{Topologies on the space of metric spaces}
\label{sec:gh-mc}
We mainly follow \cite{EJP2116,AddBroGolMie13,burago-burago-ivanov, miermont2009tessellations}. All metric spaces under consideration will be compact.
For any compact $(X, d)$ and $A_1, A_2\subseteq X$, we define the Hausdorff distance between $A_1$ and $A_2$ to be
\[
d_H(A_1, A_2):=\inf\big\{\eps>0\ :\ A_1\subseteq A_2^{\eps}\ \text{ and }\ A_2\subseteq A_1^{\eps}\big\},
\]
where $A_1^{\eps}:=\bigcup_{x\in A_1}B(x, \eps)$.

Next we recall the Gromov-Hausdorff distance $d_{\GH}$ between metric spaces.  Fix two metric spaces $X_1 = (X_1,d_1)$ and $X_2 = (X_2, d_2)$. For a subset $\ch{\fR}\subseteq X_1 \times X_2$, the distortion of $\ch{\fR}$ is defined as
\begin{equation}
	\label{eqn:def-distortion}
	\dis(\ch{\fR}):= \sup \big\{|d_1(x_1,y_1) - d_2(x_2, y_2)|: (x_1,x_2) , (y_1,y_2) \in \ch{\fR}\big\}\, .
\end{equation}
A correspondence $\ch{\fR}$ between $X_1$ and $X_2$ is a measurable subset of $X_1 \times X_2$ such that for every $x_1 \in X_1$, there exists at least one $x_2 \in X_2$ such that $(x_1,x_2) \in \ch{\fR}$ and vice-versa. The Gromov-Hausdorff (GH) distance between the two metric spaces  $(X_1,d_1)$ and $(X_2, d_2)$ is defined as
\begin{equation}
\label{eqn:dgh}
	d_{\GH}(X_1, X_2) = \frac{1}{2}\inf \set{\dis(\ch{\fR}): \ch{\fR} \mbox{ is a correspondence between } X_1 \mbox{ and } X_2}.
\end{equation}
Let $\fS_{\GH}$ denote the set of isometry equivalence classes of compact metric spaces endowed with the quotient metric induced by $d_{\GH}$, which we will continue to denote by $d_{\GH}$.


We next define the marked topology; see \cite[Sections 6.4 and 6.5]{miermont2009tessellations} for a more detailed treatment.
A marked metric space is a triple $\big(X,d,\ch{Z}\big)$, where $(X,d)$ is a compact metric space and $\ch{Z}$ is a compact subset of $X$. 
The isometry classes $\big[\big(X,d,\ch{Z}\big)\big]$ of marked spaces are defined in the obvious way, and the set of such isometry classes is denoted by $\fS_{\GH}^{\ast}$.
We put the following metric on $\fS_{\GH}^{\ast}$:
For $[\mvX_i]=\big[\big(X_i, d_i, \ch{Z_i}\big)\big]\in\fS_{\GH}^{\ast }$, $i=1,2$, define
\begin{align}\label{eqn:79}
d_{\GH}^{\ast }\big([\mvX_1], [\mvX_2]\big)
:=\inf_{\phi_1,\phi_2}\bigg\{d_H\big(\phi_1(X_1), \phi_2(X_2)\big)
+d_H\big(\phi_1(\ch{Z_1}),\phi_2(\ch{Z_2})\big)\bigg\}\, ,
\end{align}
where the infimum is taken over all isometric embeddings $\phi_i: X_i\to S$, $i=1,2$, into some metric space $S$.
(There is an equivalent definition of the Gromov-Hausdorff distance $d_{\GH}$ that is similar to \eqref{eqn:79}; see, e.g., \cite[Section 7.3.2]{burago-burago-ivanov}.)

The following result is the content of \cite[Proposition 9]{miermont2009tessellations}.
\begin{lem}\label{lem:s-seq-topology}
\begin{inparaenumaa}
\item\label{item:polish}
The space $(\fS_{\GH}^{\ast}, d_{\GH}^{\ast })$ is Polish.

\noindent\item\label{item:relative-compactness}
A collection $\big\{\big[\big(X_{\alpha}, d_{\alpha}, \ch{Z_{\alpha}}\big)\big]\ :\ \alpha\in\Lambda\big\}$ is relatively compact in $(\fS_{\GH}^{\ast}, d_{\GH}^{\ast})$ iff $\big\{[(X_{\alpha}, d_{\alpha})]\ :\ \alpha\in\Lambda\big\}$ is relatively compact in $(\fS_{\GH}, d_{\GH})$, or equivalently, iff the collection of metric spaces $\big\{(X_{\alpha}, d_{\alpha})\ :\ \alpha\in\Lambda\big\}$ is uniformly totally bounded.
\end{inparaenumaa}
\end{lem}
\ch{To ease notation, we will simply write $(X_1, d_1,Z_1)$ to denote both the marked metric space  and its equivalence class.}

A compact metric measure space $(X, d, \mu)$ is a compact metric space $(X,d)$ with an associated finite measure $\mu$ on the Borel sigma algebra of $X$.
We will use the Gromov-Hausdorff-Prokhorov (GHP) distance to compare compact metric measure spaces.
Given two compact metric measure spaces $(X_1, d_1, \mu_1)$ and $(X_2,d_2, \mu_2)$ and a measure $\pi$ on the product space $X_1\times X_2$, the discrepancy of $\pi$ with respect to $\mu_1$ and $\mu_2$ is defined as

\begin{equation}
	\label{eqn:def-discrepancy}
	D(\pi;\mu_1, \mu_2):= ||\mu_1-\pi_1|| + ||\mu_2-\pi_2||~,
\end{equation}
where $\pi_1, \pi_2$ are the marginals of $\pi$ and $||\cdot||$ denotes the total variation of signed measures. Then define the distance $d_{\GHP}(X_1, X_2)$ by
\begin{equation}
\label{eqn:dghp}
	d_{\GHP}(X_1, X_2):= \inf\bigg\{ \max\bigg(\frac{1}{2} \dis(\ch{\fR}),~D(\pi;\mu_1,\mu_2),~\pi(\ch{\fR}^c)\bigg) \bigg\},
\end{equation}
where the infimum is taken over all correspondences $\ch{\fR}$ and measures $\pi$ on $X_1 \times X_2$.


The function $d_{\GHP}$ is a pseudometric and defines an equivalence relation: $X \sim Y \Leftrightarrow d_{\GHP}(X,Y) = 0$.
Let $\fS_{\GHP}$ be the set of all equivalence classes of compact metric measure spaces.
As before, we continue to denote the quotient metric by $d_{\GHP}$.
Then by \cite{EJP2116}, $(\fS_{\GHP}, d_{\GHP})$ is a complete separable metric space.
As before, to ease notation, we will continue to use $(X, d, \mu)$ to denote both the metric space and the corresponding equivalence class.

Sometimes we will be interested in not just one but an infinite sequence of compact metric measure spaces.
Then the relevant space will be $\fS_{\GHP}^{\bN}$ equipped with the product topology inherited from $d_{\GHP}$.

\subsection{Scaling limits of component sizes at criticality}
\label{sec:erdos-scaling-limit}
\ch{As we will see in the course of our proof, a key step in understanding the geometry of the MST in the supercritical regime is obtaining the metric space scaling limit of the random graph model in the critical window.
The starting point for establishing the metric space scaling limit of critical random graph models is understanding the behavior of their component sizes.}
Aldous \cite{aldous-crit} studied the maximal components of the \erdos random graph in the critical regime and  proved the following remarkable result.
\ch{Recall the notation $\mathrm{sp}(\cdot)$ from \eqref{eqn:80}.}

\begin{thm}[\cite{aldous-crit}, Corollary 2]\label{thm:erdos-renyi-component-sizes}
Write $\cC_{i}^{\sss n,\er}(\lambda)$ for the $i$-th largest connected component of $\ER(n,\lambda)$.
Then there exists a random sequence $\mvzeta(\lambda)=\bigg(\big(\xi_i(\lambda), N_i(\lambda)\big),\ i\geq 1 \bigg)$ such that as $n \to \infty$,
\[\bigg(\Big(n^{-2/3}|\cC_{i}^{\sss n,\er}(\lambda)|\, ,\, \mathrm{sp}\big(\cC_{i}^{\sss n,\er}(\lambda)\big)\Big)\, ;\ i\geq 1\bigg)
\weakc
\mvzeta(\lambda)
\]
with respect to product topology.
\end{thm}

This convergence in fact holds w.r.t. a stronger topology.
We refer the reader to \cite{aldous-crit} for an explicit description of the limiting sequence $\mvzeta(\lambda)$.
We record here a result about the asymptotic growth of the random variables $\xi_1(\lambda)$ and  $N_1(\lambda)$.

\begin{lem}\label{lem:growth-gamma-N}
	We have, as $\lambda\to\infty$,
	\[
	\frac{\xi_1(\lambda)}{\lambda}\weakc 2,\ \ \text{ and }\ \
	\frac{N_1(\lambda)}{\lambda^3}\weakc\frac{2}{3}.
	\]
\end{lem}
The proof of this result can be found in \cite[Lemma 5.6]{AddBroGolMie13}.
(See also \cite{addarioberry-bhamidi-sen} for the analogue of this result for the multiplicative coalescent in the regime where the scaling limit is a pure-jump process.)

Theorem~\ref{thm:erdos-renyi-component-sizes} has since been generalized to a number of other random graph models. In the context of graphs with given degree sequence,  Nachmias and Peres \cite{nachmias-peres-random-regular} studied critical percolation on random regular graphs;  Riordan \cite{riordan2012phase} analyzed the configuration model with bounded degrees; Joseph \cite{joseph2014component} considered i.i.d. degrees.
A more general result was obtained in \cite{dhara-hofstad-leeuwaarden-sen}. We will state a weaker version of this result next.

For a measure $\nu$ on $\mathbb{R}$ and $p > 0$,  write $\sigma_p(\nu) =\int_\mathbb{R} |x|^p \mathrm{d}\nu$; if $\nu$ has support $\bZ_{\geq 0}$ then $\sigma_p(\nu) = \sum_{i\ge 0} i^p\nu(i)$. Recall that $\nu_n\to \nu$ w.r.t. the {\em Wasserstein distance} $W_p$ if $\nu_n \to \nu$ weakly and $\sigma_q(\nu_n) \to \sigma_q(\nu)<\infty$ for all $0 \le q \le p$; see \cite[Definition 6.8]{villani-optimal}.

\begin{ass}\label{ass:cm-deg}
Suppose $\vd=\vd^{\sss(n)}=(d_v^{\sss(n)},\ v\in [n])$ is a degree sequence for each $n\geq 1$, and 
write $\nu^n := n^{-1}\sum_{v\in[n]} \delta_{d^n_v}$ for the empirical degree distribution.
Assume the following hold as $n\to\infty$:
\begin{enumeratei}
\item 
There exists a measure $\nu$ on $\bZ_{\geq 0}$ such that $\nu^n \to \nu$ w.r.t.\ the $W_3$ distance.
\item The degree sequence is in the critical scaling window, i.e., there exists $\lambda \in \bR$ such that
\begin{align}\label{eqn:987}
\frac{\sigma_1(\nu)}{\big(\sigma_3(\nu) - 4\sigma_1(\nu)\big)^{2/3}}\cdot\left(\frac{\sigma_2(\nu^n)} {\sigma_1(\nu^n)} - 2\right)\cdot n^{1/3} \to\lambda.
\end{align}
\end{enumeratei}
\end{ass}
Note that this assumption implies that $\sigma_2(\nu) = 2\sigma_1(\nu)$.
\begin{thm}[\cite{dhara-hofstad-leeuwaarden-sen}]\label{thm:cm-component-sizes}
Consider a sequence of degree sequences $\vd=\vd^{\sss(n)}$, $n\geq 1$, satisfying Assumption~\ref{ass:cm-deg} with limiting empirical distribution $\nu$.
Let $\cC_i^n$ be the $i$-th largest connected component of $\cmnd$ (or $\cgnd$).
Then as $n \to \infty$,
\begin{align}\label{eqn:598}
\Bigg(\bigg(\frac{(\sigma_3(\nu) - 4\sigma_1(\nu))^{1/3}}{\sigma_1(\nu)\cdot n^{2/3}}\cdot \big|\cC_i^n\big|,\ \mathrm{sp}\big(\cC_i^n\big)\bigg),\ i\geq 1\Bigg)\weakc \mvzeta(\lambda)
\end{align}
with respect to product topology.
\end{thm}
This result, in a stronger form, can be found in \cite[Theorem 2 and Remark 5]{dhara-hofstad-leeuwaarden-sen}.
In \cite{dhara-hofstad-leeuwaarden-sen}, the description of the  limiting sequence is slightly different.
But it is easy to restate the result in the above form using Brownian scaling.
In the next section we will use the random sequence $\mvzeta(\lambda)$ to describe certain metric measure spaces that will appear in our proofs.

\subsection{Real trees and $\bR$-graphs}\label{sec:mm-spaces}
\ch{In this section we will first define real trees and $\bR$-graphs and introduce various notions related to them.
We will then introduce a family of random $\bR$-graphs $\cH^{(s)}$, $s\geq 0$, that act as the building blocks for the scaling limits of various critical random graph models.
Using these spaces and the sequence $\mvzeta(\lambda)$ introduced in Section~\ref{sec:erdos-scaling-limit}, we will define a sequence $\mvS(\lambda)$ of random metric measure spaces; see Construction~\ref{constr:M-D}.
As we will see in the next section, the sequence $\mvS(\lambda)$ describes (up to a multiplicative constant) the scaling limits of the critical random graph models of interest to us.
}

For any metric space $(X, d)$, a geodesic between $x_1, x_2\in X$ is an isomeric embedding $f:[0, d(x_1,x_2)]\to X$ such that $f(0)=x_1$ and $f\big(d(x_1,x_2)\big)=x_2$.
$(X, d)$ is a geodesic space if there is a geodesic between any two points in $X$.
An embedded cycle in $X$ is a subset of $X$ that is a homeomorphic image of the unit circle $S^1$.
\begin{defn}[Real trees \cite{legall-survey,evans-book}]\label{def:real-tree}
A compact geodesic metric space $(X,d)$ is called a real tree if it has no embedded cycles.
\end{defn}
\begin{defn}[$\bR$-graphs \cite{AddBroGolMie13}]\label{def:r-graph}
A compact geodesic metric space $(X,d)$ is called an $\bR$-graph if for every $x\in X$, there exists $\eps\ch{=\eps(x)}>0$ such that
$\big(B(x,\eps), d\rvert_{B(x,\eps)}\big)$ is a real tree.
A measured $\bR$-graph is an $\bR$-graph with a probability measure on its Borel $\sigma$-algebra.

The core of an $\bR$-graph $(X,d)$, denoted by $\core(X)$, is the union of all the simple arcs having both endpoints in embedded cycles of $X$. If it is non-empty, then $(\core(X), d)$ is an $\bR$-graph with no leaves.
We define $\conn(X)$ to be the set of all $x\in X$ such that $x$ belongs to an embedded cycle in $X$.
\end{defn}

Clearly, $\conn(X)\subseteq\core(X)$.
By \cite[Theorem 2.7]{AddBroGolMie13}, if $X$ is an $\bR$-graph with a non-empty core, then $(\core(X), d)$ can be represented as $(k(X), e(X), \len)$, where $(k(X), e(X))$ is a finite connected multigraph in which all vertices have degree at least $3$ and $\len: e(X)\to (0,\infty)$ gives the edge lengths of this multigraph.
We denote by $\mathrm{sp}(X)$ the number of surplus edges in $(k(X), e(X))$.
On any $\bR$-graph $(X,d)$ there exists a unique $\sigma$-finite Borel measure $\ell$, called the length measure, such that if $x_1, x_2\in X$ and $[x_1,x_2]$ is a geodesic path between $x_1$ and $x_2$ then $\ell\big([x_1,x_2]\big)=d(x_1, x_2)$.
Further, we define
\begin{align}\label{eqn:def-L}
L(X):=\sum_{e\in e(X)}\len(e)=\ell(\core(X)).
\end{align}
Note that $\ell(\conn(X))\leq\ell(\core(X))<\infty$.
If $\conn(X)\neq\emptyset$ (in which case $\ell(\conn(X))>0$), we write $\ell_{\conn(X)}$ for the restriction of the length measure to $\conn(X)$ normalized to be a probability measure, i.e.,
\begin{align*}
\ell_{\conn(X)}(\cdot)=\frac{\ell(\cdot)}{\ell(\conn(X))}~.
\end{align*}

Note that any finite connected multigraph with edge lengths, viewed as a metric space, is an $\bR$-graph.
So the above definitions make sense for any finite connected multigraph $H$.
Note the difference between $e(H)$ defined above and $E(H)$-the set of edges in $H$.
Note also that in this case, the graph theoretic $2$-core of $H$, viewed as a metric space, coincides with the space $\core(H)$ as defined above.
We will use $\core(H)$ to denote both the metric space and the graph theoretic $2$-core, and the meaning will be clear from the context.
Clearly, for any finite connected multigraph $H$ with unit edge lengths, $L(H)=|E(\core(H))|$.
Further, if $H=(V, E, \len)$ is a finite connected multigraph with edge lengths, then
\begin{align}\label{eqn:8888}
\ell(H)=\sum_{e\in E}\len(e).
\end{align} 
We will write $\ell(H)$ to mean the above even when $H$ is not connected.

Functions encoding excursions from zero can be used to construct real trees via a simple procedure.
We now describe this construction.
An \emph{excursion} on $[0,1]$ is a continuous function $h \in C([0,1], \bR)$ with $h(0)=0=h(1)$ and $h(t) \geq 0$ for $t \in (0,1)$.
Let $\cE_1$ be the space of all excursions on the interval $[0,1]$. Given an excursion $h \in \cE_1$, one can construct a real tree as follows. Define a pseudo-metric $d_h$ on $[0,1]$ as follows:
\begin{equation}
\label{eqn:d-pseudo}
	d_h(s,t):= h(s) + h(t) - 2 \inf_{u \in [s,t]}h(u), \; \mbox{ for } s,t  \in [0,1].
\end{equation}
Define the equivalence relation $s \sim t \Leftrightarrow d_h(s,t) = 0$. Let $[0,1]/\sim$ denote the corresponding quotient space and consider the metric space $\cT_h:= ([0,1]/\sim, \bar d_h)$, where $\bar d_h$ is the metric on the equivalence classes induced by $d_h$. Then $\cT_h$ is a real tree (\cite{legall-survey,evans-book}).
Let $q_h:[0,1] \to \cT_h$ be the canonical projection and write $\mu_{\cT_h}$ for the push-forward of the Lebesgue measure on $[0,1]$ onto $\cT_h$ via $q_h$. Further, we assume that $\cT_h$ is rooted at $\rho := q_h(0)$.  Equipped with $\mu_{\cT_h}$, $\cT_h$ is now a rooted compact metric measure space. Note that by construction, for any $x\in \cT_h$, the function $h$ is constant on $q_h^{-1}(x)$. Thus for each $x\in [0,1]$, we write $\hght(x) = h(q_h^{-1}(x))$ for the height of this vertex.

The Brownian continuum random tree defined below is a fundamental object in the literature of random real trees.

\begin{defn}[Aldous's Brownian continuum random tree (CRT) \cite{aldous-crt-1}]\label{def:aldous-crt}
	Let $\ve$ be a standard Brownian excursion on $[0,1]$.
The real tree $\cT_{2\ve}$ 
is called the Brownian CRT.
\end{defn}
It is well-known \cite{aldous-crt-1, aldous-crt-3} that the associated measure $\mu_{\cT_{2\ve}}$ (also called the mass measure) is non-atomic and concentrated on the \ch{set} of leaves of $\cT_{2\ve}$ almost surely.
We will now define a collection of random metric measure spaces $\cH^{(s)}$, $s\geq 2$, using the Brownian CRT.
Recall the definition of $\cmnthree$ from Section~\ref{sec:rg-models}.

\begin{constr}[The space $\cH^{(s)}$ for $s\geq 2$]\label{constr:S-k-alternate}
Fix an integer $s\geq 2$, and let $n=2(s-1)$ and $r=3(s-1)$.
\begin{enumeratea}
\item Let $\cK_{n,3}$ be distributed as $\cmnthree$ conditioned to be connected.
Label its edges arbitrarily as $(u_i, v_i)$, $1\leq i\leq r$.
\item Independently of the above, sample $(X_1,\ldots,X_r)$ from a $\mathrm{Dirichlet}(\frac{1}{2},\ldots, \frac{1}{2})$ distribution.
\item Independently of the above, sample i.i.d. Brownian CRTs $\cT_1,\ldots,\cT_r$.
For $1\leq i\leq r$,
let $\rho_i$ be the root of $\cT_i$ and $z_i$ be a point in $\cT_i$ sampled according to its mass measure.
\item For $1\leq i\leq r$,
construct the metric measure space $\cT_i'$ from $\cT_i$ by multiplying the distance between each two points by $\sqrt{X_i}$
and multiplying the measure of each Borel set by $X_i$. Denote the points in $\cT_i'$ that correspond to $\rho_i$ and $z_i$ by $\rho_i'$ and $z_i'$ respectively.
\item Form a new space $\cH^{(s)}$ from $\cK_{n,3}$ by replacing the edge $(u_i, v_i)$ by $\cT_i'$ identifying $\rho_i'$ with $u_i$ and $z_i'$ with $v_i$, $1\leq i\leq r$.
\end{enumeratea}
\end{constr}
This construction of $\cH^{(s)}$ was given in \cite[Procedure 1]{BBG-limit-prop-11}.
Note that $\core(\cH^{(s)})$ is given by the multigraph $\cK_{n,3}$ with associated edge lengths $d_{\cT_i'}(\rho_i', z_i')=\sqrt{X_i}\cdot d_{\cT_i}(\rho_i, z_i)$, $1\leq i\leq r$, where $d_{\cT_i}$ denotes the metric on $\cT_i$.

For $s=0$, we define the space $\cH^{(0)}$ to be the Brownian CRT $\cT_{2\ve}$.
The explicit construction of the space $\cH^{(1)}$ is not relevant to our proof, so we do not include it here, and instead refer the reader to \cite[Procedure 1]{BBG-limit-prop-11}.
Let us also mention here that there are two other constructions of $\cH^{(s)}$.
In Construction~\ref{constr:S-k} below, we describe a `depth-first construction' of $\cH^{(s)}$.
This construction was essentially contained in the arguments in \cite{BBG-12}.
An alternate construction that can be viewed as a `breadth-first construction' is given in \cite[Construction 2.2]{miermont-sen}.


Now, we will define a sequence 
$\mvS(\lambda)=\big(S_1(\lambda), S_2(\lambda),\ldots\big)$ 
of random metric measure spaces.
Recall the random sequence $\mvzeta(\lambda)$ from Theorem~\ref{thm:cm-component-sizes}.
\begin{constr}[The sequence $\mvS(\lambda)$]\label{constr:M-D}
Sample $\mvzeta(\lambda)=\big(\big(\xi_i(\lambda), N_i(\lambda)\big),\ i\geq 1\big)$.
For simplicity, write
$\xi_i=\xi_i(\lambda),\ \text{ and }\ N_i=N_i(\lambda).$
Conditional on $\mvzeta(\lambda)$, construct the spaces $S_i(\lambda)$ independently for $i\geq 1$, where
\[S_i(\lambda)\equald\sqrt{\xi_i}\cdot \fN^{(N_i)}.\]
Set
$\mvS(\lambda)=\big(S_1(\lambda), S_2(\lambda),\ldots\big)$.
\end{constr}
Note that the spaces $\cH^{(s)}$ and $S_i(\lambda)$, $i\geq 1$, are $\bR$-graphs (recall Definition~\ref{def:r-graph}).

\subsection{Geometry of critical random graphs}
\ch{In this section, we will state four results on the geometry and scaling limit of critical random graphs that will be pivotal in our proofs.}
\begin{thm}[Geometry of uniform connected graphs with a given surplus]\label{thm:scaling-connected}
Fix an integer $s\geq 2$.
Let $\cH_{m, s}$ be uniformly distributed over the set of all simple connected graphs on $[m]$ having surplus $s$.
Recall the notation $\big(k(\cdot),\ e(\cdot),\ \len\big)$ and $L(\cdot)$ introduced around \eqref{eqn:def-L}.
Let $r=3(s-1)$.
\begin{enumeratea}
\item\label{thm:length-convergence}
We have,
\begin{align}\label{eqn:713}
\lim_{m\to\infty}\pr\Big(\big(k(\cH_{m,s}), e(\cH_{m,s})\big)\text{\ is a\ \ }3\text{-regular multigraph}\Big)=1.
\end{align}
In particular,
\begin{align}\label{eqn:714}
\lim_{m\to\infty}\pr\big(|e(\cH_{m, s})|=r\big)=1.
\end{align}
Let $e_1^{(m)},\ldots, e_{r}^{(m)}$ (resp. $e_1,\ldots, e_{r}$) be an enumeration of $\big\{e : e\in e(\cH_{m, s})\big\}$ (resp. $\big\{e : e\in e(\cH^{(s)})\big\}$).
Then as $m\to\infty$,
\begin{align}\label{eqn:73}
\Big(\frac{1}{\sqrt{m}}\cH_{m, s}~,~ \frac{1}{\sqrt{m}}\cdot\big(\len(e_i^{(m)}),\ 1\leq i\leq r\big)\Big)
\weakc
\Big(\cH^{(s)},~ \big(\len(e_i),\ 1\leq i\leq r\big)\Big),
\end{align}
where the convergence in the first coordinate is with respect to GHP topology.
Further, for any $\alpha>0$,
\begin{align}\label{eqn:74}
\sup_{m}\ \bE\big[\exp\big(\alpha m^{-1/2}L(\cH_{m, s})\big)\big]<\infty.
\end{align}
As a consequence of \eqref{eqn:73}, for every $\eps>0$, there exists $r_{\eps}>0$ such that for all large $m$,
\begin{align}\label{eqn:71}
m^{-1/2}L(\cH_{m, s})\leq 1/r_{\eps},\ \ \text{ and }\ \
m^{-1/2}\min_{e\in e(\cH_{m, s})}\len(e)\geq r_{\eps}
\end{align}
with probability at least $1-\eps$.

\item\label{thm:pendant-tree-convergence}
Let $V_i^{(m)}$ be the set of vertices in $\cH_{m, s}$ that are connected to $\core(\cH_{m,s})$ via $e_i^{(m)}$, $1\leq i\leq r$. (The common endpoints of multiple $e\in e(\cH_{m, s})$ and their pendant subtrees are assigned to only of the $V_i^{(m)}$'s in an arbitrary way.)
Recall the real trees $\cT_i'$, $1\leq i\leq r$, from Construction~\ref{constr:S-k-alternate}.
Denote the measure on $\cH^{(s)}$ by $\mu^{(s)}$.
Then as $m\to\infty$,
\begin{align}\label{eqn:711}
\frac{1}{m}\big(|V_i^{(m)}|,\ 1\leq i\leq r\big)\weakc
\big(\mu^{(s)}(\cT_i'),\ 1\leq i\leq r\big)\sim\mathrm{Dirichlet}\big(\frac{1}{2},\ldots, \frac{1}{2}\big).
\end{align}
\end{enumeratea}
\end{thm}
\eqref{eqn:713} follows from \cite[Theorem 7]{janson1993birth}.
\eqref{eqn:714} follows from \eqref{eqn:713} and the fact $\mathrm{sp}(k(\cH_{m, s}))=s$.
The rest of the assertions can be proved by following the arguments used in \cite{BBG-12}.
An outline of the proof is given in Section~\ref{sec:appendix-uniform-connected}.

\begin{thm}[Scaling limit of $\ER(n,\lambda)$]\label{thm:scaling-erdos-renyi}
Fix $\lambda\in\bR$, and let $\cC_{i}^{\sss n,\er}(\lambda)$ denote the $i$-th largest component of $\ER(n,\lambda)$.
Then
\[n^{-1/3}\big(\cC_{1}^{\sss n,\er}(\lambda),\cC_{2}^{\sss n,\er}(\lambda),\ldots\big)\weakc\mvS(\lambda)=\big(S_1(\lambda),S_2(\lambda),\ldots\big)\]
with respect to the product topology on $\fS_{\GHP}^{\bN}$ as discussed at the end of Section~\ref{sec:gh-mc}.
\end{thm}
This result is the content of \cite[Theorem 2]{BBG-12}.
That the limiting sequence of spaces is same as $\mvS(\lambda)$ follows from the discussion around
\cite[Equation 1]{BBG-limit-prop-11}.

In \cite[Theorem 2.2]{bhamidi-sen}, the metric space scaling limit of random graphs with a critical degree sequence was established.
(See also \cite{bhamidi-broutin-sen-wang}, where a similar result for critical percolation on the supercritical configuration model was derived as an application of a more general universality principle.)
The next result gives a variant of \cite[Theorem 2.2]{bhamidi-sen}.
This result follows from arguments similar to those used in \cite{bhamidi-sen}.
A sketch of proof is given in Section~\ref{sec:appendix-given-degree-sequence}.

\begin{thm}\label{thm:scaling-given-degree-sequence}
Suppose $\big\{\vd^{\sss(n)}\big\}_{n\geq 1}$ is a sequence of degree sequences satisfying Assumption~\ref{ass:cm-deg} with limiting empirical distribution $\nu$.
Further, suppose $f:\{1,2,\ldots\}\to [0,\infty)$ satisfies
$\max\big\{f(d_v)\ :\ v\in [n]\big\}=o(n^{2/3})$ and $\sum_{k\geq 1}f(k)\nu(k)>0$.
\begin{enumeratei}
\item
Let $\cC_j^n$ be the $j$-th largest component of $\cmnd$.
View $\cC_j^n$ as a metric space in the usual way; 
further, assign  mass $f(d_v)$ to each $v\in V(\cC_j^n)$ and normalize it to make it a probability measure.
(If $\sum_{v\in\cC_j^n}f(d_v)=0$ then simply take the uniform measure on the vertices.)
Denote the resulting metric measure space by $\cC_j^{n, f}$.
Then
\[n^{-1/3}\big(\cC_1^{n, f}, \cC_2^{n, f},\ldots\big)\weakc
\frac{\sigma_1(\nu)}{\big(\sigma_3(\nu)-4\sigma_1(\nu)\big)^{2/3}}\cdot\mvS(\lambda)\]
with respect to the product topology on $\fS_{\GHP}^{\bN}$ jointly with the convergence in \eqref{eqn:598}.
\item The conclusion of part (i) continues to hold with the same limiting sequence of metric measure spaces if we replace $\cmnd$ by $\cgnd$.
\end{enumeratei}
\end{thm}
Next we state a result about the core of the components of a critical graph with given degree sequence.

\begin{thm}\label{thm:belongs-to-A-r}
Suppose $\big\{\vd^{\sss(n)}\}_{n\geq 1}$ is a sequence of degree sequences satisfying Assumption~\ref{ass:cm-deg} with limiting empirical distribution $\nu$.
Let $\cC_{1}^n$ denote the largest component of $\cmnd$, and write $\cE_{1}^n=|E(\cC_{1}^n)|$.
We will drop the superscript $n$ for convenience.
Let $\Gamma_1^{(1)}, \Gamma_1^{\sss(2)},\ldots$ be i.i.d. Exponential$(1)$ random variables independent of $\cmnd$.
\begin{enumeratea}
\item
Recall the notaion $\big(k(\cdot),\ e(\cdot),\ \len\big)$ and $L(\cdot)$ introduced around \eqref{eqn:def-L}.
Then
\[n^{-1/3}\big(L(\cC_{1}),
\min_{e\in e(\cC_{1})}\len(e)\big)
\weakc
\frac{\sigma_1(\nu)}{\big(\sigma_3(\nu)-4\sigma_1(\nu)\big)^{2/3}}\cdot\big(L(S_1(\lambda)), \min_{e\in e(S_1(\lambda))}\len(e)\big).\]
In particular, for every $\eps>0$, there exists $r_{\eps}>0$ such that for all large $n$,
\[\mathrm{sp}(\cC_{1})\leq 1/r_{\eps},\ \ \ \frac{L(\cC_{1})}{n^{1/3}}\leq 1/r_{\eps},\ \ \text{ and }\ \
\min_{e\in e(\cC_{1})}\frac{\len(e)}{n^{1/3}}\geq r_{\eps}\]
with probability at least $1-\eps$.
\item
Assign lengths $\Gamma_1^{\sss (1)},\ldots, \Gamma_1^{\sss(\cE_{1})}$ to the edges of $\ \cC_{1}$, and call the resulting graph with edge lengths $\cC_{1}^{\exp}$.
Then the conclusion in (a) continues to hold with $\cC_{1}^{\exp}$ in place of $\ \cC_{1}$.
\end{enumeratea}
\end{thm}
By Theorem~\ref{thm:cm-component-sizes}, $\mathrm{sp}(\cC_{1}^n)\weakc  N_1(\lambda)$.
The other claims in Theorem~\ref{thm:belongs-to-A-r}(a) follow from the arguments used in the proof of \cite[Theorem 2.4]{bhamidi-sen}.
The claim in (b) can be proved in an identical manner.

\section{Properties of minimal spanning trees}
\label{sec:mst-properties}
In this section we discuss various properties of MSTs and give another description of the space $\cM$ appearing in Theorem~\ref{thm:mst_complete}.
\subsection{MST and percolation}
\label{sec:prim}
Suppose $G=(V, E, w)$ is a weighted, connected, and labeled graph.
Assume that $w(e)\neq w(e')$ whenever $e\neq e'$.
We now state a useful property of the MST.




\begin{lemma}[Minimax paths property]\label{lem:mst-minimax-criterion}
Let $G=(V,E,w)$ be as above.
Then the MST $T$ of $G$ is unique.
Further, $T$ has the following property:
Any path $(x_0,\hdots,\ch{x_m})$ with $x_i\in V$ and $\{x_i,x_{i+1}\}\in E(T)$
satisfies
\[\max_i\ w\big(\{x_i,x_{i+1}\}\big)\leq \max_j\ w\big(\{x_j',x_{j+1}'\}\big)\]
for any path $(x_0',\hdots,\ch{x'_{m'}})$ with $\{x_j',x_{j+1}'\}\in E$
and $x_0=x_0'$ and $\ch{x_m=x'_{m'}}$.
In words, the maximum edge weight in the path in the MST connecting two given vertices is smallest among all paths in $G$ connecting those two vertices.

Moreover, $T$ is the only spanning tree of $G$ with the above property.
\end{lemma}
The above lemma is just a restatement of \cite[Lemma 2]{kesten-lee};
see also \cite[Proposition 2.1]{alexander-forest}.
We record the following useful observations:
\begin{obs}\label{observation:only-rank-matters}
Using Lemma~\ref{lem:mst-minimax-criterion}, we see that the MST can be constructed just from the ranks of the different edge weights.
Thus the law of the MST constructed using exchangeable edge weights that are almost surely pairwise distinct does not depend on the distribution of the weights.
\end{obs}
\begin{obs}\label{observation:percolation}
Let $G=(V, E, w)$ be a connected and labeled graph with pairwise distinct edge weights.
Let $u\in[0,\infty)$ and $\cC$ be a component of the graph $G^{u}=(V, E^u)$,
where $E^u\subseteq E$ contains only those edges $e$ for which $w(e)\leq u$.
Then the restriction of the MST of $(V, E, w)$ to $\cC$ is the MST of $\big(V(\cC), E(\cC), w|_{E(\cC)}\big)$.
This can be argued as follows:
If $v,v'\in\cC$, then there exists a path in $G$ connecting $v$ and $v'$ such that all edge weights along this path is at most $u$.
By Lemma~\ref{lem:mst-minimax-criterion}, it follows that all edge weights in the path in the MST of $(V, E, w)$ connecting $v$ and $v'$ is also smaller than $u$.
Thus the restriction of the MST of $(V, E, w)$ to $\cC$ is a spanning tree of $\cC$.
Since the restriction of the MST of $(V, E, w)$ to $\cC$ also satisfies the minimax path property, it is the MST of $\cC$ (constructed using the restriction of the weight function $w(\cdot)$ to the edges of $\cC$).
This fact is extremely useful as it can be used to connect the structure of the MST to the geometry of components of the graph under percolation.
\end{obs}
\begin{obs}\label{observation:cycle-breaking}
Let $G=(V, E, w)$ be a connected and labeled graph with pairwise distinct edge weights.
Recall the notation $\conne(\cdot)$ from Section~\ref{sec:notation}.
Let $e\in\conne(G)$ be the edge with the maximum weight among all edges in $\conne(G)$.
Then $G'=(V, E\setminus\{e\}, w')$ is connected, where $w'$ is the restriction of $w$ to $E\setminus\{e\}$.
Further, by Lemma~\ref{lem:mst-minimax-criterion}, $e$ is not contained in the MST of $G$.
Thus, the MST of $G'$ is the same as the MST of $G$.
We can use this algorithm inductively to remove edges until we are left with a tree, and this tree will be the MST of $G$.
\end{obs}

\subsection{Cycle-breaking and modified cycle-breaking}
\label{sec:cycle-breaking}
In this section we define two procedures that can be applied to $\bR$-graphs and multigraphs.
Recall the definitions related to $\bR$-graphs from Section~\ref{sec:mm-spaces}.

\begin{defn}[Cycle-breaking ($\cb$)]\label{def:algo-cb}
Let $X$ be an $\bR$-graph. If $X$ has no embedded cycles, then set $\cb(X)=X$.
Otherwise, sample $x$ from $\conn(X)$ using the measure $\ell_{\conn(X)}$,
and set $\cb(X)$ to be the completion of 
\ch{the space $X\setminus \{x\}$ endowed with the intrinsic metric inherited from the metric on $X$}. 
(Thus, $\cb(X)$ is also an $\bR$-graph.)

For $k\geq 2$, we inductively define $\cb^k(X)$ to be the space $\cb\big(\cb^{k-1}(X)\big)$.
(Thus, at the $k$-th step, if $\cb^{k-1}(X)$ has an embedded cycle, then we are using the measure to $\ell_{\conn(\cb^{k-1}(X))}$ to sample a point.)

Note that $\cb^k(X)=\cb^{\mathrm{sp}(X)}(X)$ for all $k\geq\mathrm{sp}(X)$, i.e., the spaces $\cb^k(X)$ remain the same after all cycles have been cut open.
We denote this final space (which is a real tree) by $\cb^\infty(X)$.
\end{defn}

Next we define a cycle-breaking process for discrete multigraphs.
We will use a variation of the above process.
More precisely, we will sample edges with replacement.
This will turn out to be convenient in our proof.

\begin{defn}[Cycle-breaking for discrete graphs ($\cbd$)]\label{def:algo-cbd}
Let $H=(V, E, \len)$ be a finite (not necessarily connected) multigraph with edge lengths given by the function $\len:E\to (0,\infty)$.
Set $\cbd_0(H)=H$.
For $k\geq 1$, we inductively define $\cbd_k(H)$ as follows:
Sample $e_k$ from $E$ with probability proportional to $\len(e_k)$.
If $e_k$ is not an edge in $\cbd_{k-1}(H)$, set $\cbd_{k}(H)=\cbd_{k-1}(H)$.
Otherwise, if $\cC$ is the component of $\cbd_{k-1}(H)$ containing $e_k$ and
$e_k\in\conne(\cC)$, then set $\cbd_k(H)=\cbd_{k-1}(H)\setminus e_k$;
and if $e_k\notin\conne(\cC)$, then sample a point $x$ uniformly on the edge $e_k$ and color $x$ red, and set $\cbd_{k}(H)$ to be $\cbd_{k-1}(H)$ with the point $x$ colored red.

Ignoring the colored points, the multigraphs $\cbd_k(H)$ are the same (and are all forests) for all large values of $k$.
We denote the tree (without any colored points) in this forest with the most number of vertices by $\cbd_\infty(H)$.
\end{defn}

Suppose $H$ is a finite connected multigraph with edge lengths.
Let $f_1,\ldots,f_s$ be the edges of $H$ that get removed in the process $\big(\cbd_k(H), k\geq 1\big)$. Clearly, $s=\mathrm{sp}(H)$.
For $1\leq i\leq s$, let $y_i$ be a uniformly sampled point on $f_i$.
It is easy to see that viewing $H$ as an $\bR$-graph, the completion of the space $H\setminus\{y_1,\ldots,y_s\}$ has the same distribution as $\cb^\infty(H)$.
In this coupling, $\cbd_\infty(H)$ is a subspace of $\cb^\infty(H)$, and
\begin{align}\label{eqn:cbd-cb-close}
d_H\big(\cbd_\infty(H),\ \cb^\infty(H)\big)\leq\max_{e\in E}\ \len(e).
\end{align}
Further, suppose $G_1$ (resp. $G_2$) is a finite connected graph with edge lengths and $u_1$ (resp. $u_2$) is one of its vertices.
Denote by $(G_1, u_1)\stackrel{a}{\hbox{---}}(u_2, G_2)$ the graph obtained by joining $u_1$ and $u_2$ by an edge of length $a$.
Then
\begin{align}\label{eqn:cbd-joined-by-an-edge}
\cbd_\infty\big((G_1, u_1)\stackrel{a}{\hbox{---}}(u_2, G_2)\big)
\equald
\big(\cbd_\infty(G_1), u_1\big)\stackrel{a}{\hbox{---}}\big(u_2, \cbd_{\infty}(G_2)\big).
\end{align}

We now record a useful observation that we will use in the proofs.
The proof of this result is elementary, so we omit it.

\begin{lem}\label{lem:cycle-breaking-uniformly-distributed}
Suppose $H=(V,E,\len)$ is a finite multigraph with edge lengths.
\begin{enumeratea}
\item\label{item:distinct-edges}
Assume that $\len(e)$, $e\in E$, are exchangeable random variables.
For $1\leq i\leq |E|$, let $\cE_i$ denote the $i$-th distinct edge sampled in the process $\big(\cbd_k(H),\ k\geq 1\big)$.
Then for any $j\in\{1,\ldots, |E|-1\}$ and collection of distinct edges $e_1,\ldots,e_j$, conditional on the event $\{\cE_i=e_i\text{ for }1\leq i\leq j\}$, $\cE_{j+1}$ is uniformly distributed over $E\setminus\{e_1,\ldots, e_j\}$.
\item\label{item:conne-H}
Assume that $H$ is connected and that $\len(e)$, $e\in\conne(H)$, are exchangeable random variables.
For $1\leq i\leq \mathrm{sp}(H)$, let $\cE_i'$ denote the $i$-th edge removed in the process $\big(\cbd_k(H),\ k\geq 1\big)$.
Consider $j\in\{1,\ldots, \mathrm{sp}(H)-1\}$ and a collection of edges $e_1,\ldots,e_j$ satisfying $e_i\in\conne\big(H\setminus\{e_1,\ldots,e_{i-1}\}\big)$ for all $1\leq i\leq j$.
Then conditional on the event $\{\cE_i'=e_i\text{ for }1\leq i\leq j\}$, $\cE_{j+1}'$ is uniformly distributed over $\conne\big(H\setminus\{e_1,\ldots,e_j\}\big)$.
\end{enumeratea}
\end{lem}

For any finite multigraph $H=(V, E, \len)$ having edge lengths (and possibly points colored red on its edges), we write $\shape[H]$ to denote the multigraph $(V, E)$ (without any red points).
We also define $\remove(H)$ to be the multigraph with edge lengths obtained by removing all edges of $H$ that have at least one red point on them.
We now state a lemma that connects cycle-breaking to MSTs.

\begin{lem}\label{lem:cycle-breaking-gives-mst}
Suppose $H=(V,E,\len)$ is a finite connected multigraph with random edge lengths. Assume that $\len(e)$, $e\in\conne(H)$, are exchangeable random variables.
Then $\shape[\cbd_\infty(H)]$ has the same law as the MST of $\shape[H]$ constructed by assigning exchangeable pairwise distinct weights to the edges in $\conne(H)$ and any arbitrary weights to the other edges.
\end{lem}
Note that in the setting of Lemma~\ref{lem:cycle-breaking-gives-mst}, $\shape[\cbd_\infty(H)]$ is {\bf not} the MST of the weighted graph $(\shape[H], w)$ where $w(e)=\len(e)$, even though they have the same law provided the edge lengths are almost surely pairwise distinct.

\medskip

\noindent{\bf Proof of Lemma~\ref{lem:cycle-breaking-gives-mst}:}
Let $\cE_j'$ be the $j$-th edge removed in the $\cbd$ process.
Then by \ch{Lemma}~\ref{lem:cycle-breaking-uniformly-distributed}~\eqref{item:conne-H}, $\cE_1'$ is uniformly distributed over $\conne(H)$.
In general, conditional on $\cE_i'$, $1\leq i\leq k-1$, $\cE_k'$ is uniformly distributed over $\conne\big(H\setminus\{\cE_1',\ldots,\cE_{k-1}'\}\big)$.

Now, consider edge weights $(w(e),\ e\in E)$, such that $w(e)$, $e\in\conne(H)$, are exchangeable and almost surely pairwise distinct.
Then using Observation~\ref{observation:cycle-breaking}, the MST of $(\shape[H], w)$ can be constructed by sequentially removing the edges having maximum weight among all edges whose removal do not disconnect the current graph.
By the assumptions on the weights, the edge to be removed at each step is uniformly distributed over the set of all edges whose removal do not disconnect the current graph.
In other words, the sequence of edges removed in the algorithm described in Observation~\ref{observation:cycle-breaking} has the same law as $\big(\cE_k',~ k\geq 1\big)$.
This completes the proof.
\qed

\medskip

Recall the notation $k(X)$, $e(X)$, $(\len(e),\ e\in e(X))$, $\mathrm{sp}(X)$, and $L(X)$ introduced below Definition~\ref{def:r-graph}.
For $r\in (0,1)$ define $\cA_r$ to be the set of all measured $\bR$-graphs $X$ that satisfy
\[\mathrm{sp}(X)+L(X)\leq 1/r,\ \text{ and }\ \min_{e\in e(X)}\len(e)\geq r.\]
The following theorem will allow us to prove convergence of MSTs from GHP convergence of the underlying graphs.

\begin{thm}\label{thm:mst-from-ghp-convergence}
Fix $r\in (0,1)$.
Suppose $(X,d,\mu)$ and $(X_n, d_n, \mu_n)$, $n\geq 1$, are measured $\bR$-graphs in $\cA_r$ such that $(X_n, d_n, \mu_n)\to (X,d,\mu)$ as $n\to\infty$ w.r.t. GHP topology.
\begin{enumeratea}
\item\label{item:cb-infty-convegence}
Then $\cb^{\infty}(X_n)\weakc\cb^{\infty}(X)$ as $n\to\infty$ w.r.t. GHP topology.
\item\label{item:cbd-infty-convegence}
Suppose for each $n\geq 1$, $(X_n, d_n,\mu_n)$ is the metric measure space associated with $(V_n, E_n,\len)$--a finite connected multigraph with edge lengths.
If $\max_{e\in E_n}\len(e)\to 0$ as $n\to\infty$, then $\cbd_{\infty}(X_n)\weakc\cb^{\infty}(X)$ as $n\to\infty$ w.r.t. GHP topology.
\end{enumeratea}
\end{thm}
The result in Theorem~\ref{thm:mst-from-ghp-convergence}~\eqref{item:cb-infty-convegence} is from \cite[Theorem 3.3]{AddBroGolMie13}, while the claim in \eqref{item:cbd-infty-convegence} follows from \eqref{eqn:cbd-cb-close}.

\subsection{Alternate descriptions of the space $\cM$}\label{sec:alternate-description-M}
Recall the construction of the process $\ER(n,\cdot)$ using the random variables $U_{ij}$ from Definition~\ref{def:erdos-renyi-random graph}.
Let $\cC_{1}^{\sss n,\er}(\lambda)$ be the largest component of $\ER(n,\lambda)$ and let $M_{\lambda}^{n,\er}$ be the MST of $\cC_{1}^{\sss n,\er}(\lambda)$ constructed using the random weights $U_{ij}$, $(i,j)\in E(\cC_{1}^{\sss n,\er}(\lambda))$.
Then $\lim_{\lambda\to\infty}M_{\lambda}^{n,\er}=M_{\infty}^{n,\er}$ (in fact $M_{\lambda}^{n,\er}=M_{\infty}^{n,\er}$ for all large $\lambda$), where $M_{\infty}^{n,\er}$ is the MST of $K_n$ constructed using the random weights $U_{ij}$.
Theorem~\ref{thm:mst_complete} says that $n^{-1/3}M_{\infty}^{n,\er}\weakc\cM$ as $n\to\infty$ w.r.t. GHP topology.
The natural question to ask here is whether the order in which the limits are taken can be interchanged, i.e., can we first take limit as $n\to\infty$ for fixed $\lambda$, and then let $\lambda\to\infty$?
Now, by \cite[Theorem 4.4]{AddBroGolMie13},
\begin{align}\label{eqn:700}
n^{-1/3}M_{\lambda}^{n,\er}\weakc \cb^{\infty}\big(S_1(\lambda)\big)\ \ \text{ as }\  \  n\to\infty
\end{align} 
w.r.t. GHP topology.
Then the following theorem answers the above question in the affirmative.

\begin{thm}[\cite{AddBroGolMie13}, Theorem 4.9]\label{thm:M-alternate}
As $\lambda\to\infty$,
\[\cb^{\infty}\big(S_1(\lambda)\big)\weakc \cM\]
with respect to GHP topology.
\end{thm}

The space $S_1(\lambda)$ has a random number of cycles.
The following theorem gives a derandomized version of Theorem~\ref{thm:M-alternate}.
\begin{thm}\label{thm:derandomization}
Recall the space $\cH^{(s)}$ from Construction~\ref{constr:S-k-alternate}. Then
\[\big(12 s\big)^{1/6}\cdot\cb^{\infty}\big(\cH^{(s)}\big)\weakc \cM\ \text{ as }\ s\to\infty\]
with respect to GHP topology.
\end{thm}
Theorem~\ref{thm:derandomization} plays a crucial role in our argument.
The proof of this result can be read independently of the rest and is deferred to Section~\ref{sec:derandomization}.

\section{Idea of proof }\label{sec:outline-of-proof}
\ch{In this section, we outline the proof of the fact that the claimed convergence in Theorem~\ref{thm:mst_CM} holds w.r.t. the GH topology assuming Theorem~\ref{thm:derandomization}.
We explain the ideas at a high level, and do our best to avoid getting into the technicalities.
For any finite graph $H$, we will write $H^{\exp}$ to denote the graph obtained by assigning i.i.d. $\text{Exponential}(1)$ lengths to the edges of $H$.
When $H$ is random, the edge lengths are taken to be independent of $H$.
}

\ch{Recall the notation from Construction~\ref{constr:S-k-alternate}. 
Then
\[
	\big(X_1,\ldots,X_r\big)
	\equald
	\big(\Gamma_{1/2}^{\sss(1)},\ldots, \Gamma_{1/2}^{\sss(r)}\big)/\Gamma_{r/2}\ ,
\]
where $\Gamma_{1/2}^{\sss(j)}$, $j=1,\ldots, r$, are i.i.d. $\mathrm{Gamma}(1/2,1)$ random variables, and $\Gamma_{r/2}:=\sum_{j=1}^r \Gamma_{1/2}^{\sss(j)}=\frac{r}{2}\cdot \big(1+o_P(1)\big)$ as $s\to\infty$.
Further, it is well-known that typical distance in a Brownian CRT follows a Rayleigh distribution.
Consequently, 
$d_{\cT_i}(\rho_i, z_i)\cdot\big(2\Gamma_{1/2}^{(i)}\big)^{1/2}$, $i=1,\ldots, r$, 
are i.i.d. $\text{Exponential}(1)$  random variables.
As noted below Construction~\ref{constr:S-k-alternate}, $\core(\cH^{(s)})$ can be represented by the multigraph $\cK_{n,3}$ with edge lengths given by $d_{\cT_i}(\rho_i, z_i)\cdot\sqrt{X_i}$, $i=1,\ldots, r$.
Hence, $\core(\cH^{(s)})$ is simply $\gamma_n\cdot\cK_{n,3}^{\exp}$, where $\gamma_n=r^{-1/2}\big(1+o_P(1)\big)$ as $s\to\infty$.
Now, the space $(12s)^{1/6}\cdot\cH^{(s)}$ can be obtained from $(12s)^{1/6}\cdot\core(\cH^{(s)})$ by attaching some random compact trees.
As $s\to\infty$, the maximum diameter of these trees becomes negligible.
In other words, the result in Theorem~\ref{thm:derandomization} continues to hold if we replace $\cH^{(s)}$ by $\core(\cH^{(s)})$, which can in turn be replaced by $r^{-1/2}\cdot\cK_{n,3}^{\exp}$.
Thus,
$(12s)^{1/6}\cdot r^{-1/2}\cdot \cb^{\infty}\big(\cK_{n,3}^{\exp}\big)
\weakc 
\cM$, as $n\to\infty$, 
with respect to the GH topology.
Using the relations $s\sim n/2$ and $r\sim 3n/2$ as $n\to\infty$, we conclude that
$n^{-1/3}\cdot\cb^{\infty}\big(\cK_{n,3}^{\exp}\big)
\weakc 
(0.75)^{1/3}\cdot\cM$, as $n\to\infty$, 
with respect to the GH topology.
Now using \eqref{eqn:cbd-cb-close} to go from $\cb^{\infty}$ to $\cbd_{\infty}$ and using \eqref{eqn:connected-whp}, we get
\begin{align}\label{eqn:808}
n^{-1/3}\cdot\cbd_{\infty}\big(\cmnthreeexp\big)
\weakc 
(0.75)^{1/3}\cdot\cM\, ,\ \text{ as }\ n\to\infty\, , 
\end{align}
with respect to the GH topology.
}

\ch{Now, observe that for any finite graph $H$, conditional on the event that the number of distinct edges sampled in the process $\big(\cbd_i(H^{\exp})\, ,\ 1\leq i\leq T \big)$ is $m$, the collection of these edges has the same distribution as a uniform subset of size $m$ sampled from $E(H)$.
From this, one might guess that for an appropriately chosen random $T$, 
$
\shape\big[\remove\big(\cbd_{T}(H^{\exp})\big)\big]
$
will have the same distribution as
$
\perc(H, p)
$
for a deterministic $p$, where $\perc(H, p)$ denotes the random subgraph of $H$ obtained under percolation with edge retention probability $p$.
In fact, we have the following stronger result:
Fix $t>0$, and let $R(t)$ be a $\mathrm{Poisson}\big(t\cdot \ell(H^{\exp})\big)$ random variable, where $\ell(\cdot)$ is as in \eqref{eqn:8888}.
Then
\begin{align}\label{eqn:566}
&
\bigg( 
\shape\Big[\remove\big(\cbd_{R(t)}(H^{\exp})\big)\Big],\
\remove\Big(\cbd_{R(t)}(H^{\exp})\Big)
\bigg)
\notag\\
&\hskip90pt
\equald
\bigg( 
\perc\Big(H, \frac{1}{1+t}\Big)\, ,\
\frac{1}{1+t}\cdot\Big(\perc\Big(H, \frac{1}{1+t}\Big)\Big)^{\exp}~
\bigg)\, .
\end{align}
This is the content of Lemma~\ref{lem:removed-or-marked} whose proof is rather short.
}

\ch{For fixed $\lambda\in\bR$ and $n$ large so that $2n^{1/3}>|\lambda|$, let $t_{n,\lambda}$ be such that 
$
\big(1+t_{n,\lambda}\big)^{-1}=1/2+\lambda n^{-1/3}
$.
Write $\cmnthree(\lambda)$ for $\perc\big(\cmnthree, 1/2+\lambda n^{-1/3}\big)$.
Let $\cC_1(\lambda)$ denote the largest connected component of $\cmnthree(\lambda)$.
Applying \eqref{eqn:566} with $H=\cmnthree$ and $t=t_{n,\lambda}$ will yield the following:
Let $\fG_1(n,\lambda)$ be the largest connected component of
$\remove\big(\cbd_{R(t_{n,\lambda})}\big(\cmnthreeexp\big)\big)$, 
and set 
$G_1(n,\lambda):=\shape\big[\fG_1(n,\lambda)\big]$.
Then
\begin{align}\label{eqn:665}
\big(G_1(n,\lambda)\, ,\ \fG_1(n,\lambda)\big)
\equald
\bigg(\cC_1(\lambda)\, ,\ \big(1/2+\lambda n^{-1/3}\big)\cdot\big(\cC_1(\lambda)\big)^{\exp}\bigg)\, .
\end{align}
Now, $\cmnthree(\lambda)$, conditional on its degree sequence, is distributed as a configuration model with that degree sequence.
Further, it is easy to show that the (random) degree sequence of $\cmnthree(\lambda)$ satisfies Assumption~\ref{ass:cm-deg} with limiting empirical distribution $\nu=\text{Binomial}(3, 1/2)$.
Using these observations together with Theorem~\ref{thm:scaling-given-degree-sequence}, we can show that
\begin{align}\label{eqn:569}
n^{-1/3}\cdot \cC_1(\lambda)
\weakc 
6^{1/3}\cdot S_1\big((48)^{1/3}\cdot\lambda\big)
\end{align}
w.r.t. the GH topology. 
Now, consider any self-avoiding path $\pi$ in $\cC_1(\lambda)$.
Since the edge lengths in $\big(\cC_1(\lambda)\big)^{\exp}$ are i.i.d. $\text{Exponential}(1)$ random variables, the length of $\pi$ in $\big(\cC_1(\lambda)\big)^{\exp}$ will be concentrated around the length of $\pi$ in $\cC_1(\lambda)$.
Thus, leveraging the fact that $\cC_1(\lambda)$ has only $O_P(1)$ many surplus edges, we can show that 
$
d_{\GH}\big(\cC_1(\lambda)\, ,\ \big(\cC_1(\lambda))^{\exp}\big)=o_P(n^{1/3})
$.
Combining this with \eqref{eqn:569}, \eqref{eqn:665}, Theorem~\ref{thm:mst-from-ghp-convergence}, and Theorem~\ref{thm:belongs-to-A-r} will yield
\begin{gather*}
n^{-1/3}\cbd_{\infty}\big(G_1(n,\lambda)\big)\weakc 
6^{1/3}\cdot\cb^{\infty}\Big(S_1\big((48)^{1/3}\cdot\lambda\big)\Big)\, ,\ \text{ and}\\
n^{-1/3}\cbd_{\infty}\big(\fG_1(n,\lambda)\big)\
\weakc 
(0.75)^{1/3}\cdot\cb^{\infty}\Big(S_1\big((48)^{1/3}\cdot\lambda\big)\Big) 
\end{gather*}
w.r.t. the GH topology.
Using Theorem~\ref{thm:M-alternate}, we see that for any $\bZ_{>0}$-valued sequence \chl{$\lambda_n$ that tends to infinity} sufficiently slowly,
\begin{gather}
n^{-1/3}\cbd_{\infty}\big(G_1(n,\lambda_n)\big)\
\weakc 
6^{1/3}\cdot\cM\, ,\ \text{ and}\label{eqn:709-a}\\
n^{-1/3}\cbd_{\infty}\big(\fG_1(n,\lambda_n)\big)\
\weakc 
(0.75)^{1/3}\cdot\cM\label{eqn:709}
\end{gather}
w.r.t. the GH topology.
}

\ch{Now compare \eqref{eqn:709} with \eqref{eqn:808}.
Let $B_n$ be as in \eqref{eqn:def-F-n}.
Then on the event $B_n$ (which, by \eqref{eqn:connected-whp}, occurs with high probability), 
$\cbd_{\infty}\big(\fG_1(n,\lambda_n)\big)$ is a subspace of $\cbd_{\infty}(\cmnthree^{\exp})$.
Since they have the same scaling limit,
\begin{align}\label{eqn:710}
n^{-1/3}\cdot\ind_{B_n}\cdot d_{\rH}\big(
\cbd_{\infty}\big(\fG_1(n,\lambda_n)\big)\, ,\   
\cbd_{\infty}(\cmnthree^{\exp})
\big)\weakc 0\, .
\end{align}
This follows from a general property of metric spaces; see Proposition~\ref{prop:gh-lemma-3}.
Suppose, on the event $B_n$, $\cbd_{\infty}(\cmnthree^{\exp})$ is obtained by attaching the trees
$\fT_{n,\lambda_n}^{(j)}$, $1\leq j\leq k_n(\lambda_n)$, each
via an edge to a vertex of $\cbd_{\infty}\big(\fG_1(n,\lambda_n)\big)$.
Set $k_n(\lambda_n)=0$ on $B_n^c$.
Then \eqref{eqn:710} is equivalent to the assertion that 
$
n^{-1/3}\max_{1\leq j\leq k_n(\lambda_n)}
\diam\big(\fT_{n,\lambda_n}^{(j)}\big)\weakc 0
$.
From this, it is not difficult to argue that
\begin{align}\label{eqn:712}
	n^{-1/3}\max_{1\leq j\leq k_n(\lambda_n)}
	\diam\big( \shape\big[\fT_{n,\lambda_n}^{(j)}\big] \big)
	\weakc 0\, .
\end{align}
Finally, 
using Lemma~\ref{lem:cycle-breaking-uniformly-distributed}~\eqref{item:distinct-edges}, the processes $\big(\cbd_i(\cmnthree^{\exp}),\ i\geq 1\big)$ and $\big(\cbd_t(\cmnthree),\ t\geq 1\big)$ can be coupled so that the $j$-th distinct edge sampled is the same in both processes, $1\leq j\leq 3n/2$.
In this coupling, on the event $B_n$, $\cbd_{\infty}(\cmnthree)$ is $\cbd_{\infty}\big(G_1(n,\lambda_n)\big)$ with the trees 
$\shape\big[\fT_{n,\lambda_n}^{(j)}\big]$, $1\leq j\leq k_n(\lambda_n)$, 
attached to its vertices via an edge.
This observation together with \eqref{eqn:connected-whp}, \eqref{eqn:712}, \eqref{eqn:709-a}, and Lemma~\ref{lem:cycle-breaking-gives-mst} shows that
\[
n^{-1/3}M_n
\equald 
n^{-1/3}\cbd_{\infty}(\cmnthree)
\weakc
6^{1/3}\cdot\cM
\]
as $n\to\infty$ w.r.t. the GH topology.
}

\section{Proofs of Theorems~\ref{thm:mst_CM},~\ref{thm:mst_simple},~\ref{thm:non-atomic}, and~\ref{thm:derandomization}}\label{sec:proofs}
We divide the argument into several steps.
In Section~\ref{sec:gh-convergence-cm}, we prove a weaker version of Theorem~\ref{thm:mst_CM} that only deals with convergence w.r.t. GH topology.
The proof of this result depends on several propositions whose proofs are given in Sections~\ref{sec:brw}--\ref{sec:proof-gh-lemma-3}.
The proof of Theorem~\ref{thm:mst_CM} is then completed in Section~\ref{sec:ghp-convergence-cm}.
The proof of Theorem~\ref{thm:mst_simple} is given in Section~\ref{sec:ghp-convergence-simple}.
The proof of Theorem~\ref{thm:non-atomic} is given in Section~\ref{sec:non-atomic-proof}.
Finally, the proof of Theorem~\ref{thm:derandomization} is given in Section~\ref{sec:derandomization}.

\subsection{GH convergence of the MST of $\cmnthree$}\label{sec:gh-convergence-cm}
In this section we prove the following weaker version of Theorem~\ref{thm:mst_CM}.

\begin{thm}\label{thm:gh-mst-cm}
Let $M_n$ be as in Theorem~\ref{thm:mst_CM}.
Then $n^{-1/3}\cdot M_n\weakc \ch{6^{1/3}\cdot\cM}$ with respect to GH topology.
\end{thm}
This convergence will be strengthened to GHP convergence in Section~\ref{sec:ghp-convergence-cm}.
The proof of the above theorem relies on the following four propositions.

\begin{prop}\label{lem:comparison}
For all $r\in\bN$, there exists $c>0$ small such that the following holds:
Let $G=(V,E)$ be a finite graph with maximum degree at most $r$.
Let $\Gamma_1^{(e)}$, $e\in E$, be i.i.d. $\mathrm{Exponential}(1)$ random variables.
Then for all $m\geq 1$,
\begin{align*}
\pr\Big(G\text{ contains a self-avoiding path }P\text{ with }|P|\geq m\text{ and }\sum_{e\in P}\Gamma_1^{(e)}\leq c|P|\Big)
\leq
|V|\cdot\exp(-m),
\end{align*}
where $|P|$ denotes the number of edges in the path $P$.
\end{prop}


\begin{prop}\label{prop:gh-lemma-1}
Assign i.i.d. $\mathrm{Exponential}(1)$ lengths to the edges of $\cmnthree$ and denote this multigraph with edge lengths by $\cmnthreeexp$.
Then
\[n^{-1/3}\cbd_{\infty}\big(\cmnthreeexp\big)\weakc 
\ch{\frac{1}{2}\cdot 6^{1/3}\cdot\cM=}
(0.75)^{1/3}\cdot\cM,\ \text{ as }\ n\to\infty\]
with respect to GH topology.
\end{prop}
\begin{rem}\label{rem:diam-bound}
By Lemma~\ref{lem:cycle-breaking-gives-mst},
$M_n\equald\shape\big[\cbd_{\infty}\big(\cmnthreeexp\big)\big]$.
However, conditional on $\shape\big[\cbd_{\infty}\big(\cmnthreeexp\big)\big]$, the edge lengths of $\cbd_{\infty}\big(\cmnthreeexp\big)$ are {\bf not} exchangeable, which is why Theorem~\ref{thm:gh-mst-cm} cannot be proved by just using Proposition~\ref{prop:gh-lemma-1}, and it takes quite a bit of additional work.
Note however that Proposition~\ref{prop:gh-lemma-1} implies that
\[\diam\big(\cbd_{\infty}\big(\cmnthreeexp\big)\big)=\Theta_P(n^{1/3}).\]
This observation together with Proposition~\ref{lem:comparison} implies that $\diam(M_n)=O_P(n^{1/3})$.
As noted before in \eqref{eqn:461} in the case of the complete graph, using Observation~\ref{observation:percolation} and Theorem~\ref{thm:scaling-given-degree-sequence}, it follows that $\diam(M_n)=\Omega_P(n^{1/3})$.
Thus, we get that $\diam(M_n)=\Theta_P(n^{1/3})$.
By a standard conditioning argument (see \eqref{eqn:cm-conditional-uniform} and \eqref{eqn:cm-simple-0}), this also implies that $\diam(\overline M_n)=\Theta_P(n^{1/3})$.
\end{rem}

Recall the notation $\shape[\cdot]$ and $\remove(\cdot)$ introduced right before Lemma~\ref{lem:cycle-breaking-gives-mst}.
Recall also from \eqref{eqn:8888} and the line below the meaning of $\ell(H)$ for finite multigraphs with edge lengths.
\begin{prop}\label{prop:gh-lemma-2}
Let $S_1(\cdot)$ and $\cmnthreeexp$ be as in Construction~\ref{constr:M-D} and Proposition~\ref{prop:gh-lemma-1} respectively.
For $\lambda\in\bR$ satisfying $|\lambda|<n^{1/3}/2$, let $t_{n,\lambda}$ be given by
\begin{align}\label{eqn:6699}
\frac{1}{1+t_{n,\lambda}}=\frac{1}{2}+\frac{\lambda}{n^{1/3}}~.
\end{align}
Let $R_{n,\lambda}$ be a
$\mathrm{Poisson}\big(t_{n,\lambda}\cdot \ell(\cmnthreeexp)\big)$ random variable.
Let $\fG_1(n,\lambda)$ be the largest component of
$\remove\big(\cbd_{R_{n,\lambda}}\big(\cmnthreeexp\big)\big)$.
Let $G_1(n,\lambda):=\shape\big[\fG_1(n,\lambda)\big]$.
Then for any fixed $\lambda\in\bR$, 
\begin{gather}
n^{-1/3}\cbd_{\infty}\big(G_1(n,\lambda)\big)\weakc 
6^{1/3}\cdot\cb^{\infty}\Big(S_1\big((48)^{1/3}\cdot\lambda\big)\Big),\ \ \text{ and}\label{eqn:30}\\
n^{-1/3}\cbd_{\infty}\big(\fG_1(n,\lambda)\big)\weakc (0.75)^{1/3}\cdot\cb^{\infty}\Big(S_1\big((48)^{1/3}\cdot\lambda\big)\Big)\label{eqn:31}
\end{gather}
as $n\to\infty$ with respect to the GH topology.
\end{prop}

Recall the marked topology from Section~\ref{sec:gh-mc}.
\begin{prop}\label{prop:gh-lemma-3}
Suppose $\big\{(X_n^+, d_n, X_n)\big\}_{n\geq 1}$ is a sequence of random compact marked metric spaces such that
\[X_n^+\weakc Z,\ \text{ and }\ X_n\weakc Z,\ \ \text{ as } n\to\infty\]
with respect to the GH topology for some random compact metric space $Z$.
Then $d_H(X_n, X_n^+)\weakc 0$ as $n\to\infty$.
\end{prop}
We first prove Theorem~\ref{thm:gh-mst-cm} assuming the above four propositions.
The proofs of Propositions~\ref{lem:comparison},~\ref{prop:gh-lemma-1},~\ref{prop:gh-lemma-2}, and~\ref{prop:gh-lemma-3} will be given in the next four sections.
\ch{
We will make use of the following elementary fact in the proof of Theorem~\ref{thm:gh-mst-cm}; we omit its proof.
\begin{lem}\label{lem:23}
Suppose $a_{i,j}$, $i\in\bZ_{>0}$, $j\in\bZ_{>0}\cup\{\infty\}$, and $a_{\infty, \infty}$ are elements of some metric space such that 
$\lim_{j\to\infty} a_{i,j}= a_{i,\infty}$ for every $i\in\bZ_{>0}$, and
$\lim_{i\to\infty}a_{i,\infty}=a_{\infty,\infty}$.
Then there exists a $\bZ_{>0}$-valued sequence $\big\{i_j^{\star}\big\}_{j\in\bZ_{>0}}$ with $i_j^{\star}\uparrow\infty$ such that for any $\bZ_{>0}$-valued sequence $\big\{i_j\big\}_{j\in\bZ_{>0}}$ satisfying $i_j\uparrow\infty$ and  $i_j\leq i_j^{\star}$,
$\lim_{j\to\infty}a_{i_j, j}=a_{\infty, \infty}$.
\end{lem}
}

\noindent{\bf Proof of Theorem~\ref{thm:gh-mst-cm}:}
Let $B_n$ be as in \eqref{eqn:def-F-n}.
Note that on the event $B_n$, for any $\lambda\in\bR$, the space $\cbd_{R_{n,\lambda}}(\cmnthree^{\exp})$ is simply $\fG_1(n,\lambda)$ together with some additional connected multigraphs (with edge lengths and red points) each of which is attached to a vertex of $\fG_1(n,\lambda)$ via a single edge that has at least one red point on it.
Thus, by \eqref{eqn:cbd-joined-by-an-edge}, on the event $B_n$, $\cbd_{\infty}(\cmnthree^{\exp})$ is $\cbd_{\infty}\big(\fG_1(n,\lambda)\big)$ with some additional trees,
say $\fT_{n,\lambda}^{(j)}$, $1\leq j\leq k_n(\lambda)$, each of which is attached to a vertex of $\cbd_{\infty}\big(\fG_1(n,\lambda)\big)$ via a single edge.
Define $k_n(\lambda)=0$ on $B_n^c$ for all $\lambda\in\bR$.

Using Lemma~\ref{lem:cycle-breaking-uniformly-distributed}\eqref{item:distinct-edges}, there exists a coupling of the processes $\big(\cbd_i(\cmnthree^{\exp}),\ i\geq 1\big)$ and $\big(\cbd_t(\cmnthree),\ t\geq 1\big)$ such that the $j$-th distinct edge sampled is the same in both processes, $1\leq j\leq 3n/2$.
In this coupling, on the event $B_n$, $\cbd_{\infty}(\cmnthree)$ is $\cbd_{\infty}\big(G_1(n,\lambda)\big)$ with $\shape\big[\fT_{n,\lambda}^{(j)}\big]$, $1\leq j\leq k_n(\lambda)$, attached to its vertices via an edge.

Now using \eqref{eqn:30}, \eqref{eqn:31}, Theorem~\ref{thm:M-alternate}, \ch{and Lemma~}\ref{lem:23}, it follows that there exists a $\bZ_{>0}$-valued sequence $\{\lambda^{\star}_n\}_{n\geq 1}$ with $\lambda^{\star}_n\uparrow \infty$ such that for any $\bZ_{>0}$-valued sequence $\lambda_n\uparrow\infty$ satisfying $\lambda_n\leq\lambda^{\star}_n$,
\begin{gather}
n^{-1/3}\cbd_{\infty}\big(G_1(n,\lambda_n)\big)
\weakc 6^{1/3}\cdot\cM ,\label{eqn:34}
\ \text{ and }\\
n^{-1/3}\cbd_{\infty}\big(\fG_1(n,\lambda_n)\big)
\weakc \big(0.75\big)^{1/3}\cdot\cM \label{eqn:35}
\end{gather}
with respect to GH topology.
Using \eqref{eqn:35} in conjunction with Proposition~\ref{prop:gh-lemma-1}, Proposition~\ref{prop:gh-lemma-3}, and \eqref{eqn:connected-whp}, it follows that for any $\bZ_{>0}$-valued sequence $\lambda_n\uparrow\infty$ with $\lambda_n\leq\lambda_n^{\star}$,
\begin{align}\label{eqn:36}
n^{-1/3}\max_{1\leq j\leq k_n(\lambda_n)}\diam\big(\fT_{n,\lambda_n}^{(j)}\big)\weakc 0.
\end{align}
Denoting the edge lengths of $\cmnthreeexp$ by $\Gamma_1^{(e)}$, $e\in E(\cmnthree)$, we have, for any $\eps>0$ and any $c>0$,
\begin{align*}
&\pr\Big(\max_{1\leq j\leq k_n(\lambda_n)}\diam\big(\shape\big[\fT_{n,\lambda_n}^{(j)}\big]\big)\geq\eps n^{1/3}\Big)
\leq\pr\Big(\max_{1\leq j\leq k_n(\lambda_n)}\diam\big(\fT_{n,\lambda_n}^{(j)}\big)\geq c\eps n^{1/3}\Big)\\
&\hskip35pt
+\pr\Big(\cmnthree\text{ contains a self-avoiding path }P\text{ with }|P|\geq \eps n^{1/3}
\text{ and }\sum_{e\in P}\Gamma_1^{(e)}\leq c\eps n^{1/3}\Big).
\end{align*}
Thus, using Proposition~\ref{lem:comparison} together with \eqref{eqn:36}, we get
\begin{align}\label{eqn:37}
n^{-1/3}\max_{1\leq j\leq k_n(\lambda_n)}\diam\big(\shape\big[\fT_{n,\lambda_n}^{(j)}\big]\big)\weakc 0,
\end{align}
which in turn shows that
\begin{align}\label{eqn:38}
n^{-1/3}d_H\Big(\cbd_{\infty}\big(\cmnthree\big),\ \cbd_{\infty}\big(G_1(n,\lambda_n)\big)\Big)\weakc 0.
\end{align}
Finally, by Lemma~\ref{lem:cycle-breaking-gives-mst},
\begin{align}\label{eqn:679}
M_n\equald
\cbd_{\infty}(\cmnthree).
\end{align}
The result now follows from \eqref{eqn:34} and \eqref{eqn:38}.
\qed

\subsection{Proof of Proposition~\ref{lem:comparison}}\label{sec:brw}
Fix $m\geq 1$ and $k\geq m$. 
Consider a self-avoiding path $P$ in $G$ with $|P|=k$.
Then for any $c>0$ and any $t>0$,
\begin{align}\label{eqn:4884}
\pr\Big(\sum_{e\in P}\Gamma_1^{(e)}\leq ck\Big)
=\pr\big(Z\geq k\big)
\leq e^{-tk}\bE\big[\exp(tZ)\big]
=\exp\Big(-k\big(t-c(e^t-1)\big)\Big),
\end{align}
where $Z$ is a $\mathrm{Poisson}(ck)$ random variable.
Hence,
\begin{align*}
&\pr\Big(G\text{ contains a self-avoiding path }P\text{ with }|P|\geq m\text{ and }\sum_{e\in P}\Gamma_1^{(e)}\leq c|P|\Big)\\
&\hskip20pt
\leq\sum_{k\geq m}\pr\Big(G\text{ contains a self-avoiding path }P\text{ with }|P|=k\text{ and }\sum_{e\in P}\Gamma_1^{(e)}\leq ck\Big)
\\
&\hskip40pt
\leq\sum_{k\geq m}|V|r^k\exp\Big(-k\big(t-c(e^t-1)\big)\Big)\\
&\hskip60pt
=|V|\sum_{k\geq m}\exp\Big(-k\big(t-\log r-c(e^t-1)\big)\Big)
\leq |V|e^{-m},
\end{align*}
where the second inequality uses \eqref{eqn:4884} and the fact that there are at most $|V|r^k$ many self-avoiding paths of length $k$ in $G$, and the last step follows if we choose $t$ sufficiently large and $c$ sufficiently small.
This completes the proof.

\subsection{Proof of Proposition~\ref{prop:gh-lemma-1}}\label{sec:proof-gh-lemma-1}
Recall Construction~\ref{constr:S-k-alternate}.
Let $s, n, r, \cK_{n,3}$, $\cT_1,\ldots,\cT_r$, $\cT_1',\ldots,\cT_r'$, $\rho_i, z_i$, and $(X_1,\ldots,X_r)$ be as in Construction~\ref{constr:S-k-alternate}.
Using \eqref{eqn:connected-whp}, we can assume that $\cK_{n,3}$ and $\cmnthree$ are coupled in a way so that
\begin{align}\label{eqn:39}
\lim_{n\to\infty}\pr\big(\cK_{n,3}\neq\cmnthree\big)=0.
\end{align}

Let $\big\{\Gamma_{1/2}^{\sss(j)}\big\}_{1\leq j\leq r}$ be a sequence of i.i.d. $\mathrm{Gamma}(1/2,1)$ random variables.
Then
\begin{align}\label{eqn:234}
\big(X_1,\ldots,X_r\big)
\equald
\big(\Gamma_{1/2}^{\sss(1)},\ldots, \Gamma_{1/2}^{\sss(r)}\big)/\Gamma_{r/2}\ ,
\end{align}
where $\Gamma_{r/2}=\sum_{j=1}^r \Gamma_{1/2}^{\sss(j)}$.
Note also that
\begin{align}\label{eqn:40}
r^{-1}\Gamma_{r/2}\weakc 1/2,\ \text{ as }\ r\to\infty.
\end{align}

Let $Y_i:=d_{\cT_i}(\rho_i, z_i)$.
Then $Y_i$, $1\leq i\leq r$, are i.i.d. Rayleigh random variables \cite{aldous-crt-1, aldous-crt-3} with density $f(y)=y\exp(-y^2/2)$, $y>0$.
Hence $Y_i^2\equald 2\Gamma_1$, where $\Gamma_1\sim\mathrm{Exponential}(1)$.
By \cite{malik1968exact}, for $i=1,\ldots, r$, $Z_i^2:=Y_i^2 \Gamma_{1/2}^{(i)}\equald \Gamma_1^2/2$.
Hence
\begin{align}\label{eqn:41}
\big\{\sqrt{2}Z_i\big\}_{1\leq i\leq r}\ \text{ is an i.i.d. sequence of Exponential}(1)\ \text{ random variables}.
\end{align}

Let $\fQ_{n,3}=\big(k(\cH^{(s)}), e(\cH^{(s)}), \len\big)$ be the multigraph with edge lengths that represents $\core\big(\cH^{(s)}\big)$.
As observed right below Construction~\ref{constr:S-k-alternate}, $\fQ_{n,3}$ can be constructed by assigning length $Y_i\sqrt{X_i}=Z_i/\sqrt{\Gamma_{r/2}}$ to the $i$-th edge of $\cK_{n,3}$, $1\leq i\leq r$.
There is a natural coupling between $\cb^{\infty}\big(\fQ_{n,3}\big)$ and $\cb^{\infty}\big(\cH^{(s)}\big)$ in which $\cb^{\infty}\big(\cH^{(s)}\big)$ can be obtained by attaching countably many real trees to  $\cb^{\infty}\big(\fQ_{n,3}\big)$, and the diameter of each such real tree is at most $\max_{1\leq i\leq r}\diam\big(\cT_i'\big)=\max_{1\leq i\leq r}\sqrt{X_i}\cdot\diam\big(\cT_i\big)$.
Thus, in this coupling,
\begin{align}\label{eqn:42}
\big(12 s\big)^{1/6}\cdot d_H\Big(\cb^{\infty}\big(\fQ_{n,3}\big),\ \cb^{\infty}\big(\cH^{(s)}\big)\Big)
\leq
\big(12 s\big)^{1/6}\max_{1\leq i\leq r}\sqrt{X_i}\cdot\diam\big(\cT_i\big)\weakc 0,
\end{align}
where the last step is a consequence of the facts 
$\max_{1\leq i\leq r} X_i=O_P\big(\log r/r\big)$ (which can be seen from \eqref{eqn:234} and \eqref{eqn:40}), and
$\max_{1\leq i\leq r}\diam\big(\cT_i\big)=O_P\big(\sqrt{\log r}\big)$.

Now, in the coupling used in \eqref{eqn:cbd-cb-close},
\begin{align}\label{eqn:444}
\big(12 s\big)^{1/6}\cdot d_H\Big(\cbd_{\infty}\big(\fQ_{n,3}\big),~ \cb^{\infty}\big(\fQ_{n,3}\big)\Big)
\leq 
\big(12 s\big)^{1/6}\cdot\max_{1\leq i\leq r} Z_i/\sqrt{\Gamma_{r/2}}
\weakc 0,
\end{align}
where the last step follows from \eqref{eqn:41} and \eqref{eqn:40}.
Combining \eqref{eqn:444} with \eqref{eqn:42} and Theorem~\ref{thm:derandomization}, we see that as $n\to\infty$,
\begin{align}\label{eqn:43}
\big(12 s\big)^{1/6}\cbd_{\infty}\big(\fQ_{n,3}\big)\weakc\cM\ \ \text{ w.r.t. GH topology.}
\end{align}
Finally, using \eqref{eqn:40} and the relation $r=3(s-1)$, we see that the length of the $i$-th edge in $\big(12 s\big)^{1/6}\fQ_{n,3}$ is
\[\big(12 s\big)^{1/6}\frac{Z_i}{\sqrt{\Gamma_{r/2}}}=\big(1+o_P(1)\big)\cdot\left(\frac{4}{3n}\right)^{1/3}\sqrt{2}\cdot Z_i,\]
which together with \eqref{eqn:43} implies that as $n\to\infty$,
\begin{align}\label{eqn:44}
\left(\frac{4}{3n}\right)^{1/3}\cdot\cbd_{\infty}\big(\cK_{n,3}^{\exp}\big)\weakc\cM\ \ \text{ w.r.t. GH topology,}
\end{align}
where $\cK_{n,3}^{\exp}$ is the multigraph obtained by assigning lengths $\sqrt{2}Z_i$, $1\leq i\leq r$ (which, by \eqref{eqn:41}, are i.i.d. Exponential$(1)$ random variables) lengths to the edges of $\cK_{n,3}$.
We complete the proof by combining \eqref{eqn:44} and \eqref{eqn:39}.\qed

\subsection{Proof of Proposition~\ref{prop:gh-lemma-2}}\label{sec:proof-gh-lemma-2}
As in Section~\ref{sec:outline-of-proof}, for any graph $H$ and $p\in [0,1]$, we denote by $\perc(H, p)$ the random subgraph of $H$ obtained by removing edges of $H$ independently with probability $1-p$.
The proof of Proposition~\ref{prop:gh-lemma-2} relies on the following three lemmas.

\begin{lem}\label{lem:removed-or-marked}
Let $H$ be a finite multigraph.
Let $H^{\exp}$ be the multigraph with edge lengths obtained by assigning i.i.d. $\mathrm{Exponential}(1)$ lengths to the edges of $H$.
Fix $t>0$, and let $R(t)$ be a $\mathrm{Poisson}\big(t\cdot \ell(H^{\exp})\big)$ random variable.
Then
\[
\bigg( 
\shape\Big[\remove\big(\cbd_{R(t)}(H^{\exp})\big)\Big],\
\remove\Big(\cbd_{R(t)}(H^{\exp})\Big)
\bigg)
\]
has the same distribution as
\[
\bigg( 
\perc\Big(H, \frac{1}{1+t}\Big),\
\frac{1}{1+t}\cdot\Big(\perc\Big(H, \frac{1}{1+t}\Big)\Big)^{\exp}~
\bigg)\ ,
\]
where the last graph denotes the multigraph obtained by assigning i.i.d. $\mathrm{Exponential}(1)$ lengths to the edges of $~\perc(H, 1/(1+t))$, and then multiplying the lengths by $1/(1+t)$ (or equivalently, assigning i.i.d.exponential lengths with mean $1/(1+t)$ to the edges of $~\perc(H, 1/(1+t))$).
\end{lem}

Next we state two results about the behavior of the configuration model when a uniform subset of edges of given size is removed.
Recall the notation $\cmnd$ from Definition~\ref{def:configuration-model}.

\begin{lem}\label{lem:cm-coupling}
	Suppose $\vd=(d_1,\ldots,d_n)$ is a degree sequence and $\ell_n=\sum_{v\in[n]}d_i$.
	Let $m\leq \ell_n/2$ and define $\ell_n'=\ell_n-2m$.
	Out of the $\ell_n/2$ edges of $\cmnd$, sample a subset of $m$ edges uniformly.
	Let $\cmndm$ be the graph obtained by removing those $m$ edges.
	Then
	\begin{align}\label{eqn:101}
	\big(\cmnd, \cmndm\big)\equald\big(\cQ_{n,\vd,m}^{\sss (2)},\cQ_{n,\vd,m}^{\sss (1)}\big)~,
	\end{align}
	where the pair $\big(\cQ_{n,\vd,m}^{\sss (2)},\cQ_{n,\vd,m}^{\sss (1)}\big)$ is constructed as follows:
	Start with the vertex set $[n]$ with $d_i$ many half-edges attached to vertex $i$.
	Sample $\ell_n'$ many half-edges uniformly, and construct $\cQ_{n,\vd,m}^{\sss (1)}$ by uniformly pairing up those $\ell_n'$ half-edges.
	Conditional on this step, uniformly pair the rest of the half-edges to form $\cQ_{n,\vd,m}^{\sss (2)}$.
	
	Consequently, if $p\in[0,1]$ and $m$ is a $\mathrm{Binomial}(\ell_n/2, 1-p)$ random variable independent of $\cmnd$, then 
	\begin{align}\label{eqn:333}
	\cQ_{n,\vd,m}^{\sss (1)}\equald\perc\big(\cmnd, p\big).
	\end{align}
\end{lem}
Equality in both \ch{coordinates} in \eqref{eqn:101} will be used later in Section~\ref{sec:ghp-convergence-cm}.
In the proof of Theorem~\ref{thm:gh-mst-cm} we will only need \eqref{eqn:333}, which is a consequence of $\cmndm\equald\cQ_{n,\vd,m}^{\sss (1)}$, i.e., just the equality of the second coordinates in \eqref{eqn:101}. 
The relation $\cmndm\equald\cQ_{n,\vd,m}^{\sss (1)}$ was already observed in \cite[Lemmas 3.1 and 3.2]{fountoulakis2007percolation}.
See also \cite{janson2009percolation} for a related construction.

\begin{lem}\label{lem:edge-removal-criticality}
Suppose $m=m(n)$ satisfies
\begin{align}\label{eqn:868}
n^{-2/3}\big(3n-4m\big)\to\lambda_0\ \ \text{ as }\ \ n\to\infty,
\end{align}
for some $\lambda_0\in\bR$.
Let $\cmnthreem$ be the graph obtained by removing a uniform subset of $m$ edges from $\cmnthree$.
Let $\vd'=\vd'(n):=(d_1',\ldots,d_n')$, where $d_v'$ is the (random) degree of $v$ in $\cmnthreem$.
Let $\nu$ be the $\mathrm{Binomial}(3, 1/2)$ distribution.
Then as $n\to\infty$,
\begin{gather}
\nu_n'(i):=
\frac{1}{n}\#\big\{v\in[n]\, :\, d_v'=i\big\}\weakc\nu(i),\ \ i=0,1,2,3,\label{eqn:4383}\\
\frac{1}{n}\sum_{v\in [n]}(d_v')^3\weakc\sum_{i=0}^3 i^3\nu(i),\ \ \text{ and }\ \
n^{1/3}\left(\frac{\sum_{v\in [n]} (d_v')^2}{\sum_{v\in [n]} d_v'}-2\right)\weakc \frac{\lambda_0}{3}.\label{eqn:4384}
\end{gather}
\end{lem}

\begin{rem}
\ch{
Loosely speaking, Lemma~\ref{lem:edge-removal-criticality} says that the sequence $\big\{\vd'(n)\big\}_{n=2, 4,\ldots}$ of random degree sequences satisfies Assumption~\ref{ass:cm-deg} in probability.
This in particular implies that results for configuration models that only require Assumption~\ref{ass:cm-deg} (e.g., Theorem~\ref{thm:scaling-given-degree-sequence}) apply directly to the random graphs $\cmnthreem$.
This can be argued as follows:
By Skorohod representation theorem, we can construct 
$\mvnu_n'':=\big(\nu_n''(i)\, ;\, i=0, 1, 2, 3\big)$, $n=2, 4, 6, \ldots$, on the same probability space such that $\mvnu_n''\equald\big(\nu_n'(i)\, ;\, i=0, 1, 2, 3\big)$ for $n=2, 4, 6,\ldots$, and that 
\[
\nu_n''(i)\convas\nu(i)\, ,\ i=0, 1, 2, 3,\ \ \text{ and }\ \
n^{1/3}\left(\frac{\sum_{i=0}^{3}i^2\nu_n''(i)}{\sum_{i=0}^{3}i\nu_n''(i)}-2\right)
\convas 
\frac{\lambda_0}{3}
\]
in this space.
We can further assume that $\pi_2, \pi_4, \pi_6, \ldots$ are also defined on this space, where 
(a) $\pi_n$ is a uniform permutation of $n$ elements for $n=2, 4, \ldots$;
(b) $\pi_2, \pi_4, \pi_6, \ldots$ are independent; and 
(c) $\big(\pi_n;\, n=2, 4, \ldots\big)$ is independent of $\big(\mvnu_n'';\, n=2, 4,\ldots\big)$.
For $n=2, 4, 6, \ldots$, let $\vd''(n)$ be the random sequence of length $n$ obtained by applying $\pi_n$ to the sequence 
$\big(0, \ldots, 0, 1,\ldots, 1, 2,\ldots, 2, 3,\ldots, 3\big)$ with $i$ appearing $n\nu_n''(i)$ many times, $i=0, 1, 2, 3$.
Then $\vd''(n)\equald\vd'(n)$ for $n=2, 4, 6,\ldots$, and further, in this space, the convergences in \eqref{eqn:4383} and \eqref{eqn:4384} hold almost surely if $\vd'(n)$ is replaced by $\vd''(n)$.
Conditional on $\big(\vd''(n)\, ;\, n=2, 4,\ldots\big)$, construct $H_2, H_4,\ldots$, where $H_n$ is distributed as a configuration model with degree sequence $\vd''(n)$.
(The exact way in which $H_2, H_4, \ldots$ are coupled is not important here. 
For definiteness, let us take them to be independent conditional on $\big(\vd''(n)\, ;\, n=2, 4,\ldots\big)$.)
Then Theorem~\ref{thm:scaling-given-degree-sequence} applies to the sequence of random graphs $H_n$.
Now, from the equality in the second coordinate in \eqref{eqn:101}, we see that $\cmnthreem$, conditional on $\vd'(n)$, is distributed as a configuration model with degree sequence $\vd'(n)$.
Hence, $\cmnthreem\equald H_n$, and consequently,  
Theorem~\ref{thm:scaling-given-degree-sequence} applies to the random graphs $\cmnthreem$.
}	
\end{rem}

We now prove Proposition~\ref{prop:gh-lemma-2} assuming the above three lemmas.

\medskip

\noindent{\bf Completing the proof of Proposition~\ref{prop:gh-lemma-2}:}
We first note that if $\nu$ is the $\mathrm{Binomial}(3, 1/2)$ distribution, then
\begin{align}\label{eqn:99}
\sigma_1(\nu)=3/2,\ \ \sigma_2(\nu)=3,\ \text{ and }\ \sigma_3(\nu)=27/4.
\end{align}
Next, by Lemma~\ref{lem:removed-or-marked},
\begin{align}\label{eqn:1000}
\shape\big[\remove\big(\cbd_{R_{n,\lambda}}\big(\cmnthreeexp\big)\big)\big]
\equald
\perc\big(\cmnthree, 1/2+\lambda n^{-1/3}\big).
\end{align}

Now, for any $p\in[0,1]$, the number of edges removed from $\cmnthree$ to construct $\perc(\cmnthree, p)$ is a $\mathrm{Binomial}(3n/2, 1-p)$ random variable.
In particular, when $p=1/2+\lambda n^{-1/3}$, the number of edges removed is 
\begin{align}\label{eqn:33}
m=\frac{3n}{2}\Big(\frac{1}{2}-\frac{\lambda}{n^{1/3}}\Big)+O_P(\sqrt{n}),
\end{align}
which satisfies \eqref{eqn:868} with $\lambda_0=6\lambda$.
Further, conditional on $m$, $\perc(\cmnthree, 1/2+\lambda n^{-1/3})$ is distributed as $\cmnthree^{\sss (m)}$, where the notation is as in Lemma~\ref{lem:edge-removal-criticality}.
Thus, by Lemma~\ref{lem:edge-removal-criticality}, the (random) degree sequence of $\perc(\cmnthree, 1/2+\lambda n^{-1/3})$ satisfies $\eqref{eqn:4383}$ and \eqref{eqn:4384} with limiting parameter $\lambda_0/3=2\lambda$.
Finally, using \eqref{eqn:333}, it follows that $\perc(\cmnthree, 1/2+\lambda n^{-1/3})$, conditional on its degree sequence, is distributed as a configuration model with that degree sequence.
Hence by Theorem~\ref{thm:cm-component-sizes}, Theorem~\ref{thm:scaling-given-degree-sequence},  \eqref{eqn:99}, and \eqref{eqn:1000},
\begin{gather}
\cV=\Theta_P(n^{2/3}),\ \cE=\Theta_P(n^{2/3}),\ \cS=O_P(1),\ \cD=\Theta_P(n^{1/3}),\ \text{ and} \label{eqn:48}\\
n^{-1/3}\cdot G_1(n,\lambda)\weakc 6^{1/3}\cdot S_1\big((48)^{1/3}\cdot\lambda\big)\ 
\text{ w.r.t. GH topology}, \label{eqn:49}
\end{gather}
where $\cV=|G_1(n,\lambda)|$, $\cE=|E(G_1(n,\lambda))|$, $\cS=\mathrm{sp}(G_1(n,\lambda))$, and $\cD$ denotes the diameter of $G_1(n,\lambda)$.
By Lemma~\ref{lem:removed-or-marked}, conditional on $G_1(n,\lambda)$, the lengths of the edges of $\fG_1(n,\lambda)$ are
\begin{align}\label{eqn:50}
\big(1/2+\lambda n^{-1/3}\big)\cdot\big(\Gamma_1^{\sss(1)},\ldots, \Gamma_1^{\sss(\cE)}\big),
\end{align}
where $\Gamma_1^{\sss(1)},\ldots, \Gamma_1^{\sss(\cE)}$ are i.i.d. Exponential$(1)$ random variables.

Now it is easy to see that for any two vertices in $G_1(n,\lambda)$, there are at most $2^{\cS}$ many self-avoiding paths connecting them, and the length of any such self-avoiding path is at most $6(\cS+1)\cD$.
For any such self-avoiding path $P$ and any $\eta>0$, by standard concentration inequalities,
\begin{align}\label{eqn:51}
\pr_{G_1}\big(\big|\sum_{j\in P}\Gamma_1^{\sss (j)}-1\big|\geq \big(6(\cS+1)\cD\big)^{1/2+\eta}\big)
\leq
\exp\Big(-c\big(6(\cS+1)\cD\big)^{2\eta}\Big)\, ,
\end{align}
where $\pr_{G_1}$ denotes probability conditional on $G_1(n,\lambda)$.
Let $G_1^{\exp}(n,\lambda)$ be the graph with edge lengths obtained by assigning lengths $\Gamma_1^{\sss(1)},\ldots, \Gamma_1^{\sss(\cE)}$ to the edges of $G_1(n,\lambda)$.
Then by \eqref{eqn:51},
\[
\pr_{G_1}\Big(d_{\GH}\big(G_1^{\exp}(n,\lambda),\ G_1(n,\lambda)\big)\geq\big(6(\cS+1)\cD\big)^{1/2+\eta}\Big)
\leq
\cV^2\cdot 2^{\cS}\exp\Big(-c\big(6(\cS+1)\cD\big)^{2\eta}\Big).
\]
Thus, by \eqref{eqn:48} and \eqref{eqn:49}, $n^{-1/3}\cdot G_1^{\exp}(n, \lambda)\weakc 6^{1/3}\cdot S_1\big((48)^{1/3}\cdot\lambda\big)$ w.r.t. GH topology,
which together with \eqref{eqn:50} implies
\begin{align}\label{eqn:52}
n^{-1/3}\cdot \fG_1(n,\lambda)
\weakc
\frac{1}{2}\cdot 6^{1/3}\cdot S_1\big((48)^{1/3}\cdot\lambda\big)
\end{align}
w.r.t. GH topology.
The claim now follows from \eqref{eqn:49} and \eqref{eqn:52} by using Theorem~\ref{thm:mst-from-ghp-convergence} and Theorem~\ref{thm:belongs-to-A-r}. \qed

The rest of this section is devoted to the proofs of Lemmas~\ref{lem:removed-or-marked},~\ref{lem:cm-coupling}, and~\ref{lem:edge-removal-criticality}.

\medskip

\noindent{\bf Proof of Lemma~\ref{lem:removed-or-marked}:}
Let $|E(H)|=r$.
Run two independent Poisson point processes (PPP)--a `red' PPP and a `blue' PPP, with intensities $t$ and $1$ respectively.
Let $X_1<\ldots<X_r$ be the locations of the first $r$ blue points.
Enumerate the edges of $H$ in any way, and assign length $(X_i-X_{i-1})$ to the $i$-th edge, $i=1,\ldots,r$, where $X_0=0$.
Call the resulting graph with edge lengths $H_1$.
Let $\widetilde R$ be the number of red points in $[0, \ch{X_r}]$.
Identifying the $i$-th edge of $H_1$ with the interval $[X_{i-1}, X_i]$, $i=1,\ldots, r$, place a red point on $H_1$ corresponding to the location of each of the $\widetilde R$ red points in $[0,\ch{X_r}]$.
Call the resulting graph with red points $H_2$.

Now, note that $H_1\equald H^{\exp}$.
Next, conditional on the blue PPP, $\widetilde R$ follows a 
$\mathrm{Poisson}(t \ch{X_r})\equiv\mathrm{Poisson}(t\cdot\ell(H_1))$ distribution.
Thus,
\begin{align}\label{eqn:778}
\big(H_1,\widetilde R\big)
\equald
\big(H^{\exp}, R(t)\big).
\end{align}
Finally, conditional on the blue PPP and $\widetilde R$, the locations of the red points in $[0, \ch{X_r}]$ are i.i.d. $\mathrm{Uniform}[0, \ch{X_r}]$ random variables, which implies that 
\begin{align}\label{eqn:777}
\remove\big(\cbd_{R(t)}(H^{\exp})\big)\equald\remove(H_2).
\end{align}

Now $\remove(H_2)$ can be generated in the following alternate way: 
Sample independent random variables $Z_1,\ldots, Z_r$, where $Z_i\sim\mathrm{Poisson}(t(X_i-X_{i-1}))$, $i=1,\ldots, r$. (Here $Z_i$ corresponds to the number of red points in $[X_{i-1}, X_i]$.)
Remove the $i$-th edge of $H$ iff $Z_i\geq 1$, and assign independent lengths $Y_i$ to the remaining edges, where $Y_i\equald\big((X_i-X_{i-1}) \big| Z_i=0\big)$.

Combining \eqref{eqn:777} with the facts that $\ind_{Z_i=0}$, $i\geq 1$, are i.i.d. $\mathrm{Bernoulli}(1/(1+t))$ random variables, and
$\big((X_i-X_{i-1})\big| Z_i=0\big)$ has an exponential distribution with mean $1/(1+t)$, it follows that
\[
\remove\big(\cbd_{R(t)}(H^{\exp})\big)
\equald
\frac{1}{1+t}\cdot\Big(\perc\Big(H, \frac{1}{1+t}\Big)\Big)^{\exp}.
\]
Now the result follows immediately.
\qed

\medskip

\noindent{\bf Proof of Lemma~\ref{lem:cm-coupling}:}
Let $G$ be a graph on $[n]$ with degree sequence $\vd$.
Let $G'$ be a subgraph of $G$.
Let $\vd'=(d_1',\ldots,d_n')$ be the degree sequence of $G'$.
Let $x_{ij}$ (resp. $x_{ij}'$) be the number of edges between $i$ and $j$ in $G$ (resp. $G'$), $i\neq j$, and let $x_{ii}$ (resp. $x_{ii}'$)  denote the number of loops attached to vertex $i$ in $G$ (resp. $G'$).
Using \eqref{eqn:cm-distribution} it follows that
\begin{align}\label{eqn:25}
\pr\big(\cmnd=G,\ \cmndm=G'\big)
=
\frac{1}{(\ell_n-1)!!}
\times
\frac{\prod_{i\in [n]}d_i!}{\prod_{i\in[n]}2^{x_{ii}}\prod_{i\leq j}x_{ij}!}
\times
\frac{\prod_{i\leq j}\dbinom{x_{ij}}{x_{ij}'}}{\dbinom{\ell_n/2}{m}},
\end{align}
and
\begin{align}\label{eqn:26}
\pr\big(\cQ_{n,\vd,m}^{\sss (1)}=G',\   \cQ_{n,\vd,m}^{\sss (2)}=G\big)
&=
\frac{\prod_{i\in [n]}\dbinom{d_{i}}{d_{i}'}}{\dbinom{\ell_n}{\ell_n'}}
\times
\frac{1}{(\ell_n'-1)!!}
\times
\frac{\prod_{i\in [n]}d_i'!}{\prod_{i\in[n]}2^{x_{ii}'}\prod_{i\leq j}x_{ij}'!}\\
&\hskip45pt
\times
\frac{1}{(\ell_n-\ell_n'-1)!!}
\times
\frac{\prod_{i\in [n]}(d_i-d_i')!}{\prod_{i\in[n]}2^{x_{ii}-x_{ii}'}\prod_{i\leq j}(x_{ij}-x_{ij}')!}.\nonumber
\end{align}
A direct computation shows that the right sides of \eqref{eqn:25} and \eqref{eqn:26} are equal.
This completes the proof.
\qed

\medskip

\noindent{\bf Proof of Lemma~\ref{lem:edge-removal-criticality}:}
We use the alternate construction of $\cmnthreem$  from Lemma~\ref{lem:cm-coupling}.
For each $v\in [n]$, let $f_{v,i}$ denote the $i$-th half edge attached to $v$, $i=1, 2, 3$.
Let $E_{v, i}$ denote the event that $f_{v,i}$ is one of the $3n-2m$ selected half edges.
Then
\begin{align}
\pr\big(E_{v,i}\big)=(3n-2m)/3n,\  \text{ for }\ 1\leq i\leq 3.\label{eqn:27}
\end{align}
and
\begin{align}
\pr\big(E_{v_1, i_1}\cap E_{v_2, i_2}\big)=\frac{(3n-2m)(3n-2m-1)}{3n(3n-1)},\  \text{ whenever }\ (v_1, i_1)\neq (v_2, i_2).\label{eqn:28}
\end{align}
Since $d_v'=\sum_{i=1}^3\ind\big\{E_{v,i}\big\}$,
\[
\bE\bigg[\sum_{v\in [n]}d_v'^2\big]
=
n\cdot\bE\big[d_1'^2\big]=n\big[3\times \frac{(3n-2m)}{3n}
+
6\times\frac{(3n-2m)(3n-2m-1)}{3n(3n-1)}\bigg].
\]
Using this relation, \eqref{eqn:868}, and the fact that $\sum_{v\in [n]}d_v'=3n-2m$, a direct computation shows that
\begin{align}\label{eqn:29}
\lim_{n\to\infty}n^{1/3}\left(\frac{ \bE\big[\sum_{v\in [n]}d_v'^2\big]}{\sum_{v\in [n]} d_v'}-2\right)=\frac{\lambda_0}{3}~.
\end{align}
Now it is straightforward to check that for any four distinct pairs $(v_j, i_j)$, $1\leq j\leq 4$, each of the quantities
$\cov\big(\ind\{E_{v_1, i_1}\}, \ind\{E_{v_2, i_2}\}\big)$,
$\cov\big(\ind\{E_{v_1, i_1}\cap E_{v_2, i_2}\}, \ind\{E_{v_3, i_3}\}\big)$, and
$\cov\big(\ind\{E_{v_1, i_1}\cap E_{v_2, i_2}\}, \ind\{E_{v_3, i_3}\cap E_{v_4, i_4}\}\big)$
is negative.
Thus, for any $v_1\neq v_2$, $\cov\big(d_{v_1}'^2, d_{v_2}'^2\big)<0$, which implies that
\[\var\big(\sum_{v\in [n]}d_v'^2\big)\leq\sum_{v\in[n]}\var\big(d_v'^2\big)=O(n).\]
This combined with \eqref{eqn:29} proves the second convergence in \eqref{eqn:4384}.

Next, for $v\in[n]$ and $k=0, 1, 2, 3$,
\[
\pr\big(d_v'=k\big)=\dbinom{3}{k}\dbinom{3n-3}{2m+k-3}\Big/\dbinom{3n}{2m},
\]
which together with \eqref{eqn:868} yields
\[
\lim_{n\to\infty}\ \frac{1}{n}\bE\big[\#\big\{v\in[n] : d_v'=k\big\}\big]
=\lim_{n\to\infty}\pr\big(d_1'=k\big)
=\nu(k).
\]
A little computation will show that $\var\big[\#\big\{v\in[n] : d_v'=k\big\}\big]=O(n)$ for $k=0,1,2,3$.
This proves \eqref{eqn:4383}.
Finally, the first convergence in \eqref{eqn:4384} follows from \eqref{eqn:4383}.
This completes the proof.
\qed

\subsection{Proof of Proposition~\ref{prop:gh-lemma-3}}\label{sec:proof-gh-lemma-3}
We will use the following lemmas in the proof:
\begin{lem}\label{lem:8888}
Suppose $Y_1$ and $Y_2$ are two real valued random variables defined on the same probability space such that $Y_1\leq Y_2$ almost surely. 
Suppose further that $Y_1\equald Y_2$. 
Then $Y_1=Y_2$ almost surely.
\end{lem}
This is an elementary lemma, and we omit the proof.
\begin{lem}\label{lem:9999}
Suppose $\big\{(Z_n^+, d_n, Z_n)\big\}_{n\geq 1}$ is a sequence in $\fS_{\GH}^{\ast}$ satisfying
$(Z_n^+, d_n, Z_n)\to (Z_0^+, d, Z_0)$
for some marked space $(Z_0^+, d, Z_0)$.
Then 
\[
d_H\big(Z_n^+, Z_n\big)\to d_H(Z_0^+, Z_0).
\]
\end{lem}
\noindent{\bf Proof:}
For any isometric embeddings $\phi_n:Z_n^+\to Z^{\star}$ and $\psi_n:Z_0^+\to Z^{\star}$ into some common space $Z^{\star}$, we have
\[
d_H(Z_n^+, Z_n)
\leq 
d_H\big(\phi_n(Z_n^+), \psi_n(Z_0^+)\big)
+d_H(Z_0^+,  Z_0)
+d_H\big(\psi_n(Z_0), \phi_n(Z_n)\big).
\]
Using symmetry, we see that
\begin{align}\label{eqn:7777}
\big|d_H(Z_n^+, Z_n)-d_H(Z_0^+,  Z_0)\big|
\leq 
d_H\big(\phi_n(Z_n^+), \psi_n(Z_0^+)\big)
+d_H\big(\psi_n(Z_0), \phi_n(Z_n)\big).
\end{align}
Using the fact $(Z_n^+, d_n, Z_n)\to (Z_0^+, d, Z_0)$, we can choose $\phi_n, \psi_n$ in a way so that the right side of \eqref{eqn:7777} goes to zero as $n\to\infty$.
\qed

\medskip

We will now complete the proof of Proposition~\ref{prop:gh-lemma-3}.
For any compact metric space $(X,d)$ and $\delta>0$, let $N_{\delta}(X)$ be the minimum number of closed $\delta$ balls needed to cover $X$.

Since $X_n^+\weakc Z$, the sequence $\big\{(X_n^+, d_n)\big\}_{n\geq 1}$ is relatively compact w.r.t. GH topology.
Using Lemma~\ref{lem:s-seq-topology}(b), the sequence $\big\{(X_n^+, d_n, X_n)\big\}_{n\geq 1}$ is relatively compact w.r.t. the marked topology.
Thus, there exists a subsequence $\big\{n_k\big\}_{k\geq 1}$ and a random marked space $(Z_0^+, d, Z_0)$ such that
\begin{align}\label{eqn:coupling}
\big( X_{n_k}^+, d_{n_k}, X_{n_k}\big)\weakc \big(Z_0^+, d, Z_0\big)
\end{align}
as $k\to\infty$ with respect to the marked topology.
Since $X_n^+\weakc Z$ and $X_n\weakc Z$,
we must have $Z_0^+\equald Z\equald Z_0$ as compact metric spaces.
In particular, for all $\eps>0$, 
\begin{align}\label{eqn:100}
N_{\eps}(Z_0^+)\equald N_{\eps}(Z_0).
\end{align}
Since $Z_0$ is a closed subset of $Z_0^+$, for every $\eps>0$,
$N_{\eps}\big(Z_0\big)\leq N_{\eps}\big(Z_0^+\big)$ almost surely. 
Then it follows from \eqref{eqn:100} and Lemma~\ref{lem:8888} that 
\begin{align*}
\pr\big(N_{\eps}(Z_0^+)= N_{\eps}(Z_0)\big)=1
\end{align*}
for every $\eps>0$.
This implies that 
$\pr\big(d_H(Z_0^+, Z_0)=0\big)=1$.
Thus, using Lemma~\ref{lem:9999}, we conclude that $d_H(X_{n_k}^+, X_{n_k})\probc 0$.

Now for any subsequence $\big\{m_{\ell}\big\}_{\ell\geq 1}$, using the above argument, we can extract a further subsequence $\big\{m_{{\ell}_k}\big\}_{k\geq 1}$ such that 
$d_H\big(X_{m_{{\ell}_k}}^+, X_{m_{{\ell}_k}}\big)\probc 0$ as $k\to\infty$.
Thus the claim follows.

\subsection{GHP convergence of the MST of $\cmnthree$}\label{sec:ghp-convergence-cm}
In this section we improve the convergence in Theorem~\ref{thm:gh-mst-cm} to GHP convergence, thus completing the proof of Theorem~\ref{thm:mst_CM}.
Let $\fG_1(n,\lambda)$ and $G_1(n,\lambda)$ be as in the statement of Proposition~\ref{prop:gh-lemma-2}, and let $k_n(\lambda)$ and $\fT_{n,\lambda}^{(j)}$, $1\leq j\leq k_n(\lambda)$, be as in the proof of Theorem~\ref{thm:gh-mst-cm}.
For $v\in \fG_1(n,\lambda)$, let $d_{v,\lambda}$ be the degree of $v$ in $\fG_1(n,\lambda)$, and define $d_{v,\lambda}^{\avail}:=3-d_{v,\lambda}$.
Thus, $d_{v,\lambda}^{\avail}$ denotes the number of distinct edges sampled in the process $\big(\cbd_{i}(\cmnthreeexp), 1\leq i\leq R_{n, \lambda}\big)$ that were incident to $v$, and one can picture this degree deficiency as `available half-edges' attached to $v$.

On the event $B_n$, $\cbd_{\infty}(\cmnthree)$ is $\cbd_{\infty}\big(G_1(n,\lambda)\big)$ with $\shape\big[\fT_{n,\lambda}^{(j)}\big]$, $1\leq j\leq k_n(\lambda)$, attached to its vertices via a single edge;
let $T_{n,\lambda}^{(i)}(v)$, $1\leq i\leq r_{v,\lambda}$, be the trees (arranged following some deterministic rule) attached to $v\in\cbd_{\infty}\big(G_1(n,\lambda)\big)$.
Clearly, $0\leq r_{v,\lambda}\leq d_{v,\lambda}^{\avail}$.
Thus, the collection of trees $T_{n,\lambda}^{(i)}(v)$, $1\leq i\leq r_{v,\lambda}$, $v\in\cbd_{\infty}\big(G_1(n,\lambda)\big)$, is simply $\shape\big[\fT_{n,\lambda}^{(j)}\big]$, $1\leq j\leq k_n(\lambda)$, in some order.
Recall from the proof of Theorem~\ref{thm:gh-mst-cm} that we define $k_n(\lambda)=0$ for all $\lambda\in\bR$ on $B_n^c$. Accordingly, we set $r_{v,\lambda}=0$ for all $v\in\cbd_{\infty}\big(G_1(n,\lambda)\big)$ on the event $B_n^c$.

Construct the spaces $\fM_{n, \lambda}^{\attach}$ and $\fM_{n, \lambda}^{\avail}$ by endowing $\cbd_{\infty}\big(G_1(n,\lambda)\big)$ with the tree distance and respectively assigning mass
\begin{align}
p_{v,\lambda}^{\attach}:=
\left\{
\begin{array}{l}
1/|G_1(n,\lambda)|, 
\text{ on } B_n^c,\\
\\
\frac{1}{n}\big(1+\sum_{i=1}^{r_{v,\lambda}}\big| T_{n,\lambda}^{(i)}(v) \big|\big), 
\text{ on } B_n,
\end{array}
\right.
\ \ \text{ and } \ \
p_{v,\lambda}^{\avail}:=
\left\{
\begin{array}{l}
1/|G_1(n,\lambda)|, 
\text{ if } \sum_{u\in G_1(n,\lambda)} d_{u,\lambda}^{\avail}=0,\\
\\
d_{v,\lambda}^{\avail}\big/\big(\sum_{u\in G_1(n,\lambda)} d_{u,\lambda}^{\avail}\big),
\text{ otherwise,}
\end{array}
\right.\label{eqn:p-v}
\end{align}
to $v\in\cbd_{\infty}\big(G_1(n,\lambda)\big)$.
Note that $\sum_{v\in G_1(n,\lambda)}p_{v,\lambda}^{\attach}=
\sum_{v\in G_1(n,\lambda)}p_{v,\lambda}^{\avail}=1$.
Note also that the first and the third asymptotics in \eqref{eqn:48} imply that $\pr\big(\sum_{u\in G_1(n,\lambda)} d_{u,\lambda}^{\avail}=0\big)\to 0$ as $n\to\infty$. 
Thus, the value of $p_{v,\lambda}^{\avail}$ on the event
$\big\{\sum_{u\in G_1(n,\lambda)} d_{u,\lambda}^{\avail}\geq 1\big\}$ 
is the one relevant for distributional asymptotics of $\fM_{n, \lambda}^{\avail}$.
Similarly, using \eqref{eqn:connected-whp}, only the value of $p_{v,\lambda}^{\attach}$ on $B_n$ is relevant for the asymptotic behavior of $\fM_{n, \lambda}^{\attach}$.

\ch{Throughout Section~\ref{sec:ghp-convergence-cm}, all sequences $\{\lambda_n\}_{n\geq 1}$ will be $\bZ_{>0}$-valued sequences, and we will not mention this explicitly.}

\begin{lem}\label{lem:ghp-lemma-1}
Let $\lambda_n^{\star}$ be as in the proof of Theorem~\ref{thm:gh-mst-cm}. Then for all $\lambda_n\uparrow\infty$ with $\lambda_n\leq\lambda_n^{\star}$,
\[n^{-1/3}d_{\GHP}\big(\cbd_{\infty}\big(\cmnthree\big),\ \fM_{n, \lambda_n}^{\attach}\big)\probc 0.\]
\end{lem}
\begin{lem}\label{lem:ghp-lemma-2}
There exists a sequence $\lambda_n^{\dagger}\uparrow\infty$ such that for all $\lambda_n\uparrow\infty$ with $\lambda_n\leq\lambda_n^{\dagger}$,
\[n^{-1/3}\fM_{n, \lambda_n}^{\avail}\weakc 6^{1/3}\cdot\cM\ \ \text{ w.r.t. the GHP topology.}\]
\end{lem}
\begin{lem}\label{lem:ghp-lemma-3}
There exists a sequence $\lambda_n^{\circ}\uparrow\infty$ such that for all $\lambda_n\uparrow\infty$ with $\lambda_n\leq\lambda_n^{\circ}$,
\[n^{-1/3}d_{\GHP}\big(\fM_{n, \lambda_n}^{\attach},\ \fM_{n, \lambda_n}^{\avail}\big)\probc 0.\]
\end{lem}

\noindent{\bf Completing the proof of Theorem~\ref{thm:mst_CM}:}
The result follows upon combining Lemma~\ref{lem:ghp-lemma-1}, Lemma~\ref{lem:ghp-lemma-2}, Lemma~\ref{lem:ghp-lemma-3}, and \eqref{eqn:679}.
\qed

\medskip

\noindent{\bf Proof of Lemma~\ref{lem:ghp-lemma-1}:}
On the event $B_n$, define the correspondence $C$ between $\fM_{n, \lambda_n}^{\attach}$ and $\cbd_{\infty}\big(\cmnthree\big)$ as follows:
\[C:=\big\{(v, u)\ :\ v\in\fM_{n, \lambda_n}^{\attach}\ \text{ and }\ u\in\{v\}\cup\big(\bigcup_{i=1}^{r_{v,\lambda_n}} T_{n,\lambda_n}^{(i)}(v)\big)\big\}.\]
Let $\pi$ be a measure on $\fM_{n, \lambda_n}^{\attach}\times\cbd_{\infty}\big(\cmnthree\big)$ given by $\pi(\{(v,u)\})=1/n$ for $(v,u)\in C$.
Then with this choice of $C$ and $\pi$, the claim follows immediately if we use \eqref{eqn:37}.
\qed

\medskip

\noindent{\bf Proof of Lemma~\ref{lem:ghp-lemma-2}:}
Assign mass $p_{v,\lambda}^{\avail}$ to $v\in G_1(n,\lambda)$ and call the resulting metric measure space $G_1^{\avail}(n, \lambda)$.
Using Theorem~\ref{thm:scaling-given-degree-sequence} with $f(k)=3-k$, $k=0,\ldots,3$, and the arguments used to prove \eqref{eqn:49}, we see that
$n^{-1/3} G_1^{\avail}(n, \lambda)\weakc 6^{1/3}\cdot S_1\big((48)^{1/3}\cdot\lambda\big)$
w.r.t. GHP topology.
Using Theorem~\ref{thm:belongs-to-A-r} and Theorem~\ref{thm:mst-from-ghp-convergence}, it follows that for each $\lambda\in\bR$,
\begin{align}\label{eqn:60}
n^{-1/3}\fM_{n, \lambda}^{\avail}\weakc
6^{1/3}\cdot\cb^{\infty}\big(S_1\big((48)^{1/3}\cdot\lambda\big)\big)\ \ \text{ as }\ \ n\to\infty
\end{align}
w.r.t. GHP topology.
The claim now follows from Theorem~\ref{thm:M-alternate} and Lemma~\ref{lem:23}.
\qed

\medskip

To prove Lemma~\ref{lem:ghp-lemma-3} we will make use of Lemma~\ref{lem:averages-out} stated below. Let $\lambda_n^{\star}$ be as in the proof of Theorem~\ref{thm:gh-mst-cm}.

\begin{lem}\label{lem:averages-out}
There exists a sequence $\lambda_n^{\Diamond}\uparrow\infty$ such that

\noindent {\upshape (i)} $\lambda_n^{\Diamond}\leq\lambda_n^{\star}$, 

\noindent {\upshape (ii)}  
$
\pr\big(\big|G_1(n,\lambda_n^{\Diamond})\big|>n/2\big)\to 0
$,
and 

\noindent {\upshape (iii)} 
for any 
$\lambda_n\uparrow\infty$ with $\lambda_n\leq\lambda_n^{\Diamond}$, the following holds:
For every $n$, fix an enumeration $v_1, v_2,\ldots$ of the vertices of $G_1(n,\lambda_n)$ measurable w.r.t. the $\sigma$-field generated by
$
\cbd_{\infty}\big(G_1(n,\lambda_n)\big)
$,
and define
\[
Z_n:=
\ch{\max_{1\leq j_1\leq j_2\leq |G_1(n,\lambda_n)|}\,
\bigg|\sum_{s=j_1}^{j_2}
\sum_{i=1}^{r_{v_s,\lambda_n}}
\frac{|T_{n,\lambda_n}^{(i)}(v_s)|}{n-|G_1(n,\lambda_n)|}
-
\sum_{s=j_1}^{j_2} p_{v_s,\lambda_n}^{\avail}
\bigg|
~ ,
}
\]
where $ p_{v,\lambda}^{\avail}$ is as defined in \eqref{eqn:p-v}.
Then $Z_n\probc 0$.
\end{lem}

We first prove Lemma~\ref{lem:ghp-lemma-3} assuming Lemma~\ref{lem:averages-out}.

\noindent{\bf Proof of Lemma~\ref{lem:ghp-lemma-3}:}
On the event $B_n$,
construct $\fM_{n, \lambda}^{\modi}$ by endowing $\cbd_{\infty}\big(G_1(n,\lambda)\big)$ with the tree distance and assigning mass
\[p_{v,\lambda}^{\modi}:=
\frac{\sum_{i=1}^{r_{v,\lambda}}\big| T_{n,\lambda}^{(i)}(v) \big|}{n-\big|G_1(n,\lambda)\big|}\]
to $v\in\cbd_{\infty}\big(G_1(n,\lambda)\big)$.
On $B_n^c$, set $\fM_{n, \lambda}^{\modi}=\fM_{n, \lambda}^{\attach}$.
As observed in \eqref{eqn:48}, $\big|G_1(n,\lambda)\big|=\Theta_P(n^{2/3})$.
Thus,
\[\sum_{v\in G_1(n,\lambda)}\big|p_{v,\lambda}^{\modi}-p_{v,\lambda}^{\attach}\big|=O_P(n^{-1/3}).\]
It follows that for each $\lambda\in\bR$,
$d_{\GHP}\big(\fM_{n, \lambda}^{\modi}, \fM_{n, \lambda}^{\attach}\big)\probc 0$
as $n\to\infty$.
Thus, we can choose a sequence $\lambda_n^{\oplus}\uparrow\infty$ such that 
$\pr\big(\big|G_1(n,\lambda_n^{\oplus})\big|>n/2\big)\to 0$,
and further,
for all $\lambda_n\uparrow\infty$ with $\lambda_n\leq\lambda_n^{\oplus}$, 
\begin{align}\label{eqn:61}
d_{\GHP}\big(\fM_{n, \lambda_n}^{\modi},\ \fM_{n, \lambda_n}^{\attach}\big)\probc 0.
\end{align}
Set $\lambda_n^\circ:=\min\{\lambda_n^{\oplus}, \lambda_n^{\dagger}, \lambda_n^{\Diamond}\}$, where $\lambda_n^{\dagger}$ (resp. $\lambda_n^{\Diamond}$) is as in Lemma~\ref{lem:ghp-lemma-2} (resp. Lemma~\ref{lem:averages-out}).
Fix a sequence $\lambda_n\uparrow\infty$ with $\lambda_n\leq\lambda_n^{\circ}$.

Fix $\delta>0$.
Let 
$N_{\delta}^{\sss (n)}$
be the minimum number of closed $\delta n^{1/3}$ balls needed to cover $\cbd_{\infty}\big(G_1(n,\lambda_n)\big)$.
By Lemma~\ref{lem:ghp-lemma-2}, $\big\{N_{\delta}^{\sss (n)}\big\}_{n\geq 1}
$ is tight.
Write $\dV_n$ for the set of vertices of $\cbd_{\infty}\big(G_1(n,\lambda_n)\big)$ and $d_{\infty}$ for the tree distance in $\cbd_{\infty}\big(G_1(n,\lambda_n)\big)$.
Let $A_1,\ldots,A_{N_{\delta}^{(n)}}$ be a partition of $\dV_n$ such that for $1\leq j\leq N_{\delta}^{\sss (n)}$, $d_{\infty}(v,v')\leq 2\delta n^{1/3}$ if $v, v'\in A_j$.

Let $v_1,v_2,\ldots$ be an enumeration of $\dV_n$ such that for each $j\leq N_{\delta}^{\sss (n)}$, all vertices $v\in A_j$ appear successively.
Note that this enumeration is measurable w.r.t. the $\sigma$-field generated by $\cbd_{\infty}\big(G_1(n,\lambda_n)\big)$.
By Lemma~\ref{lem:averages-out},
\begin{align}\label{eqn:62}
\max_{1\leq j\leq N_{\delta}^{\sss (n)}}\Big|\ch{\sum_{v\in A_j}}\big(p_{v,\lambda_n}^{\modi}-p_{v,\lambda_n}^{\avail}\big)\Big|\probc 0.
\end{align}
Let $\mu_n^{\modi}$ be the measure on $\fX_n:=\{A_1,\ldots, A_{N_{\delta}^{\sss (n)}}\}$ given by
$\mu_n^{\modi}(A_j)=\sum_{v\in A_j} p_{v,\lambda_n}^{\modi}$.
Define $\mu_n^{\avail}$ on $\fX_n$ analogously.
Then the total variation distance between $\mu_n^{\modi}$ and $\mu_n^{\avail}$ satisfies
\[d_{\mathrm{TV}}\big(\mu_n^{\modi}, \mu_n^{\avail}\big)
\leq
\frac{1}{2}\times N_{\delta}^{\sss (n)}\times
\max_{1\leq j\leq N_{\delta}^{(n)}}\Big|\sum_{v\in A_j}\big(p_{v,\lambda_n}^{\modi}-p_{v,\lambda_n}^{\avail}\big)\Big|\probc 0,\]
where the last step uses \eqref{eqn:62} and the fact that $\big\{N_{\delta}^{\sss (n)}\big\}_{n\geq 1}$ is tight.
Thus, for each $n$, we can construct $\fX_n$-valued random variables $X_n^{\modi}$ and $X_n^{\avail}$ distributed as $\mu_n^{\modi}$ and $\mu_n^{\avail}$ respectively such that $\pr\big(X_n^{\modi}\neq X_n^{\avail}\big)\probc 0$.
Using $X_n^{\modi}$ and $X_n^{\avail}$, there is a natural way to construct $\dV_n$-valued random variables $Y_n^{\modi}$ and $Y_n^{\avail}$ such that
$\pr\big(Y_n^{\modi}=v\big)=p_{v,\lambda_n}^{\modi}$, and
$\pr\big(Y_n^{\avail}=v\big)=p_{v,\lambda_n}^{\avail}$ for all $v\in\dV_n$, and further,
\[\pr\big(d_{\infty}(Y_n^{\modi}, Y_n^{\avail})> 2\delta n^{1/3}\big)\probc 0.\]
Since $\delta>0$ was arbitrary, we get
$
n^{-1/3}d_{\GHP}\big(\fM_{n, \lambda_n}^{\modi},\ \fM_{n, \lambda_n}^{\avail}\big)\probc 0
$,
which combined with \eqref{eqn:61} completes the proof.
\qed

\medskip

The proof of Lemma~\ref{lem:averages-out} relies on the next two lemmas.
\begin{lem}\label{lem:concentration-uniform-permutation}
	There exist universal constants $c_1, c_2>0$ such that for any $m\geq 1$ and probability vector $\mvp:=(p_1,\ldots,p_m)$, 
	\begin{align}
	\bP\bigg(\max_{j\in[m]}\bigg|\sum_{i=1}^j p_{\pi(i)}-\frac{j}{m}\bigg|\geq x\sigma(\mvp)\bigg)
	\leq\exp\big(-c_1 x\log\log x\big), \ \text{ for }\ x\geq c_2~,
	\end{align}
	where $\pi$ is a uniform permutation on $[m]$, and $\sigma(\mvp):=\sqrt{p_1^2+\ldots+p_m^2}$.
	Consequently, using the relation $\sigma(\mvp)\leq\max_j \sqrt{p_j}$, we get, for $x\geq c_2$,
	\begin{align}\label{eqn:concentration} 
	\bP\Bigg(\max_{j_1< j_2}\
	\bigg|\sum_{i=j_1+1}^{j_2} p_{\pi(i)}-\frac{j_2-j_1}{m}\bigg|\geq 2x\cdot\max_j \sqrt{p_j}\Bigg)
	\leq 
	2\exp\big(-c_1 x\log\log x\big)~.
	\end{align}
\end{lem}
This result gives a quantitative concentration inequality for the partial sums of exchangeable random variables.
The result can be found in the above form in \cite[Lemma 7.5]{bhamidi-sen}, but was essentially already contained in \cite[Lemma 4.9]{bhamidi-hofstad-sen}.


\begin{lem}\label{lem:exchangeable}
\begin{inparaenumii}
\item\label{item:exchangeable} 
Fix $\lambda\in\bR$. 
For every $v\in G_1(n,\lambda)$, append $(d_{v,\lambda}^{\avail}-r_{v,\lambda})$ many zeros to the sequence $\big(\big|T_{n,\lambda}^{(i)}(v)\big|,~1\leq i\leq r_{v,\lambda}\big)$ and let $\big(\alpha_{n,\lambda}^{(i)}(v),~1\leq i\leq d_{v,\lambda}^{\avail}\big)$ be a uniform permutation of the resulting sequence; use independent permutations for different $v\in G_1(n,\lambda)$ that are also independent of all the other random variables being considered.
Then conditional on $\shape\big[\remove\big(\cbd_{R_{n,\lambda}}(\cmnthreeexp)\big)\big]$ and $\cbd_{\infty}\big(G_1(n,\lambda)\big)$,  the \ch{family 
\[
\bigg( \alpha_{n,\lambda}^{(i)}(v)\, ;\ 1\leq i\leq d_{v,\lambda}^{\avail},\ v\in G_1(n,\lambda) \bigg)
\]
of random variables} is exchangeable.

\noindent\item\label{item:max-mass-goes-to-zero} Let $\lambda_n^{\star}$ be as in the proof of Theorem~\ref{thm:gh-mst-cm}.
Then for any $\lambda_n\uparrow\infty$ with $\lambda_n\leq\lambda_n^{\star}$,
\begin{align}\label{eqn:55}
\max\ \Big\{
\frac{|T_{n,\lambda_n}^{(i)}(v)|}{n}\ :\ 
1\leq i\leq r_{v,\lambda_n},~ v\in\cbd_{\infty}\big(G_1(n,\lambda_n)\big)
\Big\}
\probc 0.
\end{align}
\end{inparaenumii}
\end{lem}

\noindent{\bf Proof of Lemma~\ref{lem:averages-out}:}
By \eqref{eqn:48}, $|G_1(n,\lambda)|=\Theta_P(n^{2/3})$.
So, in particular, for every $\lambda\in\bR$, $\pr\big(|G_1(n,\lambda)|>n/2\big)\to 0$.
Hence, we can choose $\lambda_n^{\Diamond}\uparrow\infty$ slowly enough such that 
$\pr\big(|G_1(n,\lambda_n^{\Diamond})|>n/2\big)\to 0$ as $n\to\infty$.
We can further take $\lambda_n^{\Diamond}\leq\lambda_n^{\star}$.

Fix $\lambda_n\uparrow\infty$ with $\lambda_n\leq\lambda_n^{\Diamond}$.
Let $v_1, v_2,\ldots$ be an enumeration of the vertices of $G_1(n,\lambda_n)$ measurable w.r.t. the $\sigma$-field generated by
$
\cbd_{\infty}\big(G_1(n,\lambda_n)\big)
$.
Define 
\[p_{n,\lambda_n}^{(i)}(v_s)
:=
\frac{\alpha_{n,\lambda_n}^{(i)}(v_s)}{n-|G_1(n,\lambda_n)|}\, ,\ \ 1\leq i\leq d_{v,\lambda_n}^{\avail}, \ 1\leq s\leq |G_1(n,\lambda_n)|\, ,
\] 
where $0/0$ is interpreted as $1$.
Since $\lambda_n\leq\lambda_n^{\Diamond}$, $\pr\big(n-|G_1(n,\lambda_n)|\geq n/2\big)\to 1$ as $n\to\infty$ by our choice of $\lambda_n^{\Diamond}$.
Thus, using Lemma~\ref{lem:exchangeable} (ii) and the fact that  $\lambda_n\leq\lambda_n^{\Diamond}\leq\lambda_n^{\star}$,
\begin{align}\label{eqn:A}
\max\big\{ p_{n,\lambda_n}^{(i)}(v_s)\ :\   
1\leq i\leq d_{v,\lambda_n}^{\avail}, \ 1\leq s\leq |G_1(n,\lambda_n)| \big\}
\probc 0.
\end{align}

Now, for any $1\leq j_1\leq j_2\leq |G_1(n, \lambda_n)|$,
\[
\sum_{s=j_1}^{j_2}\sum_{i=1}^{r_{v_s,\lambda_n}}|T_{n,\lambda_n}^{(i)}(v_s)|
=\sum_{s=j_1}^{j_2}\sum_{i=1}^{d_{v_s,\lambda_n}^{\avail}}\alpha_{n,\lambda_n}^{(i)}(v_s)~,
\]
and in particular, on the event $B_n$,
\[
\sum_{s=1}^{|G_1(n, \lambda_n)|}\
\sum_{i=1}^{d_{v_s,\lambda_n}^{\avail}}\alpha_{n,\lambda_n}^{(i)}(v_s)=n-|G_1(n,\lambda_n)|~.
\]
Thus, on the event $B_n\cap\{|G_1(n,\lambda_n)|<n\}$, 
$\big(p_{n,\lambda_n}^{(i)}(v_s)~;\ 1\leq i\leq d_{v,\lambda_n}^{\avail}, \ 1\leq s\leq |G_1(n,\lambda_n)|\big)$
is a probability vector.
By \eqref{eqn:connected-whp}, 
$\pr\big(B_n\cap\{|G_1(n,\lambda_n)|<n\}\big)\to 1$ as $n\to\infty$.
Thus, the desired result follows from Lemma~\ref{lem:exchangeable} (i), \eqref{eqn:concentration}, and \eqref{eqn:A}.
\qed

\medskip

\noindent{\bf Proof of Lemma~\ref{lem:exchangeable}\eqref{item:exchangeable}:}
Consider a finite (non-random) graph $H$ and $t>0$, and let $H^{\exp}$ and $R(t)$ be as in the statement of Lemma~\ref{lem:removed-or-marked}.
Then conditional on $\shape\big[\remove\big(\cbd_{R(t)}(H^{\exp})\big)\big]=H_0$, the order in which the edges in $E(H)\setminus E(H_0)$ were sampled for the first time in the $\cbd$ process is a uniform permutation on $E(H)\setminus E(H_0)$.

Using the above observation, Lemma~\ref{lem:removed-or-marked}, and Lemma~\ref{lem:cm-coupling}, 
we can generate $\shape\big[\remove\big(\cbd_{R_{n,\lambda}}(\cmnthreeexp)\big)\big]$, $\cmnthree$, and
$\big(\alpha_{n,\lambda}^{(i)}(v);\ 1\leq i\leq d_{v,\lambda}^{\avail},\ v\in G_1(n,\lambda)\big)$
jointly as follows:
\begin{enumeratea}
\item
Sample a $\mathrm{Binomial}\big(3n/2, 1/2-\lambda n^{-1/3}\big)$ random variable. For simplicity, we denote the realization by $m$.
\item
Conditional on step (a), sample $\cQ_{n,\mvthree,m}^{\sss (1)}$ as in Lemma~\ref{lem:cm-coupling}, where $\mvthree=(3,\ldots,3)$.
By \eqref{eqn:333},
$\cQ_{n,\mvthree,m}^{\sss (1)}
\equald
\perc\big(\cmnthree, 1/2+\lambda n^{-1/3}\big)$.
Hence, using Lemma~\ref{lem:removed-or-marked},
\begin{align}\label{eqn:233}
\cQ_{n,\mvthree,m}^{\sss (1)}
\equald
\shape\big[\remove\big(\cbd_{R_{n,\lambda}}(\cmnthreeexp)\big)\big]~.
\end{align}
Thus, the largest component of $\cQ_{n,\mvthree,m}^{\sss (1)}$, say $\cC_1$, has the same law as $G_1(n,\lambda)$.
Let $d_v$ be the degree of $v\in\cC_1$.
Then each $v\in\cC_1$ has $3-d_v$ many `available' half-edges; we denote them by $f_{v, i}$, $1\leq i\leq 3-d_v, v\in\cC_1$.
\item
Conditional on steps (a) and (b), generate $\cQ_{n,\mvthree,m}^{\sss (2)}$ as in Lemma~\ref{lem:cm-coupling}.
From \eqref{eqn:101} it follows that 
\begin{align}\label{eqn:233-a}
\cQ_{n,\mvthree,m}^{\sss (2)}\equald \cmnthree
\end{align}
jointly with the equality in distribution in \eqref{eqn:233}.
Let $\cE_1,\ldots,\cE_m$ be the edges that are in $\cQ_{n,\mvthree,m}^{\sss (2)}$ but not in $\cQ_{n,\mvthree,m}^{\sss (1)}$.
For $v\in\cC_1$ and $1\leq i\leq 3-d_v$, let $\widetilde f_{v, i}$ denote the edge that was formed by pairing $f_{v, i}$ with another half-edge.
\item
If $\cQ_{n,\mvthree,m}^{\sss (2)}$ is not connected, go to the next step.
If $\cQ_{n,\mvthree,m}^{\sss (2)}$ is connected, let $\pi$ be a uniform permutation of $m$ elements independent of steps (a), (b), and (c) above.
Consider the edges $\cE_{\pi(1)},\ldots,\cE_{\pi(m)}$ sequentially in this order, and at each step, remove the edge being considered from $\cQ_{n,\mvthree,m}^{\sss (2)}$ if its removal does not disconnect the current graph.
Denote the resulting graph by $\cQ$.
Then $\cQ$ has the same law as $\shape\big[\cbd_{R_{n,\lambda}}(\cmnthreeexp)\big]$.
\item
If $\cQ_{n,\mvthree,m}^{\sss (2)}$ is not connected, define $Q^{(i)}(v)$ to be the empty graph for $1\leq i\leq 3-d_v$, $v\in \cC_1$.
If $\cQ_{n,\mvthree,m}^{\sss (2)}$ is connected, then note that $\cQ$ as constructed in (d) is simply $\cC_1$ together with some connected multigraphs each of which is connected to a vertex of $\cC_1$ by a single edge;
for $v\in\cC_1$ and $1\leq i\leq 3-d_v$, set $Q^{(i)}(v)$ to be the connected multigraph that is connected to $v$ via $\widetilde f_{v, i}$, with the convention that $Q^{(i)}(v)$ is the empty graph if $\widetilde f_{v, i}$ was removed in step (d).
Then 
\[
\bigg(\text{sort}\big(|Q^{(i)}(v)|~,\ 1\leq i\leq 3-d_v\big)~;\ v\in \cC_1\bigg)
\equald
\bigg(\text{sort}\big(\alpha_{n,\lambda}^{(i)}(v)~,\ 1\leq i\leq d_{v,\lambda}^{\avail}\big)~;\ v\in G_1(n,\lambda)\bigg)
\]
jointly with the distributional equalities in \eqref{eqn:233} and \eqref{eqn:233-a}, where $\text{sort}(\cdot)$ arranges the entries of a finite sequence in decreasing order.
\end{enumeratea}

Conditional on steps (a) and (b) above, the rest of the procedure is symmetric with respect to the available half-edges attached to the vertices of $\cC_1$.
Hence, conditional on $\cQ_{n,\mvthree,m}^{\sss (1)}$, 
the family
$\big(|Q^{(i)}(v)|\, ;\ 1\leq i\leq 3-d_v,\ v\in \cC_1\big)$
is exchangeable.
Now, the conditional distribution of $\big(|Q^{(i)}(v)|;\ 1\leq i\leq 3-d_v,\ v\in \cC_1\big)$ given $\cQ_{n,\mvthree,m}^{\sss (1)}$ is equal to the conditional distribution of the same sequence given $\cQ_{n,\mvthree,m}^{\sss (1)}$ and $\cbd_{\infty}(\cC_1)$.
Thus, $\big(|Q^{(i)}(v)|\, ;\ 1\leq i\leq 3-d_v,\ v\in \cC_1\big)$ is an exchangeable family given 
$\cQ_{n,\mvthree,m}^{\sss (1)}$ and $\cbd_{\infty}(\cC_1)$.
Thus the claim follows.
\qed

\medskip

We need the following result before proceeding to the proof of Lemma~\ref{lem:exchangeable}\eqref{item:max-mass-goes-to-zero}.
Recall the notation $\fm(\cdot ; \cdot)$ from Section~\ref{sec:notation}.
\begin{lem}\label{lem:max-attachment-mass-small}
Suppose $(Z_n, d_n, \mu_n)\to (Z,d,\mu)$ as $n\to\infty$ in $\fS_{\GHP}$.
If $\mu$ is non-atomic, then
\[
\lim_{\eps\downarrow 0}\ \limsup_{n\to\infty}\ \fm(\eps, Z_n)=0.
\]
\end{lem}
\noindent{\bf Proof:}
Using the convergence $(Z_n, d_n, \mu_n)\to (Z,d,\mu)$, it is easy to see that for every $\eps>0$ and sufficiently large $n$,
\[
\fm(\eps, Z_n)\leq \fm(2\eps, Z)+\eps.
\]
It follows easily from the compactness of $(Z, d)$ and the non-atomicity of $\mu$ that
\[
\lim_{\eps\downarrow 0}\  \fm(2\eps, Z)=0.
\]
Thus the claim follows.
\qed

\medskip

We now continue with\\
\noindent{\bf Proof of Lemma~\ref{lem:exchangeable}\eqref{item:max-mass-goes-to-zero}:}
Recall the notation from Construction~\ref{constr:S-k-alternate}.
Let $\Gamma_{r/2}$ and $\cK_{n,3}^{\exp}$ be as in the proof of Proposition~\ref{prop:gh-lemma-1}.
As observed in the proof of Proposition~\ref{prop:gh-lemma-1},  there is a natural isometric embedding of $\cK_{n,3}^{\exp}$ into $\sqrt{2 \Gamma_{r/2}}\cdot \cH^{(s)}$.
In this embedding, $\sqrt{2 \Gamma_{r/2}}\cdot \cH^{(s)}$ can be obtained by attaching countably many real trees to $\cK_{n,3}^{\exp}$.
Let $(\cK_{n,3}^{\exp}, \mu_n)$ be the measured $\bR$-graph derived by endowing $\cK_{n,3}^{\exp}$ by the measure obtained by projecting the measure from $\sqrt{2 \Gamma_{r/2}}\cdot \cH^{(s)}$ onto the attachment points in $\cK_{n,3}^{\exp}$.
Thus, the $\mu_n$ measure of the $i$-th edge of $\cK_{n,3}^{\exp}$ is $X_i$, $1\leq i\leq r$, where $(X_1,\ldots,X_r)\sim\mathrm{Dirichlet}(1/2,\ldots,1/2)$ as in Construction~\ref{constr:S-k-alternate}.

Arguing as in \eqref{eqn:42}, it is easy to show that as $n\to\infty$,
\[n^{-1/3}d_{\GHP}\big(\cb^{\infty}\big(\sqrt{2 \Gamma_{r/2}}\cdot \cH^{(s)}\big),\ \cb^{\infty}\big(\cK_{n,3}^{\exp}\big)\big)\probc 0.\]
Combining this with Theorem~\ref{thm:derandomization} and \eqref{eqn:40}, we get
\begin{align}\label{eqn:56}
n^{-1/3}\cb^{\infty}\big(\cK_{n,3}^{\exp}\big)\weakc (0.75)^{1/3}\cdot\cM\ \ \text{ w.r.t. GHP topology.}
\end{align}
Suppose $\cmnthree$ and $\cK_{n,3}$ are coupled as in \eqref{eqn:39}.
On the event $\{\cK_{n,3}=\cmnthree\}$, the tree $\cb^{\infty}\big(\cK_{n,3}^{\exp}\big)$ can be obtained by
\begin{inparaenumi}
\item
attaching each of the trees $\fT_{n,\lambda_n}^{(j)}$, $1\leq j\leq k_n(\lambda_n)$, to $\cbd_{\infty}\big(\fG_1(n,\lambda_n)\big)$ via a single edge, and then
\item
attaching some additional line segments to the space thus obtained.
\end{inparaenumi}
(Recall that in the $\cb$ process edges are cut open, while in the $\cbd$ process edges are removed. 
Because of this difference these additional line segments need to be attached.)
Thus, using \eqref{eqn:56}, Theorem~\ref{thm:non-atomic}, \eqref{eqn:36}, together with Lemma~\ref{lem:max-attachment-mass-small}, we see that 
\begin{align}\label{eqn:189}
\max_{1\leq j\leq k_n(\lambda_n)}\ind_{\{\cmnthree=\cK_{n,3}\}}\times\mu_n\big(\fT_{n,\lambda_n}^{(j)}\big)\probc 0.
\end{align}

Consider the tree among $\fT_{n,\lambda_n}^{(j)}$, $1\leq j\leq k_n(\lambda_n)$, that has the maximum number of edges (pick any one if there is more that one such tree), and let $f_1,\ldots, f_{\cE_{\max}}$ be an enumeration of its edges.
On the event $\{\cmnthree=\cK_{n,3}\}$,
\begin{align}\label{eqn:58}
\max_{1\leq j\leq k_n(\lambda_n)}\mu_n\big(\fT_{n,\lambda_n}^{(j)}\big)
\geq
\sum_{j=1}^{\cE_{\max}}\mu_n\big(f_j\big)
\geq
\sum_{i=1}^{\cE_{\max}}X_{\sss (i)}~,
\end{align}
where $X_{\sss (1)}<\ldots <X_{\sss (r)}$ are the order statistics corresponding to $(X_1,\ldots, X_r)$.
Now for any $\eps\in(0,1)$,
\begin{align}\label{eqn:59}
\sum_{i=1}^{\eps r}X_{\sss (i)}\probc 2\cdot\bE\big[\Gamma_{1/2}\ind_{\{\Gamma_{1/2}\leq Q_{\eps}\}}\big],
\end{align}
where $\Gamma_{1/2}\sim\mathrm{Gamma}(1/2, 1)$ and $\pr\big(\Gamma_{1/2}\leq Q_{\eps}\big)=\eps$.
It now follows from \eqref{eqn:189}, \eqref{eqn:58} and \eqref{eqn:59} that
\[
\cE_{\max}/r\probc 0,
\]
which in turn implies \eqref{eqn:55}.
This completes the proof.
\qed

\subsection{GHP convergence of the MST of $\cgnthree$}\label{sec:ghp-convergence-simple}
In this section we prove Theorem~\ref{thm:mst_simple}.
We first state two fundamental results about the configuration model and uniform simple graphs with prescribed degree.
\begin{enumeratea}
	\item From \eqref{eqn:cm-distribution} (see also \cite{Boll-book,mcdiarmid1998concentration}), it follows that conditional on being simple, the configuration model has the same distribution as $\cgnd$, i.e.,
	\begin{equation}
	\label{eqn:cm-conditional-uniform}
		\bP\big(\cgnd\in \cdot\big)=\bP\big(\cmnd \in \cdot\ \big|\ \cmnd \text{ is simple}\big).
	\end{equation}
    \item By \cite[Theorem 1.1]{janson2009probability}, if 
    $\sum_{v\in [n]}d_v\to\infty$ and $\sum_{v\in [n]}d_v^2=O(\sum_{v\in [n]}d_v)$, then 
	\begin{equation}
	\label{eqn:cm-simple-0}
		\liminf_{n\to\infty}\ \bP\big(\cmnd \text{ is simple}\big)>0.
	\end{equation}
\end{enumeratea}

Let $\fM_{n, \lambda}^{\avail}$ and $\fM_{n, \lambda}^{\attach}$ be as defined around \eqref{eqn:p-v}.
Define the spaces $\overline\fM_{n, \lambda}^{\avail}$ and $\overline\fM_{n, \lambda}^{\attach}$ analogously for $\cgnthree$.
Using \eqref{eqn:cm-conditional-uniform} and \eqref{eqn:cm-simple-0}, it follows that the analogues of Lemma~\ref{lem:ghp-lemma-1} and Lemma~\ref{lem:ghp-lemma-3} hold for $\cgnthree$:
Fix $\delta>0$ and consider a $\bZ_{>0}$-valued sequence $\lambda_n\uparrow\infty$ with $\lambda_n\leq\min\{\lambda_n^{\star}, \lambda_n^{\circ}\}$.
Then
\begin{align*}
&\pr\big(d_{\GHP}\big(\cbd_{\infty}\big(\cgnthree\big),\ \overline\fM_{n, \lambda_n}^{\attach}\big)>\delta n^{1/3}\big)\\
&\hskip30pt
=
\pr\big(d_{\GHP}\big(\cbd_{\infty}(\cmnthree),\ \fM_{n, \lambda_n}^{\attach}\big)>\delta n^{1/3}\ \big|\ \cmnthree \text{ is simple}\big)
\to 0,
\end{align*}
as $n\to\infty$.
Similarly
\[n^{-1/3}d_{\GHP}\big(\overline\fM_{n, \lambda_n}^{\attach},\ \overline\fM_{n, \lambda_n}^{\avail}\big)\probc 0.\]
To complete the proof, we have to show that the analogue of Lemma~\ref{lem:ghp-lemma-2} remains true for $\overline\fM_{n, \lambda_n}^{\avail}$.
Thus, it suffices to prove that for each fixed $\lambda\in\bR$,
\[
n^{-1/3}\cdot\overline\fM_{n, \lambda}^{\avail}
\weakc
6^{1/3}\cdot\cb^{\infty}\big(S_1\big((48)^{1/3}\cdot \lambda\big)\big),\ \ \text{ as }\ \ n\to\infty
\]
w.r.t. GHP topology.
Let $f:\fS_{\GHP}\to\bR$ be bounded continuous.
Then it suffices to show that as $n\to\infty$,
\[
\bE\big[f\big(n^{-1/3}\fM_{n, \lambda}^{\avail}\ \big|\ \cmnthree \text{ is simple}\big)\big]
-
\bE\big[f\big(n^{-1/3}\fM_{n, \lambda}^{\avail}\big)\big]
\to 0,
\]
or equivalently
\[
\bE\big[f\big(n^{-1/3}\fM_{n, \lambda}^{\avail}\cdot\ind\big\{\cmnthree \text{ is simple}\big\}\big)\big]
-
\bE\big[f\big(n^{-1/3}\fM_{n, \lambda}^{\avail}\big)\big]\times\pr\big(\cmnthree \text{ is simple}\big)
\to 0.
\]
This can be proved by using techniques similar to the ones \ch{used in the proof of \cite[Theorem 3]{sd-rvdh-jvl-ss};
see the argument given in \cite[Section 7]{sd-rvdh-jvl-ss}.}
We omit the details.

\subsection{Proof of Theorem~\ref{thm:non-atomic}}\label{sec:non-atomic-proof}
Let $M_{\lambda}^{n, \er}$ be as defined at the beginning of Section~\ref{sec:alternate-description-M}.
Using Observation~\ref{observation:percolation}, $M_{\lambda}^{n, \er}$ is a subtree of $M_{\infty}^{n, \er}$.
Consider the forest obtained from $M_{\infty}^{n, \er}$ by deleting all edges in $M_{\lambda}^{n, \er}$, and for every $v\in V(M_{\lambda}^{n, \er})$, let $T_{v, \lambda}^{n,\er}$ be the tree in this forest that contains $v$. 
Define $p_{v,\lambda}^{n,\er}:=|T_{v, \lambda}^{n,\er}|/n$.
We now state two lemmas that will be needed in the proof.
\begin{lem}\label{lem:exchangeable-erdos-renyi}
For every $\lambda\in\bR$, conditional on $\ER(n,\lambda)$, the family $\big(p_{v, \lambda}^{n,\er}\, ;\, v\in V(M_{\lambda}^{n,\er})\big)$ of random variables is exchangeable.
\end{lem}

\begin{rem}
\ch{
At the beginning of \cite[Section 4.4]{AddBroGolMie13}, it is stated that 
$\big(|T_{v, \lambda}^{n,\er}|,\, v\in V(M_{\lambda}^{n, \er})\big)$ is exchangeable conditional on 
$V(M_{\lambda}^{n, \er})$. 
(Here, we have translated the claim in \cite{AddBroGolMie13} using our notation.)
However, in \cite[Page 3114]{AddBroGolMie13}, the vertices of \chl{$M_{\lambda}^{n, \er}$} are \chl{relabeled so that the vertices in each element of a given cover $(B_{\lambda}^{n, i},1 \le i \le N^n_\eps)$ of $M_{\lambda}^{n, \er}$ by small-diameter sets} appear successively.
This labeling contains some information about the relative positioning of the vertices in the tree, because if two vertices are close in this arrangement of the vertices, then their tree distance is likely to be small as well.
In other words, the symmetry between the roles of a pair of vertices appearing consecutively and a pair of vertices that are far from each other in this arrangement does not follow directly.
Exchangeability of $\big(|T_{v, \lambda}^{n,\er}|,\, v\in V(M_{\lambda}^{n, \er})\big)$ is needed conditional on this labeling.
Thus, the result implicitly being used in the proof of \cite[Proposition 4.8]{AddBroGolMie13} is the following:
$\big(|T_{v, \lambda}^{n,\er}|,\, v\in V(M_{\lambda}^{n, \er})\big)$ is exchangeable conditional on 
$M_{\lambda}^{n, \er}$.
This stronger form of exchangeability follows from Lemma~\ref{lem:exchangeable-erdos-renyi}.
}
\end{rem}

\noindent{\bf Proof of Lemma~\ref{lem:exchangeable-erdos-renyi}:}
We outline the proof briefly.
Using Lemma~\ref{lem:mst-minimax-criterion}, conditional on the graph $\ER(n,\lambda)$, $M_{\infty}^{n,\er}$ can be generated as follows:

\medskip

\begin{inparaenumii}
\noindent\item 
Let $\cE_{\text{out}}$ denote the set of edges of the complete graph $K_n$ whose endpoints are in two different components of $\ER(n,\lambda)$.
Let $\cE_{\text{in}}$ denote the set of edges of $\ER(n,\lambda)$.
Construct the graph $\ER(n,\lambda)\cup\cE_{\text{out}}$.
Assign i.i.d. $\text{Uniform}[n^{-1}+\lambda n^{-4/3}, 1]$ weights to the edges in $\cE_{\text{out}}$, and independently of this, assign i.i.d. $\text{Uniform}[0, n^{-1}+\lambda n^{-4/3}]$ weights to the edges in $\cE_{\text{in}}$.
Denote the weight assigned to an edge $e$ by $w_e$.

\noindent\item 
From the graph $\ER(n,\lambda)\cup\cE_{\text{out}}$, delete all edges $e\in\cE_{\text{out}}$ that are part of a cycle $\pi$ (with no repeated edges) in $\ER(n,\lambda)\cup\cE_{\text{out}}$ and $w_e$ is the maximum among all edge weights in $\pi$.

\noindent\item 
For each $i\geq 1$, delete all edges $e\in E(\cC_i^{n,\er}(\lambda))$ that are part of a cycle $\pi$ (with no repeated edges) in $\cC_i^{n,\er}(\lambda)$ and $w_e$ is the maximum among all edge weights in $\pi$.
\end{inparaenumii}

\medskip

The marginal distribution of the resulting tree will be the same as that of $M_{\infty}^{n,\er}$.
Consider two distinct vertices $v_1, v_2\in V(\cC_1^{n,\er}(\lambda))$. 
If we interchange the values $w_{\{v_1, u\}}$ and $w_{\{v_2, u\}}$ for every vertex $u\notin V(\cC_1^{n,\er}(\lambda))$, then it is easy to check that in the above procedure, the set of edges removed in step (iii) remains the same, and the set of edges in $\cE_{\text{out}}$ that are not incident to $v_1$ or $v_2$ and are removed in step (ii) remains the same.
Further, if $\{v_1, u\}\in \cE_{\text{out}}$ was removed in step (ii) before the interchange of edge weights, then the edge $\{v_2, u\}$ will be removed in step (ii) after the interchange and vice versa. 
Consequently, the values of $p_{v_1,\lambda}^{n,\er}$ and $p_{v_2,\lambda}^{n,\er}$ would be swapped as a result of the interchange of edge weights.
This shows that conditional on $\ER(n,\lambda)$, the law of $\big(p_{v, \lambda}^{n,\er},\ v\in V(M_{\lambda}^{n,\er})\big)$ is invariant under transpositions.
We can repeat the same argument with any permutation of $V(\cC_1^{n,\er}(\lambda))$ to get the claimed exchangeability. 
\qed

\begin{lem}[Lemma 4.11 of \cite{AddBroGolMie13}]\label{lem:max-mass-small}
Let $\Delta_{n,\lambda}:=\max_{v\in V(M_{\lambda}^{n,\er})} p_{v, \lambda}^{n,\er}$.
Then for every $\delta>0$,
\[
\limsup_{\lambda\to\infty}\ \limsup_{n\to\infty}\ \pr\big(\Delta_{n,\lambda}>\delta\big)=0
\]
\end{lem}

We are now ready to prove Theorem~\ref{thm:non-atomic}.
Observe the following facts:

\begin{inparaenumaa}
\noindent\item
Fix $s\geq 2$ and let $r=3(s-1)$.
Let $e_1,\ldots, e_r$ be an enumeration of $e(\cH^{(s)})$.
Then
\[
\big(\len(e_i),\ 1\leq i\leq r\big)
\equald
\Bigg(Y_i\cdot\Bigg(\frac{\Gamma_{1/2}^{(i)}}{\sum_{j=1}^r \Gamma_{1/2}^{(j)}}\Bigg)^{1/2},\ 1\leq i\leq r\Bigg),
\]
where $Y_i$, $1\leq i\leq r$, are i.i.d. Rayleigh random variables independent of $\Gamma_{1/2}^{(i)}$, $1\leq i\leq r$, which are i.i.d. $\text{Gamma}(1/2, 1)$ random variables.
As observed in \eqref{eqn:41}, $\sqrt{2}\cdot Y_i\sqrt{\Gamma_{1/2}^{(i)}}$, $1\leq i\leq r$, are i.i.d. $\text{Exponential}(1)$ random variables.
Thus, for all $\delta>0$,
\[
\lim_{s\to\infty}\ \pr\Big(\min_{e\in e(\cH^{(s)})}\len(e)\geq s^{-\frac{3}{2}-\delta}\Big)=1.
\]

\noindent\item
By \eqref{eqn:73}, for any $s\geq 2$,
\[
\frac{1}{\sqrt{m}}\min_{e\in e(\cH_{m, s})}\len(e)\weakc\min_{e\in e(\cH^{(s)})}\len(e)
\ \text{ as }\ m\to\infty.
\]

\noindent\item
$\cC_1^{n,\er}(\lambda)$ can be generated as follows:
\begin{inparaenumii}
	\item 
Sample $|\cC_1^{n,\er}(\lambda)|$ and $\mathrm{sp}(\cC_1^{n,\er}(\lambda))$. 
Denote the realizations by $m$ and $s$ respectively.
\item 
Conditional on the previous step, generate $\cH_{m, s}$ and set this graph to be $\cC_1^{n,\er}(\lambda)$.	
\end{inparaenumii}

\noindent\item
By Lemma~\ref{lem:growth-gamma-N},
$\pr\big(\xi_1(\lambda)\leq\lambda\big)+\pr\big(N_1(\lambda)<\lambda^3/2\big)\to 0$ as $\lambda\to\infty$.

\noindent\item
By Theorem~\ref{thm:erdos-renyi-component-sizes},
$\Big(n^{-2/3}|\cC_1^{n,\er}(\lambda)|,\ \mathrm{sp}\big(\cC_1^{n,\er}(\lambda)\big)\Big)
\weakc \big(\xi_1(\lambda), N_1(\lambda)\big)$ as $n\to\infty$.
\end{inparaenumaa}

Combining the above, we see that
\begin{align}\label{eqn:643}
\limsup_{n\to\infty}\ \pr\Big(\min_{e\in e\big(\cC_1^{n,\er}(\lambda)\big)}\len(e)\leq n^{1/3}/\lambda^5\Big)
=:\eps_1(\lambda)\to 0,
\ \text{ as }\ \lambda\to\infty.
\end{align}

Using the convergences $\mathrm{sp}(\cC_1^{n,\er}(\lambda))\weakc N_1(\lambda)$ as $n\to\infty$ and $N_1(\lambda)/\lambda^3\probc 2/3$ as $\lambda\to\infty$ together with \eqref{eqn:713}, we see that
\begin{align}\label{eqn:645}
&\lim_{n\to\infty}\ \pr\Big(k\big(\cC_1^{n,\er}(\lambda)\big)\text{ is not a 3-regular multigraph}\Big)\notag\\
&\hskip100pt
=\pr\big(N_1(\lambda)\leq 1\big)
=:\eps_2(\lambda)\to 0,
\text{ as }\lambda\to\infty.
\end{align}

Let $e_i^{n,\er}(\lambda)$, $\ch{1\leq i\leq |e(\cC_1^{n,\er}(\lambda))|}$, be an enumeration of $e(\cC_1^{n,\er}(\lambda))$.
Let $V_i^{n,\er}(\lambda)$ be the set of vertices in $\cC_1^{n,\er}(\lambda)$ that are connected to $\core\big(\cC_1^{n,\er}(\lambda)\big)$ via $e_i^{n,\er}(\lambda)$. (As before, the common endpoints of multiple $e\in e(\cC_1^{n,\er}(\lambda))$ and their pendant subtrees are assigned to only of the $V_i^{n,\er}(\lambda)$'s in an arbitrary way.)
From \eqref{eqn:711} and arguments as above, 
\begin{align}\label{eqn:644}
\limsup_{n\to\infty}\ \pr\Big(\max_i \frac{|V_i^{n,\er}(\lambda)|}{|\cC_1^{n,\er}(\lambda)|}\geq (\log\lambda)^2/\lambda^3\Big)
=:\eps_3(\lambda)\to 0,
\ \text{ as }\ \lambda\to\infty.
\end{align}

Denote the complements of the events in \eqref{eqn:643},  \eqref{eqn:645}, and  \eqref{eqn:644} by $E_{n,\lambda}^{\sss (1)}, E_{n,\lambda}^{\sss (2)}$, and $E_{n,\lambda}^{\sss (3)}$ respectively.
Let $E_{n,\lambda}:=\cap_{j=1}^3 E_{n,\lambda}^{\sss (j)}$.
Note that on the event $E_{n,\lambda}$, any $\cU\subseteq V(\cC_1^{n,\er}(\lambda))$ with $\diam(\cU; \cC_1^{n,\er}(\lambda))\leq n^{1/3}/(2\lambda^5)$ can intersect $V_i^{n,\er}(\lambda)$ for at most three values of $i$.

Let $\cZ_{\lambda}^{n,\er}$ be the graph obtained by attaching, for each $v\in\cC_1^{n,\er}(\lambda)$, the tree $T_{v, \lambda}^{n,\er}$ to $\cC_1^{n,\er}(\lambda)$ via identification of the vertices labeled $v$.
Consider $\widetilde\cU\subseteq [n]$ with $\diam(\widetilde\cU; \cZ_{\lambda}^{n,\er})\leq n^{1/3}/(2\lambda^5)$.
Let 
\[
\cU=\big\{v\in\cC_1^{n,\er}(\lambda)\ :\ T_{v,\lambda}^{n,\er}\cap\widetilde\cU\neq\emptyset\big\}.
\]
Then $\diam(\cU; \cC_1^{n,\er}(\lambda))\leq n^{1/3}/(2\lambda^5)$.
Consequently, on the event $E_{n,\lambda}$,
\begin{align}\label{eqn:888}
\frac{1}{n}\cdot |\widetilde\cU|
&\leq
\sum_{v\in\cU}\frac{|T_{v,\lambda}^{n,\er}|}{n}
\leq   
3\times\max_i\Big(\sum_{v\in V_i^{n,\er}(\lambda)}p_{v,\lambda}^{n,\er}\Big)\\
&
\leq 
3\cdot\max_i\Bigg|
\sum_{v\in V_i^{n,\er}(\lambda)}
p_{v,\lambda}^{n,\er}
-\frac{|V_i^{n,\er}(\lambda)|}{|\cC_1^{n,\er}(\lambda)}
\Bigg|
+3(\log\lambda)^2\frac{1}{\lambda^3}~.
\end{align}

Arrange the vertices in $\cC_1^{n,\er}(\lambda)$ in a sequence so that for each $i$, the vertices in $V_i^{n,\er}(\lambda)$ appear consecutively.
This arrangement is measurable w.r.t. the sigma field generated by $\ER(n,\lambda)$.
By Lemma~\ref{lem:exchangeable-erdos-renyi}, conditional on this arrangement, $\big\{p_{v,\lambda}^{n,\er}\big\}_{v\in V(\cC_1^{n,\er}(\lambda))}$ is an exchangeable sequence.
Using \eqref{eqn:concentration} with $x=\Delta_{n,\lambda}^{-1/4}$, we see that
\begin{align}\label{eqn:889}
&
\eps_4^{\sss (n)}(\lambda):=
\pr\Bigg(
\max_i\Big|
\sum_{v\in V_i^{n,\er}(\lambda)}
p_{v,\lambda}^{n,\er}
-\frac{|V_i^{n,\er}(\lambda)|}{|\cC_1^{n,\er}(\lambda)}\Bigg|
\geq 
2\Delta_{n,\lambda}^{1/4}
\Bigg)\\
&\hskip40pt
\leq 
\pr\Big(\Delta_{n,\lambda}^{-1/4}\leq c_2\Big)
+
2\cdot\bE\Big[\exp\Big(-c\Delta_{n,\lambda}^{-1/4}\log\log\Delta_{n,\lambda}^{-1/4}\Big)\Big],\notag
\end{align}
where $c_2$ is as in Lemma~\ref{lem:concentration-uniform-permutation}.
Combining \eqref{eqn:888} and \eqref{eqn:889}, we see that
\[
\pr\Big(\fm\big((2\lambda^5)^{-1}\ ;\ n^{-1/3}\cZ_{\lambda}^{n,\er} \big)
\geq 
6\Delta_{n,\lambda}^{1/4}+3\cdot(\log\lambda)^2\cdot\lambda^{-3}
\Big)
\leq 
\pr(E_{n,\lambda}^c)+\eps_4^{\sss (n)}(\lambda).
\]
Since $M_{\infty}^{n,\er}$ is a subtree of $\cZ_{\lambda}^{n,\er}$, 
$\fm(\delta; M_{\infty}^{n,\er})\leq\fm(\delta; \cZ_{\lambda}^{n,\er})$ for every $\delta>0$.
Thus, using Theorem~\ref{thm:mst_complete}, we conclude that for every $\eta>0$ and $\lambda>1$,
\begin{align}\label{eqn:890}
&\pr\Big(\fm\big((4\lambda^5)^{-1}\ ;\ \cM\big)
\geq 
\eta+3\cdot(\log\lambda)^2\cdot\lambda^{-3}
\Big)\\
&\hskip50pt
\leq 
\limsup_{n\to\infty}\ \pr\big(6\Delta_{n,\lambda}^{1/4}>\eta\big)
+\sum_{i=1}^3\eps_i(\lambda)
+\limsup_{n\to\infty}\ \eps_4^{\sss (n)}(\lambda).\notag
\end{align}
The result follows upon using Lemma~\ref{lem:max-mass-small} and noting that the right side of \eqref{eqn:890} tends to zero as $\lambda\to\infty$.

\subsection{Proof of Theorem~\ref{thm:derandomization}}\label{sec:derandomization}
For $p\in(0,1)$ and $m\in\bN$, let $G^m_p$ be distributed as follows:
For any connected graph $H$ on $[m]$ having $r$ edges,
\[\pr(G^m_p=H)\propto p^r(1-p)^{-r}.\]
That is, $G^m_p$ is an \erdos random graph conditioned to be connected.
We start with the following lemma:
\begin{lem}\label{lem:surplus-leq-1}
If $p m^{3/2}\leq 1$, then
\begin{align}\label{eqn:3}
\pr\big(\mathrm{sp}(G^m_p)\geq 2\big)\leq C p^2 m^3
\end{align}
for some universal constant $C$.
Consequently, for any $\eps>0$ there exists $\lambda(\eps)>0$ such that for all $\lambda\geq\lambda(\eps)$,
\begin{align}\label{eqn:4}
\limsup_{n\to\infty}\ \pr\big(\mathrm{sp}\big(\cC_{i}^{\sss n,\er}(\lambda)\big)\geq 2\ \text{ for some }\ i\geq 2\big)\leq \eps.
\end{align}
\end{lem}

\noindent{\bf Proof:}
Suppose $\vt$ is a rooted tree on $[m]$.
For $v\in [m]$ define
\begin{gather}
R\big(\stackrel{\leftarrow\sss{(k)}}{v\ \ \ }, v, \vt\big):=\big\{u\in [m]\ :\ \overleftarrow{u}=\stackrel{\leftarrow\sss{(k)}}{v\ \ \ }\text{ and }u>\stackrel{\leftarrow\sss{(k-1)}}{v\ \ \ \ }\big\},
\ \text{ for }\ 1\leq k\leq\height(v,\vt),\ \text{ and} \label{eqn:722}
\\
R(v,\vt):= \bigcup_{k=1}^{\height(v,\vt)}R\big(\stackrel{\leftarrow\sss{(k)}}{v\ \ \ }, v, \vt\big).
\label{eqn:72}
\end{gather}
Let $g(\vt):=\sum_{v\in [m]}\big|R(v,\vt)\big|$ and $MR(\vt)=\max_{v\in [m]}\big|R(v,\vt)\big|$.
Let $T_m$ denote a  uniform rooted tree on $[m]$ and let $\widetilde T_m$ be distributed as
\begin{align}\label{eqn:1}
\bP\big(\widetilde T_m=\vt\big)=\frac{(1-p)^{-g(\vt)}\pr\big(T_m=\vt\big)}{\bE\big[(1-p)^{-g(T_m)}\big]}.
\end{align}
Then by \cite[Proposition 8]{BBG-12}, $G^m_p$ has the same law as the random graph obtained from $\widetilde T_m$ by placing an edge with probability $p$ independently between every pair of vertices $v$ and $u$, where $v\in[m]$ and $u\in R(v,\widetilde T_m)$, and then forgetting the identity of the root of $\widetilde T_m$.
In particular, $\mathrm{sp}(G^m_p)$ is distributed as $\mathrm{Binomial}(N, p)$ where $N\equald g(\widetilde T_m)$.
Hence
\begin{align}\label{eqn:2}
\pr\big(\mathrm{sp}(G^m_p)\geq 2\big)
\leq p^2 \bE\big[g(\widetilde T_m)^2\big]
\leq p^2m^2\bE\big[MR(\widetilde T_m)^2\big].
\end{align}
It follows from \eqref{eqn:1} that for any $x>0$,
\begin{align*}
\pr\Big(MR(\widetilde T_m)\geq x\sqrt{m}\Big)
&\leq \bE\Big[(1-p)^{-g(T_m)}\ind_{\big\{MR(T_m)\geq x\sqrt{m}\big\}}\Big]\\
&\leq\Big[\pr\Big(MR(T_m)\geq x\sqrt{m}\Big)\Big]^{\frac{1}{2}}\Big[\bE (1-p)^{-2g(T_m)}\Big]^{\frac{1}{2}}
\leq C e^{-C'x^2}e^{C'' p^2 m^3},
\end{align*}
where the last step uses \cite[Lemmas 13 and 14]{BBG-12}.
Using the fact that $p^2 m^3\leq 1$, it follows that $\bE\big[MR(\widetilde T_m)^2\big]\leq Cm$.
This in conjunction with \eqref{eqn:2} yields \eqref{eqn:3}.

Next note that by \cite[Theorem A.1]{janson08susceptibility}, for any $\eps\in(0,1)$, there exists $\lambda(\eps)>0$ such that for all $\lambda\geq\lambda(\eps)$,
\begin{align}\label{eqn:5}
\limsup_{n\to\infty}\ \pr\Big(\sum_{i\geq 2}\big|\cC_{i}^{\sss n,\er}(\lambda)\big|^2\geq\eps n^{4/3}\Big)\leq\eps.
\end{align}
Denote the event on the left side of \eqref{eqn:5} by $F(n,\lambda; \eps)$.
Now conditional on the component sizes of $\ER(n,\lambda)$ being equal to $m_1, m_2,\ldots, m_r$, the components are distributed as $G^{m_1}_p,\ldots, G^{m_r}_p$ with $p=n^{-1}+\lambda n^{-4/3}$.
Further, for any $\eps\in (0,1/4)$ and large $n$, $pm_i^{3/2}\leq 1$ if $m_i^2\leq\eps n^{4/3}$.
Hence, using \eqref{eqn:3},
\begin{align}\label{eqn:6}
&\pr\big(\mathrm{sp}\big(\cC_{i}^{\sss n,\er}(\lambda)\big)\geq 2\ \text{ for some }\ i\geq 2\big)\\
&\hskip10pt\leq \pr\big(F(n,\lambda; \eps)\big)
+C\ \bE\bigg[\ind{\big\{F(n,\lambda; \eps)^c\big\}}\cdot n^{-2}\sum_{i\geq 2}\big|\cC_{i}^{\sss n,\er}(\lambda)\big|^3\bigg]
\leq \pr\big(F(n,\lambda; \eps)\big)+C\eps.\nonumber
\end{align}
\eqref{eqn:6} together with \eqref{eqn:5} yields \eqref{eqn:4}.
\qed

\medskip

Our next lemma roughly states that inside the critical window, the number of surplus edges in the largest component of the \erdos random graph takes all large integer values with high probability, and during this time of the evolution, every other component is either a tree or is unicyclic.

\begin{lem}\label{lem:surplus-hits-all}
For every $\eps>0$, there exists $s_{\eps}\in\bN$ such that for all integers $s\geq s_{\eps}$,
\begin{align}\label{eqn:13}
\liminf_n\ \pr\big(\cA[s_{\eps}, s]\big)\geq 1-\eps,
\end{align}
where $\cA[s_{\eps}, s]$ denotes the event that there exist $\lambda_1\leq\lambda_2$ such that in the interval $[\lambda_1,\lambda_2]$, the process $\mathrm{sp}\big(\cC_{1}^{\sss n,\er}(\cdot)\big)$ assumes all values in $\big\{s_{\eps}, s_{\eps}+1,\ldots,s\big\}$, and $\mathrm{sp}\big(\cC_{i}^{\sss n,\er}(\lambda)\big)\leq 1$ for all $i\geq 2$ and $\lambda\in[\lambda_1,\lambda_2]$.
\end{lem}

\noindent{\bf Proof:}
For $k\geq 1$, define
\[\beta(k)=\frac{k(k+1)}{\big(k+1/6\big)\big(k+5/6\big)}.\]
We say that the ``leader changes in $\ER(n,\cdot)$ after time $\lambda$" if there exists $\lambda'>\lambda$ such that the component in $\ER(n,\lambda')$ containing $\cC_{1}^{\sss n,\er}(\lambda)$ is not $\cC_{1}^{\sss n,\er}(\lambda')$.
Fix $\eta>0$ and choose $\lambda(\eta)$ large such that the following hold:
\begin{gather}
\limsup_{n\to\infty}\ \pr\Big(\text{The leader does not change in }\ER(n,\cdot)\text{ after time }\lambda(\eta)\Big)\geq 1-\eta,
\label{eqn:9}\\
\prod_{j\geq 0}\beta\Big(\big[\lambda(\eta)^3/3\big]+j\Big)\geq 1-\eta,
\label{eqn:10}\\
\limsup_{n\to\infty}\ \pr\Big(\mathrm{sp}\big(\cC_{i}^{\sss n,\er}(\lambda(\eta))\big)\leq 1\ \text{ for all }\ i\geq 2\Big)\geq 1-\eta,\ \text{ and}
\label{eqn:7}\\
\limsup_{n\to\infty}\ \pr\Big(\lambda(\eta)^3/3<
\mathrm{sp}\big(\cC_{1}^{\sss n,\er}(\lambda(\eta))\big)<\lambda(\eta)^3 \Big)\geq 1-\eta.
\label{eqn:8}
\end{gather}
\eqref{eqn:9} uses \cite[Theorem 7]{luczak1990component} (see also \cite{addarioberry-bhamidi-sen}).
\eqref{eqn:10} uses the fact that $\prod_{k\geq 1}\beta(k)>0$.
\eqref{eqn:7} uses  \eqref{eqn:4}.
\eqref{eqn:8} uses Lemma~\ref{lem:growth-gamma-N} and Theorem~\ref{thm:erdos-renyi-component-sizes}.

Let $F^n_k$ denote the event that there exists $\lambda\in\bR$ such that the random graph $\ER(n,\lambda)$ has exactly one component with surplus $k+1$ and the surplus of every other component is at most one. Then by \eqref{eqn:7} and \eqref{eqn:8},
\begin{align}\label{eqn:11}
\liminf_{n\to\infty}\ \pr\Big(\bigcup\displaystyle_{k=[\lambda(\eta)^3/3]}^{[\lambda(\eta)^3]}F^n_k\Big)\geq 1-2\eta.
\end{align}
By \cite[Theorem 5.28]{janson2000random}, for $[\lambda(\eta)^3/3]\leq k\leq [\lambda(\eta)^3]$ and any $s\geq [\lambda(\eta)^3]+1$,
\begin{align}\label{eqn:12}
\lim_{n\to\infty}\pr\Big(\bigcap_{j=k+1}^{s} F^n_j\ \Big|\ F^n_k\setminus\big(\bigcup_{j=[\lambda(\eta)^3/3]}^{k-1} F^n_j\big)\Big)
=\prod_{j=k+1}^{s}\beta(j).
\end{align}
Note that
\begin{align*}
\pr\Big(\bigcap_{j=[\lambda(\eta)^3]}^{s} F^n_j\Big)
&
\geq
\sum_{k=[\lambda(\eta)^3/3]}^{[\lambda(\eta)^3]}
\pr\Big(\big(\bigcap_{j=k}^{s} F^n_j\big)\setminus\big(\bigcup_{j=[\lambda(\eta)^3/3]}^{k-1} F^n_j\big)\Big)
\\
&\geq\sum_{k=[\lambda(\eta)^3/3]}^{[\lambda(\eta)^3]}
\pr\Big(\bigcap_{j=k+1}^{s} F^n_j\ \Big|\ F^n_k\setminus\big(\bigcup_{j=[\lambda(\eta)^3/3]}^{k-1} F^n_j\big)\Big)\cdot
\pr\Big(F^n_k\setminus\big(\bigcup_{j=[\lambda(\eta)^3/3]}^{k-1} F^n_j\big)\Big)
\\
&\geq
\pr\Big(\bigcup_{k=[\lambda(\eta)^3/3]}^{[\lambda(\eta)^3]}F^n_k\Big)
\times\min\displaystyle_{\star}\,
\pr\Big(\bigcap_{j=k+1}^{s} F^n_j\ \Big|\ F^n_k\setminus\big(\bigcup_{j=[\lambda(\eta)^3/3]}^{k-1} F^n_j\big)\Big),
\end{align*}
where $\min_{\star}$ is minimum taken over $[\lambda(\eta)^3/3]\leq k\leq [\lambda(\eta)^3]$.
Thus, \eqref{eqn:11}, \eqref{eqn:12}, and \eqref{eqn:10} give
\[\liminf_{n\to\infty}\ \pr\Big(\bigcap_{j=[\lambda(\eta)^3]}^{s} F^n_j\Big)\geq (1-2\eta)(1-\eta).\]
Combining this with \eqref{eqn:7} and \eqref{eqn:8}, we see that
\begin{align}
&
\liminf_{n\to\infty}\
\pr\Big(\Big\{\bigcap_{j=[\lambda(\eta)^3]}^{s} F^n_j\Big\}
\bigcap\Big\{\mathrm{sp}\big(\cC_{i}^{\sss n,\er}\big(\lambda(\eta)\big)\big)\leq 1\ \text{ for all }\ i\geq 2\Big\}\\
&\hskip120pt
\bigcap\Big\{\lambda(\eta)^3/3< \mathrm{sp}\big(\cC_{1}^{\sss n,\er}(\lambda(\eta))\big)<\lambda(\eta)^3\Big\}\Big)
> 1-5\eta.
\end{align}
Thus, for all large $n$, in the process $\big(\ER(n,\lambda), \lambda\geq\lambda(\eta)\big)$, with probability at least $1-5\eta$, the surplus of the component containing $\cC_{1}^{\sss n,\er}(\lambda(\eta))$ assumes all values in $\big\{[\lambda(\eta)^3]+1,\ldots, s\big\}$, and during this part of the evolution, the surplus of the other components remains at most one.
By \eqref{eqn:9}, the component containing $\cC_{1}^{\sss n,\er}(\lambda(\eta))$ remains the largest component after time $\lambda(\eta)$ with probability at least $1-\eta$.
Thus, \eqref{eqn:13} follows if we take $s_{\eps}=[\lambda(\eta)^3]+1$ with $\eta=\eps/6$.
\qed

\medskip

Let $\cH_{n,s}$ be as in Theorem~\ref{thm:scaling-connected} and let $L(\cdot)$ be as in \eqref{eqn:def-L}.
Define $\widetilde\cH_{n,s}$ and $\widetilde\cH^{(s)}$ via
\[\bE\big[f(\widetilde\cH_{n,s})\big]=\frac{\bE\big[f(\cH_{n,s})L(\cH_{n,s})\big]}{\bE\big[L(\cH_{n,s})\big]}
\ \text{ and }\
\bE\big[f(\widetilde\cH^{(s)})\big]=\frac{\bE\big[f(\cH^{(s)})L(\cH^{(s)})\big]}{\bE\big[L(\cH^{(s)})\big]}\]
for every bounded measurable $f:\fS_{\GHP}\to\bR$.
For $s\geq 2$ define
\begin{align}\label{eqn:def-tau-k}
\tau_s=\inf\big\{\lambda\ :
& \ \mathrm{sp}(\cC^{\star})=s\text{ for some component }\cC^{\star}\text{ of }\ER(n, \lambda)\\
&\ \text{ and } \mathrm{sp}(\cC)\leq 1\text{ for every other component }\cC\text{ of }\ER(n, \lambda)\big\}.\nonumber
\end{align}
If $\tau_s<\infty$, define $\cC_{\tau_s}^{n, \star}$ to be the (unique) component of $\ER(n, \tau_s)$ with $\mathrm{sp}(\cC_{\tau_s}^{n, \star})=s$.
If $\tau_s=\infty$, define $\cC_{\tau_s}^{n, \star}$ to be the one-point space.

\begin{lem}\label{lem:core-biased}
Fix $s\geq 2$.
Let $U_{ij}$, $1\leq i<j\leq n$, be the i.i.d. $\mathrm{Uniform}[0,1]$ random variables used in the construction of $\ER(n,\cdot)$.
Let $\widetilde\cH_{n,s}$ be independent of $(U_{ij}$, $1\leq i<j\leq n)$.
Define $M_{\tau_s}^{n,\star}$ to be the MST of $\cC_{\tau_s}^{n, \star}$ constructed using the weights $U_{ij}$ if $\tau_s<\infty$ and $|\cC_{\tau_s}^{n, \star}|\geq\log n$, and set $M_{\tau_s}^{n,\star}=\cbd_{\infty}\big(\widetilde\cH_{n,s}\big)$ otherwise.
Then as $n\to\infty$,
\[\big(|M_{\tau_s}^{n,\star}|\big)^{-1/2}M_{\tau_s}^{n,\star}\weakc\cb^{\infty}(\widetilde\cH^{(s)})\ \ \
\text{ w.r.t. GHP topology.}\]
\end{lem}

\noindent{\bf Proof:}
For convenience, we will assume that the random vector $(U_{ij}, 1\leq i<j\leq n)$ is given by the identity map on the canonical probability space $[0,1]^{{n}\choose{2}}$ endowed with the ${n}\choose{2}$-fold product of the uniform measure on $[0,1]$.

For any subgraph $H$ of the complete graph on $[n]$, define the event
\[F_H:=\big\{\tau_s<\infty,\ \ER(n, \tau_s)\setminus\cC_{\tau_s}^{n, \star}=H\big\}.\]
Fix any $H$ with $\pr(F_H)>0$.
Then $\pr\big(F_H\cap\{\cC_{\tau_s}^{n, \star}=H_1\}\big)>0$ for any connected graph $H_1$ with
\begin{align}\label{eqn:66}
V(H_1)=[n]\setminus V(H),\ \text{ and }\ \mathrm{sp}(H_1)=s.
\end{align}
Now for any $H_1$ satisfying \eqref{eqn:66}, the realizations $(u_{ij}\, ;\ 1\leq i<j\leq n)$ of the random variables $U_{ij}$ for which $F_H\cap\{\cC_{\tau_s}^{n, \star}=H_1\}$ holds are given by
\begin{align*}
F_H\cap\{\cC_{\tau_s}^{n, \star}=H_1\}
&=\Big\{
\max\big\{u_{ij}\ :\ (i,j)\in E(H_1)\cup E(H)\big\}\\
&\hskip15pt
=\max\big\{u_{ij}\ :\ (i,j)\in E(\core(H_1))\big\}
<\min\big\{u_{ij}\ :\ (i,j)\notin E(H_1)\cup E(H)\big\}
\Big\},
\end{align*}
and for any such realization $(u_{ij})$, we have $\big(u_{\pi(i,j)}\, ;\ 1\leq i<j\leq n\big)\in F_H\cap\{\cC_{\tau_s}^{n, \star}=H_1\}$
for any permutation $\pi$ of $\{(i,j)\ :\ 1\leq i<j\leq n\}$ satisfying $\pi(i,j)=(i,j)$ for all $(i,j)\notin E(\core(H_1))$.
Hence, conditional on $\tau_s<\infty$ and $\ER(n, \tau_s)$, the random variables $U_{i,j}$, $(i,j)\in\core(\cC_{\tau_s}^{n, \star})$, are exchangeable.
Using Lemma~\ref{lem:cycle-breaking-gives-mst}, we see that the following equality of conditional distributions hold for any $m\geq\log n$:
\begin{align}\label{eqn:70}
\big(M_{\tau_s}^{n,\star}\ \big|\ \tau_s<\infty,\ |\cC_{\tau_s}^{n, \star}|=m\big)
\equald
\big(\cbd_{\infty}\big(\cC_{\tau_s}^{n,\star}\big)\ \big|\ \tau_s<\infty,\ |\cC_{\tau_s}^{n, \star}|=m\big)
\end{align}

Next, for any two graphs $G_1, G_2$ on $[n]$, write $\pr^{\er}(G_1, G_2)$ to denote the probability that there exist $\lambda_1\leq\lambda_2$ such that $\ER(n,\lambda_1)=G_1$ and $\ER(n,\lambda_2)=G_2$.
Thus, if $G_1$ is a subgraph of $G_2$, then
\begin{align}\label{eqn:6565}
\pr^{\er}(G_1, G_2)=
\frac{1}{N!}\cdot|E(G_1)|!\cdot\big(|E(G_2)|-|E(G_1)|\big)!\cdot\big(N-|E(G_2)|\big)!\, ,
\end{align}
where $N={{n}\choose{2}}$.
Now for any $H_1$ satisfying \eqref{eqn:66},
\begin{align*}
&\pr\big(\cC_{\tau_s}^{n, \star}=H_1\ \big|\ F_H\big)
=\frac{1}{\pr(F_H)}\pr\big(\tau_s<\infty,\ \ER(n, \tau_s)=H_1\cup H\big)\\
&\hskip60pt=\frac{1}{\pr(F_H)}
\sum_{e\in E(\core(H_1))}
\pr^{\er}\big((H_1\setminus e)\cup H,\ H_1\cup H\big)
\propto
\big|E(\core(H_1))\big|=L(H_1),
\end{align*}
where in the penultimate step we have used \eqref{eqn:6565} to deduce that the summands are the same for any $H_1$ satisfying \eqref{eqn:66}.
Thus, for any $m\geq\log n$, the conditional distribution of $\cC_{\tau_s}^{n, \star}$ given $\tau_s<\infty$ and $|\cC_{\tau_s}^{n, \star}|=m$ satisfies
\begin{align}\label{eqn:67}
\big(\cC_{\tau_s}^{n, \star}\ \big|\ \tau_s<\infty,\ |\cC_{\tau_s}^{n, \star}|=m\big)
\equald
\widetilde\cH_{m,s}.
\end{align}
Now for any bounded continuous $f:\fS_{\GHP}\to\bR$,
\begin{align*}
\lim_{m\to\infty}\bE\big[f\big(\frac{1}{\sqrt{m}}\widetilde\cH_{m,s}\big)\big]
=\lim_{m\to\infty}\frac{\bE\big[f\big(\frac{1}{\sqrt{m}}\cH_{m,s}\big)L(\cH_{m,s})\big]}{\bE\big[L(\cH_{m,s})\big]}
=
\frac{\bE\big[f\big(\cH^{(s)}\big)L(\cH^{(s)})\big]}{\bE\big[L(\cH^{(s)})\big]}
=\bE\big[f\big(\widetilde\cH^{(s)}\big)\big],
\end{align*}
where the second step uses \eqref{eqn:73} and \eqref{eqn:74}.
Hence $m^{-1/2}\widetilde\cH_{m,s}\weakc\widetilde\cH^{(s)}$ as $m\to\infty$ w.r.t. GHP topology.
Using Theorem~\ref{thm:mst-from-ghp-convergence} and \eqref{eqn:71}, it follows that as $m\to\infty$,
\[m^{-1/2}\cbd_{\infty}\big(\widetilde\cH_{m,s}\big)\weakc\cb^{\infty}\big(\widetilde\cH^{(s)}\big)\]
w.r.t. GHP topology.
Now the claim follows from \eqref{eqn:70} and \eqref{eqn:67}.
\qed

\medskip

\noindent{\bf Proof of Theorem~\ref{thm:derandomization}:}
Fix $0<\eps<1/2$.
For $s\geq 3$, define $\lambda_s$ by the relation $2\lambda_s^3=3s$.
Define
\begin{align}\label{eqn:65}
E_s:=\bigg\{\cA[s,s]
\bigcap\big\{\text{the leader does not change in }\ER(n,\cdot)\text{ after time }\lambda_s/2\big\}\bigg\},
\end{align}
where $\cA[\cdot\, , \cdot]$ is as in Lemma~\ref{lem:surplus-hits-all}.
Using \eqref{eqn:13} and \eqref{eqn:9}, choose $s_1$ large so that
\begin{align}\label{eqn:14}
\limsup_{n\to\infty}\ \pr\big(E_s^c\big)\leq\eps\text{ for all }s\geq s_1.
\end{align}

Let $M_{\lambda}^{n,\er}$ denote the MST of $\cC_{1}^{n,\er}(\lambda)$ constructed using the same i.i.d. $\mathrm{Uniform}[0,1]$ random variables $U_{ij}$ used to construct the process $\ER(n,\cdot)$.
If the leader does not change after time $\lambda$, then using Observation~\ref{observation:percolation}, we see that $M_{\lambda}^{n,\er}$ is a subtree of $M_{\lambda'}^{n,\er}$ for any $\lambda'>\lambda$.
Thus, using \cite[Lemma 4.5]{AddBroGolMie13}, we can choose $s_2$ large enough so that
\begin{align}\label{eqn:15}
\limsup_{n\to\infty}\ \pr\Big(d_H\big(M_{\lambda}^{n,\er}, M_{\lambda'}^{n,\er}\big)>\eps n^{1/3}\Big)\leq\eps
\ \text{ for all }\ \lambda,\ \lambda'\geq\lambda_{s_2}.
\end{align}
Next, define
$\overline\lambda_s=\lambda_s(1+\eps)$ and $\underline\lambda_s=\lambda_s(1-\eps)$,
and using Lemma~\ref{lem:growth-gamma-N}, Theorem~\ref{thm:erdos-renyi-component-sizes},
\eqref{eqn:4}, and \eqref{eqn:65}, choose $s_3$ large enough so that
\begin{align}\label{eqn:16}
\limsup_{n\to\infty}\ \pr\big(F_{s,\eps}^c\big)\leq 2\eps
\ \text{ for all }\ s\geq s_3\, ,
\end{align}
where
\begin{align*}
&F_{s,\eps}:=
E_s\,
\bigcap
\Big\{\big|\cC_{1}^{\sss n,\er}(\underline\lambda_s)\big|\geq 2\underline\lambda_s(1-\eps)n^{2/3}
\ \ \ \text{ and }\ \ \
\big|\cC_{1}^{\sss n,\er}(\overline\lambda_s)\big|\leq 2\overline\lambda_s(1+\eps)n^{2/3}\Big\} \\
&\hskip50pt
\ \bigcap\ 
\Big\{2\leq \mathrm{sp}\big(\cC_{1}^{\sss n,\er}(\underline\lambda_s)\big)\leq s-1
\ \ \ \text{ and }\ \ \
\mathrm{sp}\big(\cC_{1}^{\sss n,\er}(\overline\lambda_s)\big)\geq s+1\Big\}\\
&\hskip100pt
\ \bigcap\ 
\Big\{\forall i\geq 2\, , \ \ \mathrm{sp}\big(\cC_{i}^{\sss n,\er}(\underline\lambda_s)\big)\leq 1
\ \ \ \text{ and }\ \ \
\mathrm{sp}\big(\cC_{i}^{\sss n,\er}(\overline\lambda_s)\big)\leq 1\Big\}\, .
\end{align*}
Set $s_0:=\max\big\{s_2, s_3\big\}$.
From now on, we will only consider $s\geq s_0$.

Let $\tau_s$ be as in \eqref{eqn:def-tau-k}.
If $\tau_s<\infty$, let $M_{\tau_s}^{n,\er}$ be the MST of $\cC_{1}^{\sss n,\er}(\tau_s)$ constructed using the edge weights $U_{ij}$.
If $\tau_s=\infty$, set $\cC_{1}^{\sss n,\er}(\tau_s)$ to be the complete graph $K_n$, and let $M_{\tau_s}^{n,\er}=M_{\infty}^{n, \er}$--the MST of $K_n$ constructed using the edge weights $U_{ij}$.
Note that on the event $F_{s,\eps}$,
$\underline\lambda_s<\tau_s<\overline\lambda_s$,
$\cC_{1}^{n,\er}(\tau_s)=\cC_{\tau_s}^{n, \star}$, and
$M_{\tau_s}^{n,\er}=M_{\tau_s}^{n, \star}$, where the notation is as in Lemma~\ref{lem:core-biased}.
Thus, writing $\cL(\cdot)$ and $d_{\mathrm{PR}}(\cdot,\cdot)$ to denote the law of a random metric measure space and the Prokhorov distance between two measures respectively, it follows from Lemma~\ref{lem:core-biased} that
\begin{align}\label{eqn:23}
\limsup_{n\to\infty}\
d_{\mathrm{PR}}
\bigg(\cL\bigg(\frac{(12s)^{1/6}}{\big(|\cC_{1}^{\sss n,\er}(\tau_s)\big|\big)^{1/2}}M_{\tau_s}^{n,\er}\bigg), \
\cL\bigg(\big(12 s\big)^{1/6}\cdot\cb^{\infty}\big(\widetilde\cH^{(s)}\big)\bigg)\bigg)
\leq 2\eps.
\end{align}

Next note that on $F_{s,\eps}$, 
$M_{\underline\lambda_s}^{n, \er}\subseteq M_{\tau_s}^{n, \er} \subseteq M_{\overline\lambda_s}^{n, \er}$.
On $F_{s,\eps}$, for every $i\in M_{\tau_s}^{n, \er}$, let
\[V_i:=\big\{j\in M_{\overline\lambda_s}^{n, \er}\ :\
\text{the path connecting }j\text{ and }i\text{ in }M_{\overline\lambda_s}^{n, \er}\text{ intersects } M_{\tau_s}^{n, \er}\text{ only at }i\big\}.\]
Note that $i\in V_i$.
Let $C$ be the correspondence between 
$M_{\tau_s}^{n, \er}$ and $M_{\overline\lambda_s}^{n, \er}$ given by 
$C=\{(i,j)\ :\ i\in M_{\tau_s}^{n, \er},\ j\in V_i\}$.
Define a measure $\pi$ on $M_{\tau_s}^{n, \er}\times M_{\overline\lambda_s}^{n, \er}$ via
$\pi\{(i, j)\}=1/|\cC_{1}^{\sss n,\er}(\overline\lambda_s)|$ for $(i, j)\in C$.
Then on the event $F_{s,\eps}$,
\begin{align}\label{eqn:17}
\dis(C)\leq \frac{1}{2}d_{\rH}\big(M_{\underline\lambda_s}^{n, \er}, M_{\overline\lambda_s}^{n, \er}\big)\, ,\ \text{ and }\
\pi(C^c)=0.
\end{align}
Further, writing $\mu_1$ and $\mu_2$ for the uniform probability measures on $M_{\tau_s}^{n, \er}$ and $M_{\overline\lambda_s}^{n, \er}$ respectively, on the event $F_{s,\eps}$,
\begin{align}
D(\pi; \mu_1,\mu_2)
&\leq\sum_{i\in\cC_{1}^{n,\er}(\tau_s)}\bigg|\frac{1}{\big|\cC_{1}^{n,\er}(\tau_s)\big|}
-\frac{|V_i|}{\big|\cC_{1}^{n,\er}(\overline\lambda_s)\big|}\bigg|\nonumber\\
&\leq
\sum_{i\in\cC_{1}^{n,\er}(\tau_s)}
\left(\frac{1}{|\cC_{1}^{n,\er}(\tau_s)|}
-\frac{1}{\big|\cC_{1}^{n,\er}(\overline\lambda_s)\big|}\right)
+\sum_{i\in\cC_{1}^{n,\er}(\tau_s)}
\frac{|V_i|-1}{\big|\cC_{1}^{n,\er}(\overline\lambda_s)\big|}\nonumber\\
&
\leq
\sum_{i\in\cC_{1}^{n,\er}(\tau_s)}
\left(\frac{\big|\cC_{1}^{n,\er}(\overline\lambda_s)\big|
-\big|\cC_{1}^{n,\er}(\underline\lambda_s)\big|}{|\cC_{1}^{n,\er}(\underline\lambda_s)|\cdot \big|\cC_{1}^{n,\er}(\overline\lambda_s)\big|}\right)
+
\frac{\big|\cC_{1}^{n,\er}(\overline\lambda_s)\big|
-\big|\cC_{1}^{n,\er}(\underline\lambda_s)\big|}{\big|\cC_{1}^{n,\er}(\overline\lambda_s)\big|}
\end{align}
\begin{align}
&\hskip50pt\leq
2\times\frac{\big|\cC_{1}^{n,\er}(\overline\lambda_s)\big|
-\big|\cC_{1}^{n,\er}(\underline\lambda_s)\big|}{\big|\cC_{1}^{n,\er}(\underline\lambda_s)\big|}
\leq
2\frac{(1+\eps)^2-(1-\eps)^2}{(1-\eps)^2}\leq 32\eps
\, ,
\label{eqn:18}
\end{align}
where the last step uses $\eps<1/2$.
By \eqref{eqn:15}, \eqref{eqn:16}, and \eqref{eqn:18}, it follows that
\begin{align}\label{eqn:19}
\limsup_{n\to\infty}\ \pr\Big(d_{\GHP}\big(n^{-\frac{1}{3}}\cdot M_{\tau_s}^{n, \er}, \ \
n^{-\frac{1}{3}}\cdot M_{\overline\lambda_s}^{n, \er}\Big)>32\eps\Big)\leq 3\eps.
\end{align}
Note that
\begin{align}\label{eqn:20}
\Big\{2\underline\lambda_s(1-\eps)n^{2/3}
\leq
\big|\cC_{1}^{\sss n,\er}(\tau_s)\big|
\leq
2\overline\lambda_s(1+\eps)n^{2/3}
\Big\}
\supseteq
F_{s,\eps}\, .
\end{align}
Hence, on the event $F_{s,\eps}$,
\begin{align}\label{eqn:21}
d_{\GHP}\left(\frac{(12s)^{1/6}}{\big(|\cC_{1}^{\sss n,\er}(\tau_s)\big|\big)^{1/2}}M_{\tau_s}^{n, \er}, \
\frac{1}{n^{1/3}}M_{\tau_s}^{n, \er}\right)
&\leq
\diam\big(M_{\tau_s}^{n, \er}\big)\left|
\frac{(12s)^{1/6}}{\big(|\cC_{1}^{\sss n,\er}(\tau_s)\big|\big)^{1/2}}
-\frac{1}{n^{1/3}}\right|\nonumber\\
&\leq\frac{1}{n^{1/3}}\times \diam\big(M_{\infty}^{n, \er}\big)\times 3\eps\, ,
\end{align}
where the last step uses \eqref{eqn:20} and the relation $2\lambda_s^3=3s$.
By Theorem~\ref{thm:mst_complete}, the sequence of random variables $\big(n^{-1/3}\diam\big(M_{\infty}^{n, \er}\big)\, ;\ n\geq 1\big)$ is tight.
It thus follows from \eqref{eqn:16}, \eqref{eqn:19}, and \eqref{eqn:21} that
\begin{align}\label{eqn:22}
\limsup_{n\to\infty}\ \pr\bigg(d_{\GHP}
\left(\frac{(12s)^{1/6}}{\big(|\cC_{1}^{\sss n,\er}(\tau_s)\big|\big)^{1/2}}M_{\tau_s}^{n, \er}, \
\frac{1}{n^{1/3}}M_{\overline\lambda_s}^{n, \er}\right)
\geq\sqrt{\eps}\bigg)
=:\delta_{\eps}
\end{align}
satisfies
$\delta_{\eps}\downarrow 0$ as $\eps\downarrow 0$.

By \eqref{eqn:700},
$n^{-1/3}M_{\overline\lambda_s}^{n, \er}\weakc \cb^{\infty}\big(S_1(\overline\lambda_s)\big)$ as $n\to\infty$ w.r.t. GHP topology.
Combining this with \eqref{eqn:22} and \eqref{eqn:23}, we see that
\begin{align}\label{eqn:24}
d_{\mathrm{PR}}
\left(
\cL\bigg(\big(12s\big)^{1/6}\cdot\cb^{\infty}\big(\widetilde\cH^{(s)}\big)\bigg),\
\cL\bigg(\cb^{\infty}\big(S_1(\overline\lambda_s)\big)\bigg)
\right)
\leq 2\eps+\delta_{\eps}+\sqrt{\eps}.
\end{align}
Finally, by Theorem~\ref{thm:M-alternate},  $\cb^{\infty}\big(S_1(\overline\lambda_s)\big)\weakc\cM$ as \ch{$s\to\infty$} w.r.t. GHP topology.
Combining this observation with \eqref{eqn:24} implies that
\[\big(12s\big)^{1/6}\cdot\cb^{\infty}\big(\widetilde\cH^{(s)}\big)\weakc\cM\ \ \text{ as }\ \ s\to\infty\]
w.r.t. GHP topology.
Now the proof is completed by using Lemma~\ref{lem:tilde-equivalent} stated below.
\qed

\medskip

\begin{lem}\label{lem:tilde-equivalent}
For any bounded measurable $f:\fS_{\GHP}\to\bR$,
\[\bE\big[f\big(\big(12s\big)^{1/6}\cdot\cb^{\infty}\big(\widetilde\cH^{(s)}\big)\big)\big]
-
\bE\big[f\big(\big(12s\big)^{1/6}\cdot\cb^{\infty}\big(\cH^{(s)}\big)\big)\big]\to 0\ \ \text{ as }\ \ s\to\infty\, .\]
\end{lem}

\noindent{\bf Proof:}
Let $r=3(s-1)$.
Let $(X_1,\ldots, X_r)$ be as in Construction~\ref{constr:S-k-alternate} and $Y_i, Z_i$, $1\leq i\leq r$, and $\Gamma_{r/2}$ be as in the proof of Proposition~\ref{prop:gh-lemma-1}.
Then $\Gamma_{r/2}\sim\mathrm{Gamma}(r/2, 1)$, and as observed in \eqref{eqn:41}, $\sqrt{2}Z_i$, $1\leq i\leq r$, are i.i.d. $\mathrm{Exponential}(1)$ random variables.
Hence
\begin{align}\label{eqn:68}
L(\cH^{(s)})=\sum_{i=1}^r Y_i\sqrt{X_i}=\sum_{i=1}^r\frac{Z_i}{\sqrt{\Gamma_{r/2}}}=\sqrt{r}\cdot\big(1+o_P(1)\big).
\end{align}
Further, for any $s\geq 3$,
\begin{align}\label{eqn:69}
\bE\Big[\frac{L(\cH^{(s)})^2}{r}\Big]
\leq \frac{1}{r}\cdot\Big(\bE\big[\big(\sum_{i=1}^r Z_i\big)^4\big]\Big)^{1/2}\cdot\Big(\bE\big[\Gamma_{r/2}^{-2}\big]\Big)^{1/2}\leq C
\end{align}
for a universal constant $C$.
It follows from \eqref{eqn:68} and \eqref{eqn:69} that
\begin{align}\label{eqn:77}
\lim_{s\to\infty}\frac{1}{\sqrt{r}}\cdot\bE\big[L(\cH^{(s)})\big]= 1,
\end{align}
which in turn implies that
$r^{-1/2}\Big(L(\cH^{(s)})-\bE\big[L(\cH^{(s)})\big]\Big)\probc 0$
as $s\to\infty$.
Now
\[\bE\Big[\big(L(\cH^{(s)})-\bE\big[L(\cH^{(s)})\big]\big)^2\Big]\leq\bE\big[L(\cH^{(s)})^2\big]\leq Cr\]
by using \eqref{eqn:69}.
It thus follows that
\begin{align}\label{eqn:78}
\lim_{s\to\infty}r^{-1/2}\bE\big|L(\cH^{(s)})-\bE\big(L(\cH^{(s)})\big)\big|=0.
\end{align}
Hence
\begin{align*}
&\Big|\bE\Big(f\big(\big(12s\big)^{1/6}\cdot\cb^{\infty}\big(\widetilde\cH^{(s)}\big)\big)\Big)
-
\bE\Big(f\big(\big(12s\big)^{1/6}\cdot\cb^{\infty}\big(\cH^{(s)}\big)\big)\Big)\Big|
\\
&\hskip100pt
\leq\|f\|_{\infty}\cdot\frac{\bE\big|L(\cH^{(s)})-\bE\big(L(\cH^{(s)})\big)\big|}{\bE\big[L(\cH^{(s)})\big]}\to 0
\end{align*}
as $s\to\infty$, where the last step follows from \eqref{eqn:77} and \eqref{eqn:78}.
\qed

\section{Discussion}\label{sec:disc}

Here we briefly discuss universality of the scaling limit of the MST and related open problems.

\noindent{\bf (a) Universality of MST scaling limit for models exhibiting mean-field behavior:}
The geometry of the MST of an underlying discrete structure is closely related to the geometry of the structure under critical percolation.
The behavior under critical percolation of several models exhibiting mean-field behavior is well-understood.
In \cite{BBG-12}, the metric space scaling limit of the critical \erdos random graph was established.
Soon after this work, an abstract universality principle was developed in \cite{SBSSXW14,bhamidi-broutin-sen-wang} which was used to establish \erdos type scaling limits for a wide array of critical random graph models including the configuration model under critical percolation, various models of inhomogeneous random graphs, and the Bohman-Frieze process.
In \cite{bhamidi-sen}, the metric space scaling limit of random graphs with critical degree sequence having finite third moment was established.
Further, existing literature suggests that the components of the high-dimensional discrete torus \cite{hofstad-sapozhnikov,heydenreich-hofstad-1,heydenreich-hofstad-2} and the hypercube \cite{hofstad-nachmias} under critical percolation, and the critical quantum random graph model \cite{dembo-levit-vadlamani} also share the \erdos scaling limit.
It is believed that the scaling limit of the MST of each of these models exists and has the same law as $\cM$ up to a scaling factor.

\ch{We briefly discuss here how such a result might be proved for general random graphs with given degree sequences. 
	\begin{ass}\label{ass:cm-deg-general}
		Suppose $\vd=\vd^{\sss(n)}=(d_v^{\sss(n)},\ v\in [n])$ is a degree sequence for each $n\geq 1$, and write $\nu^n := n^{-1}\sum_{v\in[n]} \delta_{d^n_v}$ for the empirical degree distribution.
		Assume the following hold:
		\begin{enumeratei}
			\item\label{item:1}
			There exists a measure $\nu$ on $\bZ_{\geq 0}$ such that $\nu^n \to \nu$  as $n\to\infty$ w.r.t.\ the $W_3$ distance.
			\item\label{item:2}
			We have $\nu^n(0)=0$ for all $n\geq 1$, $\lim_{n\to\infty}n^{1/2}\cdot\nu^n(1)=0$, and $\nu(2)=0$.
		\end{enumeratei}
	\end{ass}
}
\ch{Assumption~\ref{ass:cm-deg-general}~\eqref{item:2} ensures that $\pr\big(\fB_n\big)\to 1$, where $\fB_n$ denotes the event that $\cmnd$ is connected \cite[Theorem 2.2]{federico-hofstad}.
Note that Assumption~\ref{ass:cm-deg-general}~\eqref{item:2} in particular implies that $\sigma_2(\nu)\geq 3\sigma_1(\nu)$, which is stronger than the condition for supercriticality of $\cmnd$, namely, $\sigma_2(\nu)>2\sigma_1(\nu)$.}
\ch{
	\begin{conj}\label{conj:1}
		Let $M^{\vd}$ (resp. $\overline{M}^{\vd}$) denote the MST of $\cmnd$ (resp. $\cgnd$), and let $\cM$ be as in Theorem~\ref{thm:mst_complete}.
		Then under Assumption~\ref{ass:cm-deg-general},
		\begin{align}\label{eqn:47}
			n^{-1/3}\cdot M^{\vd}\weakc \beta(\nu)\cdot\cM\ \ \text{ as }\ \ n\to\infty
		\end{align}
		with respect to the GHP topology, where the constant $\beta(\nu)$ is given by the following prescription:
		Let $D$ and $Y$ be random variables such that $D$ has law $\nu$ and, conditional on $D$, $Y$ is $\mathrm{Binomial}\big(D, p\big)$-distributed, 
		where $p=\sigma_1(\nu)/\big(\sigma_2(\nu)-\sigma_1(\nu)\big)$.
		Then
		\[
		\beta(\nu):=
		\bE\big[Y\big]\cdot\Big(\bE\big[Y^3\big]-4\cdot\bE\big[Y\big]\Big)^{-2/3}\, .
		\]
		Further, \eqref{eqn:47} continues to hold if we replace $M^{\vd}$ by $\overline{M}^{\vd}$.
	\end{conj}
}\ch{
	In the context of the $3$-regular configuration model, $\nu=\delta_{\{3\}}$.
	A simple calculation shows that $\beta\big(\delta_{\{3\}}\big)=6^{1/3}$, which is exactly the constant in Theorem~\ref{thm:mst_CM}.
	Let $U_e$, $e\in E(\cmnd)$, be i.i.d. $\mathrm{Uniform}[0, 1]$ random variables conditional on $\cmnd$.
	Suppose $M^{\vd}$ is constructed using these edge weights.
	For $\lambda\geq 0$, let 
	\[
	\fp^n_{\lambda}:=
	\frac{\sigma_1(\nu^n)}{\sigma_2(\nu^n)-\sigma_1(\nu^n)}+\frac{\lambda}{n^{1/3}}\, .
	\]
	Let $\cmnd(\lambda)$ be the graph with vertex set $[n]$ and edge set 
	$\big\{e\in E(\cmnd)\, :\, U_e\leq \fp^n_{\lambda} \big\}$.
	Write $\cC_1^{\vd}(\lambda)$ for the largest connected component of $\cmnd(\lambda)$, and let $M^{\vd}_{\lambda}$ denote the MST of $\cC_1^{\vd}(\lambda)$ constructed using the edge weights $U_e$,  $e\in E\big(\cC_1^{\vd}(\lambda)\big)$.
	Using Observation~\ref{observation:percolation}, on the event $\fB_n$, $M^{\vd}_{\lambda}$ is the restriction of $M^{\vd}$ to $\cC_1^{\vd}(\lambda)$.
	On $\fB_n$, consider the forest obtained by removing from $M^{\vd}$ the vertices in $M^{\vd}_{\lambda}$ and all edges incident to the vertices in $M^{\vd}_{\lambda}$; 
	let $\fT^{\vd}_{i, \lambda}$, $i=1,\ldots, k^{\vd}_{\lambda}$, denote the trees in this forest, and set
	$X^{\vd}_{\lambda}:=\max_i |\fT^{\vd}_{i, \lambda}|$.
	On $\fB_n^c$, set $X^{\vd}_{\lambda}=0$.}

\ch{
	To prove Conjecture~\ref{conj:1}, it is enough to prove the following two estimates:
For all $\eps>0$,
		\begin{gather}\label{eqn:45}
			\lim_{\lambda\to\infty}\limsup_{n\to\infty}\, \pr\bigg(
			\fB_n\cap
			\bigg\{d_{\rH}\big(M^{\vd}_{\lambda},\, M^{\vd} \big)>\eps n^{1/3} \bigg\}
			\bigg)
			=0\, ,\ \ \text{ and}\\
            \label{eqn:46}
			\lim_{\lambda\to\infty}\limsup_{n\to\infty}\, 
			\pr\big( X^{\vd}_{\lambda}>\eps n\big)=0\, .
		\end{gather}
Note that on the event $\fB_n$, $M^{\vd}_{\lambda}$ is a subtree of $M^{\vd}$, so it makes sense to measure the Hausdorff distance between them; moreover, this distance is bounded from above by the maximum of the diameters of the trees $\fT^{\vd}_{i, \lambda}$, $i=1,\ldots, k^{\vd}_{\lambda}$.} 
\ch{In the setting of the complete graph, results analogous to \eqref{eqn:45} and \eqref{eqn:46} were established in \cite{addarioberry-broutin-reed} and \cite[Lemma 4.11]{AddBroGolMie13} respectively.
	For inhomogeneous random graphs in the heavy-tailed regime, the analogue of \eqref{eqn:45} is proved in \cite[Theorem 6.1]{bhamidi-sen-geometry}.
}

\ch{
	The bound in \eqref{eqn:45} together with Lemma~\ref{lem:cycle-breaking-gives-mst}, Theorems~\ref{thm:scaling-given-degree-sequence},~\ref{thm:belongs-to-A-r}, and~\ref{thm:mst-from-ghp-convergence} would imply that \eqref{eqn:47} holds with respect to the GH topology. 
	This could be strengthened to GHP convergence with the help of the following exchangeability result:
	For each $v\in V(\cC^{\vd}_1(\lambda))$, let $d^{\sss(n)}_{v,\lambda}$ denote the degree of $v$ in $\cC^{\vd}_1(\lambda)$, and define
	$d_{v,\lambda}^{\avail}:=d^{\sss(n)}_{v}-d^{\sss(n)}_{v,\lambda}$.
	On the event $\fB_n$, for every $v\in V(M^{\vd}_{\lambda})$, let $\fT^{\vd}_{i, \lambda}(v)$, $1\leq i\leq r_{v,\lambda}$, be the trees among $\fT^{\vd}_{i, \lambda}$, $1\leq i\leq k^{\vd}_{\lambda}$, that are attached to $v$ in $M^{\vd}$.
	For every $v\in V(M^{\vd}_{\lambda})$, append $(d_{v,\lambda}^{\avail}-r_{v,\lambda})$ many zeros to the sequence $\big(\big|\fT^{\vd}_{i, \lambda}(v)\big|,~1\leq i\leq r_{v,\lambda}\big)$ and let $\big(\alpha_{\vd,\lambda}^{(i)}(v),~1\leq i\leq d_{v,\lambda}^{\avail}\big)$ be a uniform permutation of the resulting sequence; use independent permutations for different $v\in V(M^{\vd}_{\lambda})$ that are also independent of all the other random variables being considered.
	On $\fB_n^c$, set $\alpha_{\vd, \lambda}^{\sss (i)}(v):=0$ for  $v\in V(\cC_1^{\vd}(\lambda))$ and $i=1,\ldots, d_{v, \lambda}^{\avail}$.
	Then conditional on $\cmnd(\lambda)$ and $M^{\vd}_{\lambda}$, the family
	\begin{align}\label{eqn:488}
		\big(\alpha_{\vd, \lambda}^{\sss (i)}(v)\, ;\, 1\leq i\leq d_{v, \lambda}^{\avail}, v\in \cC_1^{\vd}(\lambda)\big)\
		\text{ of random variables is exchangeable}\,.
	\end{align}
}

\ch{
	The proof of \eqref{eqn:488} is similar to those of Lemma~\ref{lem:exchangeable}~\eqref{item:exchangeable} and Lemma~\ref{lem:exchangeable-erdos-renyi}.
	We outline the argument here for the readers' convenience. 
	We can generate $\cmnd(\lambda), M^{\vd}_{\lambda}$, and 
	$
	\big(\alpha_{\vd, \lambda}^{\sss (i)}(v)\, ;\, 1\leq i\leq d_{v, \lambda}^{\avail}, v\in \cC_1^{\vd}(\lambda)\big)
	$
	jointly as follows:
	\begin{enumeratea}
		\item
		Sample a $\mathrm{Binomial}\big(\sum_v d_{v}^{\sss (n)}/2\, ,\ 1-\fp_{\lambda}^n\big)$ random variable. For simplicity, we denote the realization by $m$.
		\item
		Consider the vertex set $[n]$ with $d_v^{\sss (n)}$ many half-edges attached to the vertex $v$.
		Conditional on step (a), sample $\cQ_{n,\vd,m}^{\sss (1)}$ as in Lemma~\ref{lem:cm-coupling}.
		By \eqref{eqn:333},
		\begin{align}\label{eqn:223}
			\cQ_{n,\vd,m}^{\sss (1)}
			\equald
			\cmnd(\lambda)~.
		\end{align}
		Thus, the largest component of $\cQ_{n,\vd,m}^{\sss (1)}$, say $\cC_1^{\bullet}$, has the same law as $\cC^{\vd}_1(\lambda)$.
		Let $d_v^{\bullet}$ be the degree of $v\in\cC_1^{\bullet}$.
		Then each $v\in\cC_1^{\bullet}$ has $d_v^{\sss(n)}-d_v^{\bullet}$ many `available' half-edges; 
		we denote them by $f_{v, i}$, $1\leq i\leq d_v^{\sss(n)}-d_v^{\bullet},\, v\in\cC_1^{\bullet}$.
		\item
		Assign i.i.d. $\text{Uniform}[0, \fp_{\lambda}^n]$ weights to the edges of $\cQ_{n,\vd,m}^{\sss (1)}$. 
		Let $M_1^{\bullet}$ be the MST of $\cC_1^{\bullet}$ constructed using these edge weights.
		Then
		\begin{align}\label{eqn:224}
			M_1^{\bullet}\equald M^{\vd}_{\lambda}
		\end{align}
		jointly with the equality in distribution in \eqref{eqn:223}.
		\item
		Conditional on steps (a), (b), and (c), sample $\cQ_{n,\vd,m}^{\sss (2)}$ as in Lemma~\ref{lem:cm-coupling}, i.e., by uniformly pairing the previously unpaired half-edges attached to the vertices in $\cQ_{n,\vd,m}^{\sss (1)}$.
		Then 
		$\cQ_{n,\vd,m}^{\sss (2)}
		\equald
		\cmnd$
		jointly with \eqref{eqn:223} and \eqref{eqn:224}.
		Let $\cE_{\text{new}}:=E(\cQ_{n,\vd,m}^{\sss (2)})\setminus E(\cQ_{n,\vd,m}^{\sss (1)})$.
		\item
		If $\cQ_{n,\vd,m}^{\sss (2)}$ is not connected, go to the next step. If $\cQ_{n,\vd,m}^{\sss (2)}$ is connected, assign i.i.d. $\text{Uniform}[\fp_{\lambda}^n,\, 1]$ weights to the edges in $\cE_{\text{new}}$.
		Let $w_e$ denote the weight assigned to the edge $e$.
		Delete all edges $e\in\cE_{\text{new}}$ \chl{for which there exists a cycle $\pi$ (with no repeated edges) in $\cQ_{n,\vd,m}^{\sss (2)}$ such that $e$ is an edge in $\pi$ and $w_e$ is the maximum among all edge weights in $\pi$.}
		This will yield a connected multigraph which we denote by $\cQ$.
		\item
		If $\cQ_{n,\vd,m}^{\sss (2)}$ is not connected, define $Q^{(i)}(v)$ to be the empty graph for $1\leq i\leq d_v^{\sss(n)}-d_v^{\bullet}$, $v\in \cC_1^{\bullet}$.
		If $\cQ_{n,\vd,m}^{\sss (2)}$ is connected, then note that $\cQ$ as constructed in (e) is simply $\cC_1^{\bullet}$ together with some connected multigraphs each of which is connected to a vertex of $\cC_1^{\bullet}$ by a single edge;
		for $v\in\cC_1^{\bullet}$ and $1\leq i\leq d_v^{\sss(n)}-d_v^{\bullet}$, set $Q^{(i)}(v)$ to be the connected multigraph that is connected to $v$ via the edge formed by pairing $f_{v, i}$ with another half edge in step (d), with the convention that $Q^{(i)}(v)$ is the empty graph if this edge was removed in step (e).
		Then 
		\[
		\bigg(\text{sort}\big(|Q^{(i)}(v)|~,\ 1\leq i\leq d_v^{\sss(n)}-d_v^{\bullet}\big)~;\ v\in \cC_1^{\bullet}\bigg)
		\equald
		\bigg(\text{sort}\big(\alpha_{\vd,\lambda}^{(i)}(v)~,\ 1\leq i\leq d_{v,\lambda}^{\avail}\big)~;\ v\in \cC^{\vd}_1(\lambda)\bigg)
		\]
		jointly with the distributional equalities in \eqref{eqn:223} and \eqref{eqn:224}, where $\text{sort}(\cdot)$ arranges the entries of a finite sequence in decreasing order.
	\end{enumeratea}
	Conditional on steps (a), (b), and (c) above, the rest of the procedure is symmetric with respect to the available half-edges attached to the vertices of $\cC_1^{\bullet}$.
	Hence, conditional on $\cQ_{n,\vd,m}^{\sss (1)}$ and $M_1^{\bullet}$, 
	the family
	$\big(|Q^{(i)}(v)|~;\ 1\leq i\leq d_v^{\sss(n)}-d_v^{\bullet},\, v\in \cC_1^{\bullet}\big)$
	is exchangeable.
	Thus, \eqref{eqn:488} follows.
}

\ch{
	The GH convergence in \eqref{eqn:47} can be lifted to GHP convergence using \eqref{eqn:46}, \eqref{eqn:488}, and arguments similar to the ones used in this paper.
	To carry out this argument, one would need to consider the metric measure space obtained by assigning mass $d_{v,\lambda}^{\avail}$ to each $v\in V(\cC^{\vd}_1(\lambda))$ and normalizing it to make it a probability measure.
	This does not quite fit into the framework of Theorem~\ref{thm:scaling-given-degree-sequence}~(i) where the same function $f$ is used for all vertices.
	This was done to keep the statement of that theorem simple.
	However, the proof of Theorem~\ref{thm:scaling-given-degree-sequence} outlined in Section~\ref{sec:appendix-given-degree-sequence} goes through without any change for the measure being considered in this setting.
}

\vskip3pt

\noindent{\bf (b) MST scaling limit in the heavy-tailed regime:}
This regime seems more interesting.
Consider scale-free random graphs on $n$ vertices where the tail of the empirical degree distribution $\nu_n$ asymptotically decays like $\nu_n([x,\infty))\sim x^{1-\tau}$ for some $\tau\in (3,4)$.
(In particular, the degree distribution asymptotically has infinite third moment and finite second moment.)
It is predicted \cite{braunstein2003optimal, braunstein2007optimal} that typical distance on the MST of such graphs scale like $n^{\frac{\tau-3}{\tau-1}}$.
In this regime, the scaling limit at criticality was first established in \cite{bhamidi-hofstad-sen} for inhomogeneous random graphs, and in \cite{SB-SD-vdH-SS} for random graphs with given degree sequences.
The recent preprint \cite{broutin-duquesne-wang} studies scaling limits of critical inhomogeneous random graphs in greater generality. 
The works \cite{CG17, GHS17} study scaling limits of critical random graphs with i.i.d. heavy-tailed degree sequences and alternate constructions of the limiting spaces.
Very recently in \cite{bhamidi-sen-geometry}, the scaling limit of the MST on the giant component in a supercritical inhomogeneous random graph with tail expoenent $\tau\in (3, 4)$ has been established.
Almost surely, the limiting space in \cite{bhamidi-sen-geometry} is compact, 
every point in this space either has degree one (leaf), or two, or infinity (hub), both
the set of leaves and the set of hubs are dense in this space, and the Minkowski dimension
of this space equals $(\tau-1)/(\tau-3)$.
We expect this space to be the candidate for the scaling limit of the MST of a wide array of heavy-tailed random graphs under some general assumptions.

\appendix
\section{} \label{sec:appendix}
Our aim in this section is to briefly describe the ideas needed to prove Theorems~\ref{thm:scaling-connected} and~\ref{thm:scaling-given-degree-sequence}.

\subsection{Sketch of proof of Theorem~\ref{thm:scaling-connected}}\label{sec:appendix-uniform-connected}
Suppose $\vt$ is a rooted tree with vertices labeled by $[m]$ and let $R(\cdot,\vt)$ be as in \eqref{eqn:72}.
For $s\geq 1$, let
\begin{align}\label{eqn:75}
A_s(\vt):=\big\{\big(v_1, u_1,\ldots, v_s,u_s\big) \ :
&\ 1\leq v_1\leq\ldots\leq v_s\leq m,\
u_i\in R(v_i,\vt),
\\
&\
\text{ if } i< j \text{ and }v_i=v_j\text{ then }u_i<u_j\big\}.
\end{align}
Note that $s!\times|A_s(\vt)|\leq |A_1(\vt)|^s$.
Let $T_m$ denote a uniform rooted labeled tree on $[m]$, and let $\ltms$ be distributed as
\begin{align}\label{eqn:76}
\bP\big(\ltms=\vt\big)=\frac{\pr\big(T_m=\vt\big)\cdot |A_s(\vt)|}{\bE\big(|A_s(T_m)|\big)}.
\end{align}
Then we have the following decomposition of $\cH_{m, s}$:

\begin{thm}\label{thm:uniform-connected-generate}
Fix $s\geq 1$. Sample $\ltms$, and
conditional on the realization, sample $\big(\lv_{1,m}, \lu_{1,m},\ldots, \lv_{s,m},\lu_{s,m}\big)$ from $A_s(\ltms)$ uniformly.
Place an edge between $\lv_{i,m}$ and $\lu_{i,m}$ for $1\leq i\leq s$, and then forget about the root of $\ltms$.
Call the resulting graph $\lhms$. 
Then $\lhms\equald\cH_{m, s}$.
\end{thm}

This can be seen as follows: 
Consider a simple, connected, rooted graph $G$ on $[m]$ with $\mathrm{sp}(G)=s$.
Let $\vt$ be the tree constructed by following a depth-first exploration of $G$ starting at its root, and let $v_i, u_i$, $1\leq i\leq s$, be the endpoints of the $s$ edges that need to be added to $\vt$ to recover $G$.
We can arrange $v_1, u_1,\ldots, v_s, u_s$ in a unique way so that the resulting sequence becomes an element of $A_s(\vt)$.
It thus follows that the set of simple, connected, rooted graphs on $[m]$ having $s$ surplus edges is in bijective correspondence with the set 
\begin{align}\label{eqn:90}
\big\{
(\vt, v_1, u_1,\ldots, v_s, u_s)\ :\ \vt\text{ rooted tree on }[m], (v_1, u_1,\ldots, v_s, u_s)\in A_s(\vt)
\big\} .
\end{align}
Then one can show that if we root $\cH_{m, s}$ at a uniform vertex, then its corresponding element in the set \eqref{eqn:90} will be distributed as
$\big(\ltms, \lv_{1,m},\lu_{1,m},\ldots,\lv_{s,m},\lu_{s,m}\big)$.
We omit the details as similar ideas have already been used in \cite{BBG-12, SBSSXW14, bhamidi-sen}.

For any tree $\vt$ on $[m]$ rooted at $\rho$, endow the children of each vertex in $\vt$ with the linear order induced by their labels.
Let $\rho=w_0, w_1,\ldots, w_{m-1}$ be the vertices of $\vt$ in order of appearance in a depth-first exploration of $\vt$ using the above order.
Let $\Ht_{\vt}:[0,m]\to\bR$ be the height function of $\vt$ given by $\Ht_{\vt}(m)=0$, and
\[
\Ht_{\vt}(x)=\height(w_{\lfloor x\rfloor},\vt),\ \ \ x\in[0,m).
\]
The following lemma is a collection of some standard results about $T_m$:
\begin{lem}\label{lem:uniform-connected}
\begin{inparaenumi}
\noindent\item\label{item:convergence}
The following convergences hold:
\begin{gather}
m^{-1/2}\Ht_{T_m}\big(m\cdot\big)\weakc 2\ve(\cdot),
\ \ \ \ \text{ and}\label{eqn:2323}\\
m^{-1/2}\max_{v\in [m]}\big|2\big|R(v, T_m)\big|-\height(v, T_m)\big|\probc 0,\label{eqn:2324}
\end{gather}
where the convergence in \eqref{eqn:2323} is w.r.t. the Skorohod $J_1$ topology.

\medskip

\noindent\item\label{item:height-bound}
For all $m\geq 1$, $\pr\big(\height(T_m)\geq x\sqrt{m}\big)\leq c x^3\exp\big(-x^2/2\big)$.

\medskip

\noindent\item\label{item:exploration-process-bound}
For all $x\geq 0$ and $m\geq 1$, 
\[\pr\big(\max_{v\in [m]}|R(v, T_m)|\geq x\sqrt{m}\big)\leq c\exp(-c'x^2).\]
Using the bounds $|A_s(T_m)|\times s!\leq |A_1(T_m)|^s$ and $|A_1(T_m)|\leq m\cdot\max_{v\in [m]}|R(v, T_m)|$, we further have
\[\pr\big(|A_s(T_m)|\geq x m^{3s/2}\big)\leq c\exp\big(-c'x^{2/s}\big)\]
for any $s\geq 1$, $x\geq 0$, and $m\geq 1$.

\medskip

\noindent\item\label{item:difference-negligible}
For any $s\geq 1$,
$m^{-3s/2}\big(|A_1(T_m)|^s-|A_s(T_m)|\times s!\big)\probc 0$.

\end{inparaenumi}
\end{lem}
Lemma~\ref{lem:uniform-connected}\eqref{item:convergence} follows from \cite{marckert-mokkadem}.
\eqref{item:height-bound} follows from \cite[Corollary 1]{luczak1995number}.
\eqref{item:exploration-process-bound} is the content of \cite[Lemma 13]{BBG-12}.
The proof of \eqref{item:difference-negligible} is similar to that of \cite[Lemma 7.3~(iii)]{bhamidi-sen}.

\medskip

\noindent{\bf Sketch of proof of \eqref{eqn:74}:}
In view of Theorem~\ref{thm:uniform-connected-generate}, $s\cdot\big(\hght(\ltms)+1\big)$ dominates $L(\cH_{m, s})$ stochastically for any $s\geq 1$.
Thus, \eqref{eqn:74} follows from Lemma~\ref{lem:uniform-connected} \eqref{item:height-bound} and \eqref{item:exploration-process-bound}.
\qed

\medskip

To prove the other assertions in Theorem~\ref{thm:scaling-connected} it will be convenient to work with two slightly different spaces $\cH_{m, s}^{\circ}$ and $\cH_{m, s}^{\dagger}$ which we define next.
Recall the notation $R(\cdot,\cdot,\cdot)$ from \eqref{eqn:722}. Sample $T_m^{\circ}$ according to distribution 
\begin{align}\label{eqn:927}
\pr\big(T_m^{\circ}=\vt\big)=\frac{\pr(T_m=\vt)\cdot|A_1(\vt)|^s}{\bE\big[|A_1(T_m)|^s\big]}~,\ \ \ \
\vt\text{ rooted tree on }[m].
\end{align}
Conditional on $T_m^{\circ}$, sample an i.i.d. sequence of triples $(v_{i,m}^{\circ}, u_{i,m}^{\circ}, f_{i,m}^{\circ})$, $1\leq i\leq s$, where
\begin{gather*}
\pr\big(v_{i,m}^{\circ}=v\mid T_m^{\circ}\big)=|R(v, T_m^{\circ})|\big/|A_1(T_m^{\circ})|,\ \ \ 
v\in[m],\\
\pr\big(u_{i,m}^{\circ}=u\mid T_m^{\circ},v_{i,m}^{\circ}\big)
=|R(u, v_{i,m}^{\circ}, T_m^{\circ})|\big/|R(v, T_m^{\circ})|,\ \ \ 
u\in\big\{\stackrel{\longleftarrow{\sss(k)}}{(v_{i,m}^{\circ}) }\ :\ 1\leq k\leq\height(v_{i,m}^{\circ})\big\},\ \ \text{ and}\\
\pr\big(f_{i,m}^{\circ}=f\mid T_m^{\circ},v_{i,m}^{\circ}, u_{i,m}^{\circ}\big)
=1\big/|R(u_{i,m}^{\circ}, v_{i,m}^{\circ}, T_m^{\circ})|,\ \ \
f\in R(u_{i,m}^{\circ}, v_{i,m}^{\circ}, T_m^{\circ}). 
\end{gather*}
Let $\cH_{m, s}^{\dagger}$ (resp. $\cH_{m, s}^{\circ}$) be the space obtained by adding an edge between $v_{i,m}^{\circ}$ and $f_{i,m}^{\circ}$ (resp. between $v_{i,m}^{\circ}$ and $u_{i,m}^{\circ}$) for $1\leq i\leq s$, and then forgetting about the root of $T_m^{\circ}$. 
It follows from Lemma~\ref{lem:uniform-connected} \eqref{item:exploration-process-bound} and \eqref{item:difference-negligible} that the total variation distance between the laws of $\overline\cH_{m, s}$ (as defined in Theorem~\ref{thm:uniform-connected-generate}) and $\cH_{m, s}^{\dagger}$ tends to zero as $m\to\infty$.
It thus follows from Theorem~\ref{thm:uniform-connected-generate} that there exists a coupling of $\cH_{m, s}$ and $\cH_{m, s}^\dagger$ such that
\begin{align}\label{eqn:88}
\pr\big(\cH_{m, s}\neq \cH_{m, s}^\dagger\big)\to 0,\ \ \text{ as }\ \ m\to\infty.
\end{align}

We will now recall an alternate construction of $\cH^{(s)}$ which is essentially given in \cite{BBG-12}; see also the discussion below \cite[Equation (1)]{BBG-limit-prop-11}. We first introduce some notation. For any $f:[0,1]\to\bR, x\in[0,1]$, and $h>0$, let
\[
\prev(x,h;f)=\sup\big\{y\in[0,x) : f(y)=h\big\}, \ \ \ \text{ and } \ \ \ 
\nxt(x,h;f)=\inf\big\{y\in(x,1] : f(y)<h\big\},
\]
where $\sup\{\ \}=-\infty$ and $\inf\{\ \}=\infty$ by convention.

\begin{constr}[Alternate construction of $\cH^{(s)}$]\label{constr:S-k}
Fix an integer $s\geq 2$.
	\begin{enumeratea}
		\item Sample $\ve^{\circ}$ with law given by
		\[
		\bE\big[f(\ve^{\circ})\big]=\frac{\bE\big[f(\ve)\big(\int_0^1 \ve(t)dt\big)^s\big]}{\bE\big[\big(\int_0^1 \ve(t)dt\big)^s\big]}.
		\]
		\item Conditional on $\ve^{\circ}$, sample i.i.d. points $y_1^{\circ},\ldots, y_s^{\circ}$ having density $\ve^{\circ}(y)\big/\int_0^1 \ve^{\circ}(t)dt $.
		\item Conditional on the above, sample $h_1^{\circ},\ldots,h_s^{\circ}$ independently, where $h_i^{\circ}\sim\mathrm{Unif}[0,\ve^{\circ}(y_i^\circ)]$. 
		Set $x_i^\circ=\prev(y_i^\circ, h_i^\circ; \ve^\circ)$.
		\item Form the quotient space  $\cT_{\ve^\circ}/\sim$, where $\sim$ is the equivalence relation under which $q_{\ve^\circ}(x_i^\circ)\sim q_{\ve^\circ}(y_i^\circ)$, $1\leq i\leq s$.
	\end{enumeratea}
Then $\cH^{(s)}\equald 2\cdot\big(\cT_{\ve^\circ}/\sim\big)$.
\end{constr}

Now observe that $\cH_{m, s}^{\circ}$ has a similar alternate construction:  
First sample $T_m^\circ$ as in \eqref{eqn:927}.
Let $w_0,\ldots,w_{m-1}$ be the vertices of $T_m^\circ$ in order of appearance in a depth-first exploration of $T_m^\circ$.
Let $\Ht^\circ$ be the height function of $T_m^\circ$.
Conditional on $T_m^\circ$, sample i.i.d. random variables $y_{1,m}^\circ,\ldots,y_{s,m}^\circ$, where 
\[
\pr\big(y_{i,m}^\circ=j\mid T_m^\circ\big)=|R(w_j, T_m^\circ)|\big/|A_1(T_m^\circ)|,\ \ \ 
1\leq j\leq m-1.
\]
Conditional on the above, sample $h_{1,m}^\circ,\ldots, h_{s,m}^\circ$ independently via
\[
\pr\big(h_{i,m}^\circ=\Ht^\circ(y_{i,m}^\circ)-k\ \big|\ T_m^\circ,y_{1,m}^\circ,\ldots,y_{s,m}^\circ\big)
=
\frac{|R\big(\stackrel{\leftarrow\sss{(k)}}{v\ \ \ }, v, T_m^\circ \big)|}{|R(v, T_m^\circ)|},\ \ \
1\leq k\leq \Ht^\circ(y_{i,m}^\circ),
\]
where $v=w_{y_{i,m}^\circ}$.
Let $x_{i,m}^\circ=\prev(y_{i,m}^\circ, h_{i,m}^{\circ}; \Ht^\circ)-1$.
Then $\cH_{m, s}^\circ$ has the same distribution as the space obtained by placing an edge in $T_m^\circ$ between $w_{y_{i,m}^\circ}$ and $w_{x_{i,m}^\circ}$ for $1\leq i\leq s$.

\medskip

\noindent{\bf Sketch of proof of \eqref{eqn:73}:}
Using Lemma~\ref{lem:uniform-connected} \eqref{item:convergence} and \eqref{item:exploration-process-bound}, it can be shown that the following convergences hold jointly:
\begin{align}\label{eqn:666}
\frac{1}{\sqrt{m}}\Ht^\circ\big(m\cdot\big)\weakc 2\ve^\circ(\cdot), \ \ \text{ and }\ \  
\Big(\frac{x_{i,m}^\circ}{m}, \frac{y_{i,m}^\circ}{m}, \frac{h_{i,m}^\circ}{\sqrt{m}}\Big)
\weakc
\big(x_i^\circ, y_i^\circ, 2h_i^\circ\big),\ \ 1\leq i\leq s
\end{align}
as $m\to\infty$. 
Using Construction~\ref{constr:S-k} and the above alternate construction of $\cH_{m, s}^\circ$, it is now routine to prove the assertion in \eqref{eqn:73} for $\cH_{m, s}^\circ$, from which it follows that the same is true for $\cH_{m, s}^\dagger$.
The desired result now follows from \eqref{eqn:88}.
\qed

\medskip

Let $y_{(i),m}^\circ$ (resp. $y_{(i)}^\circ$), $1\leq i\leq s$, be $y_{i,m}^\circ$ (resp. $y_i^\circ$), $1\leq i\leq s$, arranged in an increasing order.
For $1\leq i\leq s-1$ define $z_{i,m}^\circ$ and $z_i^\circ$ via
\begin{gather*}
z_{i,m}^\circ=\min\Big\{t\in[y_{(i),m}^\circ, y_{(i+1),m}^\circ]\ :\ \Ht^\circ(t)=\min\big\{ \Ht^\circ(a) :  y_{(i),m}^\circ\leq a\leq y_{(i+1),m}^\circ \big\}  \Big\},\ \ \text{ and}\\
\ve^\circ(z_i^\circ)=\inf\big\{ \ve^\circ(t) : y_{(i)}^\circ\leq t\leq  y_{(i+1)}^\circ \big\}.
\end{gather*}
Further, define
\begin{gather*}
x_{i,m}^{\circ, +}=\nxt\big(x_{i,m}^\circ, h_{i,m}^\circ+1;\ \Ht^\circ\big),\ \ 
x_i^{\circ, +}=\nxt\big(x_i^\circ, h_i^\circ; \ve^\circ\big),\ \ 1\leq i\leq s,\\
z_{i,m}^{\circ, -}=\prev\big(z_{i,m}^\circ, \Ht^\circ(z_{i,m}^\circ)-1;\ \Ht^\circ\big)-1,\ \ 
z_i^{\circ, -}=\prev\big(z_i^\circ, \ve^\circ(z_i^\circ); \ve^\circ\big),\ \ 1\leq i\leq s-1,\\
z_{i,m}^{\circ, +}=\nxt\big(z_{i,m}^\circ, \Ht^\circ(z_{i,m}^\circ);\ \Ht^\circ\big)-1,\ \ 
z_i^{\circ, +}=\nxt\big(z_i^\circ, \ve^\circ(z_i^\circ); \ve^\circ\big),\ \ 1\leq i\leq s-1.
\end{gather*}

\medskip

\noindent{\bf Sketch of proof of \eqref{eqn:711}:}
From \eqref{eqn:666} it follows that the following convergence holds jointly with the convergence in $\eqref{eqn:666}$: As $m\to\infty$,
\begin{align}\label{eqn:998}
\frac{x_{i,m}^{\circ, +}}{m}\weakc x_i^{\circ, +},\ 1\leq i\leq s,\  \text{ and }\ 
\frac{1}{m}\big(z_{i,m}^\circ,\ z_{i,m}^{\circ, -},\  z_{i,m}^{\circ, +}\big)\weakc \big(z_i^\circ, z_i^{\circ, -}, z_i^{\circ, +}\big),\ 1\leq i\leq s-1.
\end{align}
Arrange $x_i^\circ, x_i^{\circ, +}, y_i^\circ$, $1\leq i\leq s$, and 
$z_i^\circ, z_i^{\circ, -}, z_i^{\circ, +}$, $1\leq i\leq s-1$, 
(resp. $x_{i,m}^\circ, x_{i,m}^{\circ, +}, y_{i,m}^\circ$, $1\leq i\leq s$, and 
$z_{i,m}^\circ, z_{i,m}^{\circ, -}, z_{i,m}^{\circ, +}$, $1\leq i\leq s-1$)
in increasing order as $a_1,\ldots,a_{6s-3}$ (resp. as $a_{1,m},\ldots,a_{6s-3, m}$).
Let 
\[
\Delta_j=a_{j+1}-a_j,\ \ \text{ and }\ \ 
\Delta_{j,m}=a_{j+1,m}-a_{j,m},\ \ 1\leq j\leq 6s-4.
\]
Then it follows from \eqref{eqn:998} and the second convergence in \eqref{eqn:666} that
\begin{align}\label{eqn:111}
\big(\Delta_{j,m},\ 1\leq j\leq 6s-4\big)
\weakc
\big(\Delta_j,\ 1\leq j\leq 6s-4\big),\ \ \text{ as }\ \ m\to\infty
\end{align}
jointly with \eqref{eqn:666} and \eqref{eqn:998}.

Recall the notation used in \eqref{eqn:711}, and note that there exists a partition $\cP=\{\cP_1,\ldots,\cP_r\}$ of $[6s-4]$ that depends only on the realizations of $\ve^\circ$ and $x_i^\circ,y_i^\circ$, $1\leq i\leq s$, such that
\begin{align}\label{eqn:999}
\big(\mu^{(s)}\big(\cT_i'\big), 1\leq i\leq r\big)
\equald
\big(\sum_{j\in\cP_i}\Delta_j,\ 1\leq i\leq r\big).
\end{align}
Further, it follows from \eqref{eqn:666} that for large $m$, the vector consisting of the numbers of vertices in $\cH_{m, s}^\circ$ that are connected to the different elements of $e(\cH_{m, s}^\circ)$ is given by $\big(\sum_{j\in\cP_i}\Delta_{j,m},\ 1\leq i\leq r\big)$, where the common endpoints of multiple $e\in e(\cH_{m, s}^\circ)$ and the vertices in their pendant subtrees have been accounted for in $\sum_{j\in\cP_i}\Delta_{j,m}$ for exactly one value of $i$ in a specific way.
Using \eqref{eqn:111} and \eqref{eqn:999}, we get the analogue of \eqref{eqn:711} for $\cH_{m, s}^\circ$ for the above specific way of assigning the common endpoints of multiple $e\in e(\cH_{m, s}^\circ)$ and the vertices in their pendant subtrees to the different terms $\sum_{j\in\cP_i}\Delta_{j,m}$.

This together with \eqref{eqn:88} would complete the proof if we can show that the sizes of the pendant subtrees of the common endpoints of multiple $e\in e(\cH_{m, s}^\circ)$ are asymptotically negligible. 
This negligibility claim follows from the following facts:
\begin{enumerateA}
	\item $Y_{i,m}=o_P(m)$, $1\leq i\leq s$, where $Y_{i,m}$ denotes the number of descendants of $v_{i,m}^\circ$ in $T_m^\circ$.
	\item $X_{i,m}=o_P(m)$, $1\leq i\leq s$, where $X_{i,m}$ denotes the number of descendants of $u_{i,m}^\circ$ in $T_m^\circ$ that are not in the subtree that contains $v_{i,m}^\circ$.
	\item For every $\eps>0$,
	\[
	\pr\big(\exists v\in T_m^\circ : v\text{ has at least three subtrees in }T_m^\circ\text{ each of size }\geq\eps m\big)
	\to 0\ \ \text{ as }\ \ m\to\infty.
	\]
\end{enumerateA}

(A) and (C) follow from \eqref{eqn:666} and the facts that $q_{\ve^\circ}(y_i^\circ)$ is almost surely a leaf in $\cT_{\ve^\circ}$ and that $\cT_{\ve^\circ}$ is almost surely binary. The proof of (B) is also routine.

\subsection{Sketch of proof of Theorem~\ref{thm:scaling-given-degree-sequence}}\label{sec:appendix-given-degree-sequence}
Assume that for each $m\geq 1$, $\mvk^{\sss(m)}=(k_i^{\sss(m)}, i\ge 0)$, where $k_i^{\sss(m)}$ are nonnegative integers satisfying $\sum_{i\geq 0}k_i^{\sss(m)}=m$ and $\sum_{i\geq 0}ik_i^{\sss(m)}=m-1$.
Then there exist trees on $m$ vertices in which for each $i\geq 0$, there are exactly $k_i^{\sss(m)}$ many vertices with $i$ many children.
We call $\mvk^{\sss(m)}$ the child sequence of such a tree.
Assumption~\ref{ass:cm-deg} gives the criterion for graphs with given degree sequences to be critical.
The following assumption gives the analogous criterion for plane trees with given child sequences.

\begin{ass}\label{ass:ecd}
There exists a pmf $(p_0, p_1,\ldots)$ with
\[ p_0>0,\quad  \sum_{i\ge 1}i p_i=1,\quad\text{and }\sum_{i\ge 1}i^2 p_i<\infty\]
such that
\[\frac{k_i^{\sss(m)}}{m}\to p_i\ \text{ for }\ i\ge 0,\ \text{ and }\
\frac{1}{m}\sum_{i\ge 0}i^2 k_i^{\sss(m)}\to\sum_{i\ge 1}i^2 p_i.\]
We will write $\sigma^2 = \sum_i i^2 p_i - 1$ for the variance associated with the pmf $(p_0, p_1,\ldots)$.
\end{ass}

Let $\bT_{\mvk^{\sss(m)}}$ be the set of plane trees with child sequence $\mvk^{\sss(m)}$.
Let $\cT_{\mvk^{\sss(m)}}$ be a uniform element of $\bT_{\mvk^{\sss(m)}}$ endowed with the tree distance and the uniform probability measure on $m$ vertices and viewed as a metric measure space.
Broutin and Marckert \cite{broutin-marckert} showed that under Assumption~\ref{ass:ecd}, $\sigma m^{-1/2}\cT_{\mvk^{\sss(m)}}\weakc\cT_{2\ve}$
w.r.t. GHP topology.
The following variant of this result follows from \cite[Lemma 7.4 and Lemma 7.6]{bhamidi-sen}:

\begin{lem}\label{lem:plane-tree-available-measure}
Suppose $\mvk^{\sss(m)}$ satisfies Assumption~\ref{ass:ecd}.
Further, suppose $f_m:\{0, 1,\ldots\}\to [0,1]$ is such that
\[\sum_{i\geq 0}k_i^{\sss(m)} f_m(i)=1,\ \ \text{ and }\ \ \lim_{m\to\infty}\ \max_{i: k_i^{\sss(m)}>0}\ f_m(i)=0.\]
Let $\cT_{\mvk^{\sss(m)}}^{f_m}$ be a uniform element of $\bT_{\mvk^{\sss(m)}}$ endowed with the tree distance and the measure that assigns probability $f_m(i)$ to any node that has $i$ children, $i\geq 0$.
Then
\[\sigma m^{-1/2}\cdot\cT_{\mvk^{\sss(m)}}^{f_m}\weakc\cT_{2\ve}\ \ \text{ w.r.t. GHP topology}.\]
\end{lem}

Now we can prove Theorem~\ref{thm:scaling-given-degree-sequence} using the above lemma and the techniques used in the proof of  \cite[Theorem 2.2]{bhamidi-sen}.

\section*{Acknowledgments}
The authors thank two anonymous referees whose careful reading and detailed comments led to significant improvements in the paper.
LAB was supported in part by an NSERC Discovery Grant and Discovery Accelerator Supplement, and by an FRQNT Team Grant, during the preparation of this research. 
SS was partially supported by the Infosys Foundation, Bangalore and by MATRICS grant MTR/2019/000745 from SERB.

\bibliographystyle{plain}
\bibliography{regular_bib}

\end{document}